\documentclass[11pt,twoside]{amsart}
\usepackage{amssymb,amsmath}

\newcommand{\Z}{{\mathbb{Z}}}
\newcommand{\N}{{\mathbb{N}}}
\newcommand{\Q}{{\mathbb{Q}}}
\newcommand{\R}{{\mathbb{R}}}
\newcommand{\fT}{{\mathbb{T}}}
\newcommand{\bDt}{{\mathbf{\Delta}}}
\newcommand{\bC}{{\mathbf{C}}}
\newcommand{\bD}{{\mathbf{D}}}
\newcommand{\bH}{{\mathcal{H}}}
\newcommand{\cA}{{\mathcal{A}}}

\newcommand{\cD}{{\mathcal{D}}}
\newcommand{\cB}{{\mathcal{B}}}

\newcommand{\cL}{{\mathcal{L}}}
\newcommand{\cR}{{\mathcal{R}}}
\newcommand{\cLR}{{\mathcal{LR}}}
\newcommand{\cM}{{\mathcal{E}}}

\newcommand{\fC}{{\mathfrak{C}}}
\newcommand{\fR}{{\mathfrak{R}}}
\newcommand{\fS}{{\mathfrak{S}}}
\newcommand{\fI}{{\mathfrak{I}}}
\newcommand{\fX}{{\mathfrak{X}}}
\newcommand{\ba}{{\mathbf{a}}}

\newcommand{\Irr}{{\operatorname{Irr}}}
\newcommand{\Ind}{{\operatorname{Ind}}}

\renewcommand{\leq}{\leqslant}
\renewcommand{\geq}{\geqslant}

\renewcommand{\atop}[2]{\genfrac{}{}{0pt}{}{#1}{#2}}

\newtheorem{thm}{Theorem}[section]
\newtheorem{lem}[thm]{Lemma}
\newtheorem{cor}[thm]{Corollary}
\newtheorem{prop}[thm]{Proposition}
\newtheorem{conj}[thm]{Conjecture}

\theoremstyle{definition}
\newtheorem{defn}[thm]{Definition}
\newtheorem{exmp}[thm]{Example}

\theoremstyle{remark}
\newtheorem{rem}[thm]{Remark}

\begin{document}

\title{Relative Kazhdan--Lusztig cells}

\author{Meinolf Geck}
\address{Institut Girard Desargues, bat. Jean Braconnier, Universit\'e Lyon 1,
21 av Claude Bernard, F--69622 Villeurbanne Cedex, France}
\email{geck@igd.univ-lyon1.fr}
\date{february 12, 2005}

\begin{abstract} In this paper, we study the Kazhdan--Lusztig cells of a 
Coxeter group $W$ in a ``relative'' setting, with respect to a parabolic 
subgroup $W_I \subseteq W$. This relies on a factorization of the
Kazhdan--Lusztig basis $\{\bC_w\}$ of the corresponding (multi-parameter) 
Iwahori--Hecke algebra with respect to $W_I$. We obtain two applications to 
the ``asymptotic case'' in type $B_n$, as introduced by Bonnaf\'e--Iancu: 
we show that $\{\bC_w\}$ is a ``cellular basis'' in the sense of 
Graham--Lehrer, and we construct an analogue of Lusztig's canonical 
isomorphism from the Iwahori--Hecke algebra to the group algebra of
the underlying Weyl group of type $B_n$.
\end{abstract}

\maketitle

\pagestyle{myheadings}

\markboth{Geck}{Relative Kazhdan--Lusztig cells}

\section{Introduction} \label{sec0}
Let $W$ be a Coxeter group and $L\colon W \rightarrow {\Z}_{\geq 0}$ be a 
weight function, in the sense of Lusztig \cite{Lusztig03}. This gives rise 
to various pre-order relations on $W$, usually denoted by $\leq_{\cL}$, 
$\leq_{\cR}$ and $\leq_{\cLR}$. Let $\sim_{\cL}$, $\sim_{\cR}$ and 
$\sim_{\cLR}$ be the corresponding equivalence relations. The equivalence 
classes are called the left, right and two-sided cells of $W$, 
respectively. They were first defined by Kazhdan and Lusztig \cite{KaLu} in 
the case where $L$ is the length function on $W$ (the ``equal parameter 
case''), and by Lusztig \cite{Lusztig83} in general. They play a fundamental 
role, for example, in the representation theory of finite or $p$-adic 
groups of Lie type; see Lusztig \cite{LuBook}, \cite{Lu1} and the survey 
in \cite[Chap.~0]{Lusztig03}. 

The definition of the above relations relies on the construction of the 
Kazhdan--Lusztig basis $\{\bC_w \mid w \in W\}$ in the associated 
Iwahori--Hecke algebra $\bH$. This paper arose from an attempt to show 
that the basis $\{\bC_w\}$ is a ``cellular basis'' in the sense of 
Graham--Lehrer \cite{GrLe}, in the case where $W=W_n$ is of type $B_n$ 
with diagram and weight function given by
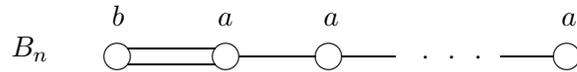
\begin{center}
\begin{picture}(250,32)
\put(  0, 10){$B_n$}
\put( 40, 10){\circle{10}}
\put( 44,  7){\line(1,0){33}}
\put( 44, 13){\line(1,0){33}}
\put( 81, 10){\circle{10}}
\put( 86, 10){\line(1,0){29}}
\put(120, 10){\circle{10}}
\put(125, 10){\line(1,0){20}}
\put(155,  7){$\cdot$}
\put(165,  7){$\cdot$}
\put(175,  7){$\cdot$}
\put(185, 10){\line(1,0){20}}
\put(210, 10){\circle{10}}
\put( 38, 22){$b$}
\put( 78, 22){$a$}
\put(118, 22){$a$}
\put(208, 22){$a$}
\end{picture}
\end{center}
where $a,b$ are positive integers such that $b/a$ is ``large'' with respect
to $n$. This is the ``asymptotic case'' studied by Bonnaf\'e--Iancu \cite{BI}.

After a number of intermediate results, this goal will be achieved in 
Section~\ref{sec-repB}. Those intermediate results concern 
properties of left, right and two-sided cells which are important in 
their own right. In fact, combining the results in this paper with the 
results of Bonnaf\'e--Iancu \cite{BI}, Bonnaf\'e \cite{BI2} and Geck--Iancu 
\cite{GeIa05}, we have that {\bf (P1)--(P14)} from Lusztig's list of 
conjectures in \cite[Chap.~14]{Lusztig03}, as well as a weak version 
of {\bf (P15)}, hold in the ``asymptotic case'' in type $B_n$. The weak 
version of {\bf (P15)} is sufficient, for example, to establish the 
existence of an analogue of Lusztig's canonical isomorphism \cite{Lu0}
for the two-parameter Iwahori--Hecke algebra of type $B_n$. (These things
will be discussed at the end of this paper, in Section~7.)

The main and unifying idea of this paper is to combine the existing theory
(due to Lusztig in general, and to Bonnaf\'e and Iancu as far as type
$B_n$ is concerned) with a detailed analysis of the decomposition of the 
Kazhdan--Lusztig basis of a Coxeter group with respect to a parabolic 
subgroup, based on the author's article \cite{myind}. 

Here is the first property that we consider in this paper. It has been
conjectured by Lusztig \cite[14.2]{Lusztig03} that we always have the 
following implication for elements $x,y$ in a Coxeter group $W$:
\begin{equation*} 
x\leq_{\cL} y \quad \mbox{and} \quad x\sim_{\cLR} y \quad\Rightarrow \quad x 
\sim_{\cL} y. \tag{$\spadesuit$}
\end{equation*}
This is known to hold in the equal parameter case when $W$ is a finite or 
affine Weyl group\footnote{Other situations where ($\spadesuit$) is known 
to hold include the quasi-split case discussed in \cite[Chap.~16]{Lusztig03}
(which is derived from the equal parameter case), and explicitly worked 
examples like the infinite dihedral group in \cite[Chap.~17]{Lusztig03} 
or type $F_4$ in \cite{my04}.}; see Lusztig \cite{Lu1}. However, although 
all the notions involved in the above statement are completely elementary, 
the proof is surprisingly complicated: it relies on a  geometric 
interpretation of the Kazhdan--Lusztig basis of $\bH$ and some deep results 
from algebraic geometry; see Springer \cite{Spr} and Lusztig \cite{Lu1}. A 
somewhat different proof is given by Lusztig \cite{Lu0} for finite Weyl 
groups (relying on the connection between cells and primitive ideals in
universal envelopping algebras via the main conjecture in \cite{KaLu}); 
in that article, ($\spadesuit$) is used to construct a canonical isomorphism 
from $\bH$ to the group algebra of $W$. The property ($\spadesuit$) also 
plays an important role in Lusztig's study \cite{LuBook} of representations 
of reductive groups over finite fields.

In Section~4, we develop the formulation of a relative version of 
($\spadesuit$), taking into account the presence of a parabolic subgroup 
$W_I \subseteq W$. (The original version of ($\spadesuit$) corresponds to 
the case where $W_I=W$.) The tools for dealing with this relative setting 
are provided by \cite{myind}; we recall the basic ingredients, with some 
refinements, in Section~3. We conjecture that the relative version of
($\spadesuit$) holds for all $W,L$ and all choices of $W_I \subseteq W$. 
In Section~4, we prove our conjecture for finite and for affine Weyl groups 
in the equal parameter case. The method is inspired by Lusztig's proof of 
($\spadesuit$) in \cite[Chap.~15]{Lusztig03}. The additional complication 
arising from the presence of $W_I$ is dealt with by Lemma~\ref{relcellA0}, 
which reduces to a triviality if $W_I=W$. 

A priori, we do not have any geometric interpretation of the 
Kazhdan--Lusztig basis in the general case of unequal parameters. (Note, 
however, that there is a conjectural geometrical interpretation by Lusztig 
\cite[Chap.~27]{Lusztig03} for certain values of the parameters.) So the
above methods and results will not apply in type $B_n$ with parameters
as specified as above. In Theorem~\ref{relsn}, we do prove ($\spadesuit$) 
in this case, by reduction to the relative version of ($\spadesuit$) for 
the symmetric group $\fS_n$. Thus, eventually, the proof of ($\spadesuit$) 
in type $B_n$ also rests on the geometric interpretation of the 
Kazhdan--Lusztig basis for $\fS_n$. The proof of that reduction argument 
occupies almost all of Section~\ref{sec-cellu}; this relies once more on 
the results in \cite{myind}, and on the results of Bonnaf\'e--Iancu \cite{BI}
and Bonnaf\'e \cite{BI2} on the left cells and two-sided cells, respectively. 
At one point in the proof, we also use an idea of Dipper--James--Murphy 
\cite{DiJa95} to deal with the action of the generator with parameter~$b$ 
in the above diagram. 

In Section~6, we go on to study the representations carried by the
left cells in the ``asymptotic case'' in type $B_n$. The main result, 
Theorem~\ref{eqcellsB}, shows that two left cells which afford the same
character actually give rise to exactly the same representation (and not 
only equivalent ones). Again, the proof relies on the techniques in
\cite{myind}, concerning the ``induction'' of cells. An analogous result 
for the left cell representations of the symmetric group has already been 
obtained by Kazhdan--Lusztig in their original article \cite{KaLu} where 
they introduced left cells and the corresponding representations. 

Combining the main results of Bonnaf\'e--Iancu \cite{BI} and Bonnaf\'e
\cite{BI2} with Theorem~\ref{relsn} and Theorem~\ref{eqcellsB} in 
this paper, we immediately get that $\{\bC_w\}$ is a ``cellular basis''
in the ``asymptotic case'' in type $B_n$; see Corollary~\ref{grle}. 

As a further application of our results, we can exhibit a new basis in 
the Iwahori--Hecke algebra of type $B_n$ whose structure constants are 
integers. (In fact, the structure constants are $0$, $1$.) This uses
an idea of Neunh\"offer \cite{max} concerning an explicit Wedderburn
decomposition in terms of the Kazhdan--Lusztig basis. We show that the 
subring generated by that basis is nothing but Lusztig's ring $J$; we 
also obtain an analogue of Lusztig's homomorphism from the Iwahori--Hecke 
algebra into $J$; see Section~7. As an application, this gives rise to 
a ``canonical'' homomorphism from the generic Iwahori--Hecke algebra of 
type $B_n$ into the group algebra of the underlying Weyl group. An
explicit example is worked out in Example~\ref{expB2}. In the equal 
parameter case, such a homomorphism was first constructed by Lusztig 
\cite{Lu0}.

We close this introduction with the remark that the results in Sections~2--4 
hold for general Coxeter groups and may be of independent interest. The 
applications to type $B_n$, to be found in Section~5--7, depend on the 
two articles by Bonnaf\'e--Iancu \cite{BI} and Bonnaf\'e \cite{BI2}
(where the left cells and the two-sided cells are determined), but are
otherwise self-contained. 

\section{The basic set-up} \label{sec0a}
We begin by recalling the basic definitions concerning Kazhdan--Lusztig 
cells in the general multi-parameter case.  Let $W$ be a Coxeter group, 
with generating set $S$. (We assume that $S$ is a finite set, but the 
group $W$ may be finite or infinite.) In \cite{Lusztig03}, the parameters 
of the corresponding Iwahori--Hecke algebra are specified by an 
integer-valued weight function. Following a suggestion of Bonnaf\'e 
\cite{BI2}, we can slightly modify Lusztig's definition so as to include 
the more general setting in \cite{Lusztig83} as well (where the parameters
may be contained in a totally ordered abelian group). So let $\Gamma$ be 
an abelian group (written additively) and assume that there is a total 
order $\leq$ on $\Gamma$ compatible with the group structure. (In the 
setting of \cite{Lusztig03}, we take $\Gamma=\Z$ with the natural order.)

Let $A={\Z}[\Gamma]$ be the free abelian group with basis $\{e^\gamma
\mid \gamma \in \Gamma\}$. There is a well-defined ring structure on $A$
such that $e^\gamma e^{\gamma'}=e^{\gamma+ \gamma'}$ for all $\gamma,\gamma'
\in \Gamma$.  (Hence, if $\Gamma=\Z$, then $A$ is nothing but the ring of 
Laurent polynomials in an indeterminate~$e$.) We write $1=e^0 \in A$. 
Given $a\in A$ we denote by $a_\gamma$ the coefficient of $e^\gamma$, so 
that $a= \sum_{\gamma \in \Gamma} a_\gamma e^\gamma$. We let $A_{\geq 0}:=
\langle e^\gamma \mid \gamma \geq 0\rangle_{\Z}$; similarly, we define
$A_{>0}$, $A_{\leq 0}$ and $A_{<0}$. We say that a function
\[ L \colon W\rightarrow\Gamma\] 
is a weight function if $L(ww')= L(w)+L(w')$ whenever we have $l(ww')=
l(w)+l(w')$ where $l\colon W\rightarrow {\N}$ is the usual length function.
(We denote $\N=\{0,1,2,\ldots\}$.) We assume throughout that $L(s)>0$
for all $s\in S$. Let $\bH=\bH(W,S,L)$ be the generic Iwahori--Hecke 
algebra over $A$ with parameters $\{q_s \mid s\in S\}$ where $q_s:=e^{L(s)}$
for $s\in S$.  The algebra $\bH$ is free over $A$ with basis $\{T_w\mid w
\in W\}$, and the multiplication is given by the rule
\[ T_sT_w=\left\{\begin{array}{cl} T_{sw} & \quad \mbox{if $l(sw)>l(w)$},\\
T_{sw}+(q_s-q_s^{-1})T_w & \quad \mbox{if $l(sw)<l(w)$},\end{array}
\right.\]
where $s\in S$ and $w\in W$. (Note that the above elements $T_w$ are 
denoted $\tilde{T}_w$ in \cite{Lusztig83}.) 

For any $a\in A$, we define $\bar{a}:=\sum_{\gamma \in \Gamma} 
a_\gamma e^{-\gamma}$. We extend the map $a \mapsto \bar{a}$ to a ring 
involution $\bH\rightarrow \bH$, $h \mapsto \overline{h}$, by the formula
\[ \overline{\sum_{w \in W} a_w T_w}=\sum_{w \in W} \bar{a}_w
T_{w^{-1}}^{-1} \qquad (a_w \in A).\] 
Now we have a corresponding {\em Kazhdan--Lusztig basis} of $\bH$, which 
we denote by $\{\bC_w\mid w\in W\}$\footnote{Note that this basis is 
denoted by $C_w'$ in \cite{Lusztig83} and by $c_w$ in \cite{Lusztig03}.}.
The basis element $\bC_w$ is uniquely determined by the conditions that 
\[\overline{\bC}_w=\bC_w\qquad \mbox{and}\qquad \bC_w\equiv T_w \quad\bmod 
\bH_{<0},\]
where $\bH_{<0}:=\sum_{w\in W} A_{<0}\, T_w$; see \cite[Prop.~2]{Lusztig83}
and \cite[Theorem~5.2]{Lusztig03}. 

\subsection{Multiplication rules} \label{multrule} 
For any $x,y\in W$, we write
\[ \bC_x\,\bC_y=\sum_{z\in W} h_{x,y,z} \, \bC_z \qquad \mbox{where 
$h_{x,y,z} \in A$ for all $x,y,z\in W$}.\]
An easy induction on $l(x)$ shows that $T_xT_y$ is a linear combination of 
basis elements $T_z$ where $l(z)\leq l(x)+l(y)$. This also implies that 
\[ h_{x,y,z} \neq 0 \qquad \Rightarrow \qquad l(z)\leq l(x)+l(y).\]
We have the following more explicit formula for $s\in S$, $y\in W$ 
(see \cite[\S 6]{Lusztig83}):
\[ \bC_s\,\bC_y = \left\{\begin{array}{ll} \displaystyle{\bC_{sy}+
\sum_{\atop{z\in W}{sz<z<y}} M_{z,y}^s \bC_z} &\quad \mbox{if
$sy>y$},\\(q_s+q_s^{-1})\,\bC_y &\quad \mbox{if $sy<y$},\end{array}\right.\]
where $M_{z,y}^s=\overline{M}_{z,y}^s \in A$ is determined as in 
\cite[\S 3]{Lusztig83} and $\leq$ denotes the Bruhat--Chevalley order. In 
particular, we have 
\[ h_{s,y,z} \neq 0 \quad \Rightarrow \quad z=y>sy \quad \mbox{or}\quad
z=sy>y\quad \mbox{or}\quad sz<z<y<sy.\]

\subsection{The Kazhdan--Lusztig pre-orders} \label{klpre} 
As in \cite[\S 8]{Lusztig03}, we write $x \leftarrow_{\cL} y$ if there 
exists some $s\in S$ such that $h_{s,y,x}\neq 0$, that is, $\bC_x$ occurs 
in $\bC_s\, \bC_y$ (when expressed in the $\bC$-basis). The Kazhdan--Lusztig 
left pre-order $\leq_{\cL}$ is the relation on $W$ generated by 
$\leftarrow_{\cL}$, that is, we have $x\leq_{\cL} y$ if there exists a 
sequence $x=x_0,x_1, \ldots,x_k=y$ of elements in $W$ such that $x_{i-1}
\leftarrow_{\cL} x_i$ for all~$i$. The equivalence relation associated 
with $\leq_{\cL}$ will be denoted by $\sim_{\cL}$ and the corresponding 
equivalence classes are called the {\em left cells} of $W$.  

Similarly, we can define a pre-order $\leq_{\cR}$ by considering 
multiplication by $\bC_s$ on the right in the defining relation. The 
equivalence relation associated with $\leq_{\cR}$ will be denoted by 
$\sim_{\cR}$ and the corresponding equivalence classes are called the 
{\em right cells} of $W$.  We have 
\[ x \leq_{\cR} y \quad \Leftrightarrow \quad x^{-1} \leq_{\cL} y^{-1}.\]
This follows by using the antiautomorphism $\flat\colon \bH\rightarrow 
\bH$ given by $T_w^\flat=T_{w^{-1}}$; we have $\bC_w^\flat=\bC_{w^{-1}}$ 
for all $w\in W$; see \cite[5.6]{Lusztig03}. Thus, any statement concerning 
the left pre-order relation $\leq_{\cL}$ has an equivalent version for the
right pre-order relation $\leq_{\cR}$, via $\flat$.  Finally, we define a 
pre-order $\leq_{\cLR}$ by the condition that $x\leq_{\cLR} y$ if there 
exists a sequence $x=x_0,x_1,\ldots, x_k=y$ such that, for each $i \in 
\{1,\ldots,k\}$, we have $x_{i-1} \leq_{\cL} x_i$ or $x_{i-1}\leq_{\cR} x_i$.
The equivalence relation associated with $\leq_{\cLR}$ will be denoted by
$\sim_{\cLR}$ and the corresponding equivalence classes are called the 
{\em two-sided cells} of $W$.  

\subsection{Left cell representations} \label{leftrep} 
Let $\fC$ be a left cell or, more generally, a union of left cells of $W$. 
We define an $\bH$-module by $[\fC]_A:={\fI}_{\fC}/\hat{\fI}_{\fC}$, where 
\begin{align*}
{\fI}_{\fC} &:=\langle \bC_w\mid w\leq_{\cL} 
z\mbox{ for some $z \in\fC$}\rangle_A,\\
\hat{\fI}_{\fC} &:=\langle \bC_w\mid w \not\in \fC, w\leq_{\cL} z
\mbox{ for some $z \in\fC$}\rangle_A.
\end{align*}
Note that, by the definition of the pre-order relation $\leq_{\cL}$, these
are left ideals in $\bH$. Now denote by $c_x$ ($x\in \fC$) the residue
class of $\bC_x$ in $[\fC]_A$. Then the elements $\{c_x\mid x\in \fC\}$ 
form an $A$-basis of $[\fC]_A$ and the action of $\bC_w$ ($w \in W$) is 
given by the formula
\[ \bC_w.c_x=\sum_{y \in \fC} h_{w,x,y}\, c_y.\]
Assume now that $\fC$ is a finite set and write $\fC=\{x_1,\ldots,x_d\}$.
Let $\{c_1,\ldots,c_d\}$ be the corresponding standard basis of 
$[\fC]_A$, where $c_i=c_{x_i}$ for all~$i$. Then we obtain a matrix 
representation
\[ \fX_{\fC}\colon \bH \rightarrow M_{d}(A)\quad \mbox{where} \quad
\fX_{\fC}(\bC_w)=\big( h_{w,x_j,x_i}\big)_{1\leq i,j \leq d}\]
for any $w\in W$. Thus, $h_{w,x_j,x_i}$ is the $(i,j)$-coefficient of the
matrix $\fX_{\fC}(\bC_w)$. 

Although this will not play a role in this paper, we mention that, for 
various reasons, it is sometimes more convenient\footnote{Here is a simple 
example to illustrate this point: Let $\fC=\{1\}$ be the left cell consisting 
of the identity element of $W$. Then $[\fC]_A$ affords the representation 
$T_s \mapsto -q_s^{-1}$ ($s\in S$) and $[\fC]_A^\delta$ affords the 
representation $T_s \mapsto q_s$ ($s\in S$). Specializing $q_s \mapsto 1$, 
we obtain the sign and the unit representation of $W$, respectively. It is 
sometimes more natural to associate the unit representation with the left 
cell $\{1\}$; so one should work with $[\fC]_A^\delta$ in this case. 
Especially, this can be seen in \cite[Chap.~21]{Lusztig03} where Lusztig
works with $[\fC]_A^\delta$ throughout.} to twist the action of $\bH$ on 
$[\fC]_A$ by the $A$-algebra automorphism 
\[ \delta \colon \bH\rightarrow \bH,\qquad T_s \mapsto -T_s^{-1} \quad
(s \in S).\]
We shall often write $h^\delta$ instead of $\delta(h)$. As in 
\cite[21.1]{Lusztig03}, we define  a new $\bH$-module by taking the
same underlying $A$-module as before, but where the action is given by
the formulas
\[ \bC_w^{\, \delta}.c_x=\sum_{y \in \fC} h_{w,x,y}\, c_y \qquad (w\in W,
x \in \fC).\]
We denote this new $\bH$-module by $[\fC]_A^\delta$. It is readily checked 
that $[\fC]_A^\delta=\delta({\fI}_{\fC})/\delta (\hat{\fI}_{\fC})$. 

\addtocounter{thm}{3}
\begin{rem} \label{otherC} We have a unique ring involution $j \colon \bH
\rightarrow \bH$ such that $j(e^\gamma)=e^{-\gamma}$ for $\gamma \in \Gamma$
and $j(T_w)=(-1)^{l(w)}T_{w}$ for $w\in W$. Then $j$ commutes with
$\delta$ and the composition $j \circ \delta$ is nothing but the involution
$h \mapsto \overline{h}$ on $\bH$; see \cite[\S 6]{Lusztig83}. Thus, we have 
\[ \delta(\bC_w)=\delta(\overline{\bC}_w)=j(\bC_w)\qquad\mbox{for any 
$w \in W$}.\]
This observation can be used to obtain formulas for $\delta(h)$ ($h \in \bH$)
which would be difficult to compute using the definition of $\delta$. For 
example, we obtain 
\[ \delta(\bC_w)=j(\bC_w)=(-1)^{l(w)}T_w+\sum_{\atop{y\in W}{y<w}} 
(-1)^{l(y)} \,\overline{P}_{y,w}^{\,*} \,T_y\]
for any $w\in W$. 
\end{rem} 

We shall be interested in the following property. 

\begin{defn} \label{eqcells} Let $\fC$ and $\fC_1$ be left cells or, more 
generally, be unions of left cells of $W$.  We write $\fC \approx\fC_1$, 
if there exists a bijection $\fC\stackrel{\sim}{\rightarrow} {\fC}_1$, 
$x \mapsto x_1$, such that the following condition is satisfied:
\begin{equation*}
h_{w,x,y}=h_{w,x_1,y_1}\qquad\mbox{for all $w\in W$ and all $x,y\in\fC$}.
\tag{$\heartsuit$}
\end{equation*}
This means that the $\bH$-modules $[\fC]_A$ and $[\fC_1]_A$ 
are not only isomorphic, but even the action of any $\bC_w$ ($w\in W$) is 
given by exactly the same formulas with respect to the standard bases of
$[\fC]_A$ and $[\fC]_A$, respectively. A similar remark applies, of course, 
to the $\bH$-modules $[\fC]_A^\delta$ and $[\fC_1]_A^\delta$. Note 
that, in order to verify that ($\heartsuit$) holds, it is enough to consider 
the case where $w=s\in S$ (since the elements $\bC_s$, $s \in S$, generate
$\bH$ as an $A$-algebra).
\end{defn}

\begin{exmp} \label{expA} Let $W=\fS_n$ be the symmetric group, with
generating set $S=\{s_1, \ldots,s_{n-1}\}$ where $s_i=(i,i+1)$ for 
$1\leq i \leq n-1$.  Let $\Gamma=\Z$ with its natural order, and set
$q:=e^1$. Then $A={\Z}[\Gamma]={\Z}[q,q^{-1}]$ is the ring of Laurent
polynomials in an indeterninate~$q$. Let $L\colon \fS_n \rightarrow \Z$
be the weight function given by $L(s_i)=1$ for $1\leq i \leq n-1$, and 
denote by $H(\fS_n)$ the corresponding Iwahori--Hecke algebra over $A$. 
Thus, we have the following diagram specifying the generators, relations 
and parameters:
\begin{center}
\begin{picture}(300,44)
\put( 40, 20){$A_{n-1}$}
\put( 45, 03){$\{q_s\}$:}
\put(101, 20){\circle{10}}
\put(106, 20){\line(1,0){29}}
\put(140, 20){\circle{10}}
\put(145, 20){\line(1,0){20}}
\put(175, 17){$\cdot$}
\put(185, 17){$\cdot$}
\put(195, 17){$\cdot$}
\put(205, 20){\line(1,0){20}}
\put(230, 20){\circle{10}}
\put( 96, 32){$s_1$}
\put( 98, 03){$q$}
\put(136, 32){$s_2$}
\put(137, 03){$q$}
\put(222, 32){$s_{n-1}$}
\put(227, 03){$q$}
\end{picture}
\end{center}
The classical Robinson--Schensted correspondence associates with
each element $\sigma \in \fS_n$ a pair of standard tableaux $(A(\sigma),
B(\sigma))$ of the same shape. For any partition $\nu$ of $n$, we set
\[ \fR_\nu:=\{\sigma \in \fS_n \mid \mbox{ $A(\sigma)$, $B(\sigma)$ have
shape $\nu$}\}.\]
Thus, we have $\fS_n=\coprod_{\nu } \fR_\nu$ where $\nu$ runs
over all partitions of~$n$. Then the following hold.
\begin{itemize}
\item[(a)] {\em For a fixed standard tableau $T$, the set $\{\sigma \in
\fS_n \mid B(\sigma)=T\}$ is a left cell of $\fS_n$ and $\{\sigma \in
\fS_n \mid A(\sigma)=T\}$ is a right cell of $\fS_n$. Furthermore, all
left cells and all right cells arise in this way}.
\item[(b)] {\em Let $\fC,\fC_1$ be left cells and assume that $\fC\subseteq 
\fR_\nu$, $\fC_1\subseteq \fR_{\nu_1}$. Then we have $\fC\approx\fC_1$ 
if and only if $\nu=\nu_1$. The required bijection from $\fC$ onto $\fC_1$ 
can be explicitly described in terms of the ``star'' operation defined in
\cite[\S 4]{KaLu}}.
\end{itemize}
These statements were first proved by Kazhdan--Lusztig \cite[\S 5]{KaLu}. 
(See also Ariki \cite{Ar}.) It is actually shown there that the bijection 
$x \mapsto x_1$ is determined by the condition that $x\in \fC$ and $x_1\in
\fC_1$ lie in the same right cell. In Proposition~\ref{uniqueeq}, we will 
see that this property automatically follows from some general principles. 
\end{exmp}

\begin{lem} \label{cellaut} Let $\varphi \colon W \rightarrow W$ be
a group automorphism such that $\varphi(S)=S$ and $q_{\varphi(s)}=q_s$
for all $s\in S$. Let $\fC,\fC_1$ be left cells of $W$. Then $\varphi(\fC),
\varphi(\fC_1)$ are left cells and we have $\fC \approx \fC_1$ if and only
if $\varphi(\fC) \approx \varphi(\fC_1)$.
\end{lem}

\begin{proof} Our assumptions imply that $\varphi$ induces an $A$-algebra 
automorphism $\tilde{\varphi} \colon \bH \rightarrow \bH$ such that 
$\tilde{\varphi}(T_w) =T_{\varphi(w)}$ for all $w\in W$. It is readily 
checked that $\tilde{\varphi}$ commutes with the involution $h \mapsto 
\overline{h}$ of $\bH$. This implies that 
\[ \tilde{\varphi}(\bC_w)=\bC_{\varphi(w)} \qquad \mbox{for all $w\in W$}.\]
Consequently, we also have $h_{x,y,z}=h_{\varphi(x),\varphi(y),\varphi(z)}$ 
for all $x,y,z\in W$. By the definition of left cells, this yields that 
$\varphi$ preserves the partition of $W$ into left cells and that
we have $\fC \approx \fC_1$ if and only if $\varphi(\fC) \approx 
\varphi(\fC_1)$.
\end{proof}

\begin{lem} \label{cellw0} Assume that $W$ is finite and let $w_0\in W$
be the unique element of maximal length. Let $\fC$ and $\fC_1$ be left
cells such that $\fC \approx \fC_1$. Then we also have $\fC w_0\approx
\fC_1 w_0$ and $w_0\fC \approx w_0\fC_1$. (Note that $\fC w_0$, $\fC_1w_0$ 
and $w_0\fC$, $w_0\fC_1$ are left cells; see \cite[Cor.~11.7]{Lusztig03}.)
\end{lem}

\begin{proof} First we prove that $\fC w_0\approx \fC_1w_0$. Let 
$\fC\stackrel{\sim}{\rightarrow} {\fC}_1$, $x \mapsto x_1$, be a bijection 
such that ($\heartsuit$) holds; see Definition~\ref{eqcells}. In particular, 
this means that $h_{s,x,y}= h_{s,x_1,y_1}$ for all $s\in S$ and $x,y\in \fC$. 

Now recall the formula for multiplication by $\bC_s$ from (\ref{multrule}). 
That formula shows that, for any $s\in S$ and any $x\in \fC$, we have $sx<x$ 
if and only if $sx_1<x_1$. Furthermore, by \cite[Prop.~11.6]{Lusztig03}, we 
have
\[M_{xw_0,yw_0}^s=-(-1)^{l(x)+l(y)} M_{y,x}^s\qquad\mbox{if $sy<y<x<sx$}.\]
Hence we obtain 
\[h_{s,xw_0,yw_0}=h_{s,x_1w_0,y_1w_0} \quad \mbox{for all $s\in S$ and 
$x,y\in \fC$}.\]
Consequently, ($\heartsuit$) holds for the bijection $\fC w_0 \rightarrow 
\fC_1w_0$, $xw_0 \mapsto x_1w_0$.  Now consider the group automorphism
$\varphi\colon W\rightarrow W$ given by $\varphi(w)=w_0ww_0$. It is
well-known that $\varphi(S)=S$. Furthermore, since $s\in S$ and
$\varphi(s)$ are conjugate, we have $q_s=q_{\varphi(s)}$. Hence,
Lemma~\ref{cellaut} shows that $w_0\fC w_0$ and $w_0\fC_1w_0$ are left 
cells such that $w_0\fC w_0 \approx w_0\fC_1 w_0$. Hence the previous 
argument shows that  $w_0\fC=(w_0\fC w_0)w_0\approx (w_0\fC_1w_0)w_0=w_0
\fC_1$.
\end{proof}

We close this section with some results which show that, under suitable
hypotheses, a bijection $\fC\stackrel{\sim}{\rightarrow} {\fC}_1$, 
$x \mapsto x_1$, satisfying ($\heartsuit$) automatically respects the 
right cells of $W$. (These results will also play an important role in  
Section~\ref{sec-jB}.) Let us assume throughout that $W$ is a finite group.
Since the group $\Gamma$ is totally ordered, $A={\Z}[\Gamma]$ is easily seen 
to be an integral domain. Let $K$ be the field of fractions of ${\R}[\Gamma]
\supseteq A$. By extension of scalars, we obtain a $K$-algebra 
$\bH_K=K\otimes_{A}\bH$. 

\begin{rem} \label{semis} The algebra $\bH_K$ is split semisimple.
\end{rem}

\begin{proof} The fact that $\bH_K$ is semisimple relies on two ingredients:
firstly, ${\R}W$ (the group algebra of $W$ over $\R$) is known to be split 
semisimple and, secondly, ${\R}[\Gamma]\otimes_A \bH$ specializes to ${\R}W$, 
via the ring homomorphism $\theta \colon {\R}[\Gamma] \rightarrow \R$ such that 
$\theta(e^\gamma)=1$ for all $\gamma \in \Gamma$. Then it remains to use
known results on splitting fields; see \cite[\S 9]{ourbuch} and the
references there. For more details, see \cite[Remark~3.1]{GeIa05}.
\end{proof}

Let $\Irr(\bH_K)$ be the set of irreducible characters of $\bH_K$. We write 
this set in the form
\[ \Irr(\bH_K)=\{\chi_\lambda \mid \lambda\in \Lambda\},\]
where $\Lambda$ is some finite indexing set. The algebra $\bH_K$ is 
{\em symmetric} with respect to the trace function $\tau \colon \bH_K 
\rightarrow K$ defined by $\tau(T_1)=1$ and $\tau(T_w)=0$ for $1\neq w\in W$;
see \cite[\S 8.1]{ourbuch}. The fact that $\bH_K$ is split semisimple yields 
that
\[ \tau=\sum_{\lambda\in \Lambda} \frac{1}{c_\lambda} \, \chi_\lambda
\qquad \mbox{where $0 \neq c_\lambda \in {\R}[\Gamma]$};\]
see \cite[\S 7.2 and 9.3.5]{ourbuch}. The elements $c_\lambda$ are called 
the {\em Schur elements}. 

For any $\lambda \in \Lambda$, let us denote by $\fX_\lambda \colon H_K
\rightarrow M_{d_\lambda}(K)$ a matrix representation with character 
$\chi_\lambda$. Let $\fX_\lambda^{ij}(h)$ denote the $(i,j)$-coefficient
of $\fX_\lambda(h)$ for any $h \in \bH_K$. By Wedderburn's theorem, the 
algebra $\bH_K$ is abstractly isomorphic to the direct sum of the
matrix algebras $M_{d_\lambda}(K)$ ($\lambda \in \Lambda$). Since 
$\bH_K$ is symmetric, this isomorphism can be described explicitly:

\begin{prop}[Explicit Wedderburn decomposition] \label{wedder}  Let
$\cB$ be any basis of $\bH_K$  and ${\cB}^\vee=\{b^\vee \mid b \in\cB\}$ 
the dual basis with respect to $\tau$.  We set 
\[ E_{\lambda}^{ij}=\frac{1}{c_\lambda}\sum_{b \in \cB} \fX_\lambda^{ji}
(b)\, b^\vee \qquad \mbox{for any $\lambda \in \Lambda$, $1\leq i,j \leq 
d_\lambda$}.\]
Then $\fX_\lambda(E_\lambda^{ij}) \in M_{d_\lambda}(K)$ is the matrix with 
$(i,j)$-coefficient $1$ and coefficient $0$ otherwise.  Furthermore, if 
$\mu\neq \lambda$, we have $\fX_\mu^{kl}(E_\lambda^{ij})=0$ for all 
$1 \leq k,l \leq d_\mu$. In particular, the elements 
\[\{E_\lambda^{ij} \mid \lambda \in \Lambda, 1\leq i,j \leq d_\lambda\}\]
form a basis of $\bH_K$.
\end{prop}

(For a proof, see \cite[Prop.~7.2.7]{ourbuch}, for example.)

We want to apply the above result to the Kazhdan--Lusztig basis $\cB:=
\{\bC_w \mid w\in W\}$. The dual basis can be described as follows. We set 
\[ \bD_{z^{-1}}:=(-1)^{l(z)+l(w_0)} \bC_{zw_0}^{\,\delta}\, T_{w_0} \qquad
\mbox{for any $z \in W$}.\]
where $w_0 \in W$ is the unique element of maximal length. Then we have
\[ \tau(\bC_w\bD_{z^{-1}})=\left\{\begin{array}{cl} 1 & \quad 
\mbox{if $w=z$},\\ 0 & \quad \mbox{if $w\neq z$};\end{array}\right.\]
see \cite[Prop.~11.5]{Lusztig03}. Hence we have $\bC_w^\vee=\bD_{w^{-1}}$ for 
all $w\in W$. In particular, the structure constants of $\bH$ can be 
expressed by 
\[ h_{x,y,z}=\tau(\bC_x\bC_y\bD_{z^{-1}})\qquad\mbox{for all $x,y,z\in W$}.\]
This immediately yields that 
\[ \bC_x\bD_{y^{-1}}=\sum_{w \in W} h_{w,x,y}\, \bD_{w^{-1}}\qquad
\mbox{for any $x,y \in W$}.\]
The following two results were observed by Neunh\"offer in his thesis 
\cite[Kap.~VI, \S 4]{max}. For any left cell $\fC$, denote by 
$\chi_{\fC}$ the character afforded by the left cell module 
$[\fC]_K:= K\otimes_A [\fC]_A$ of $\bH_K$.

\begin{lem}[Neunh\"offer] \label{neun1} Let $\fC$ be a left cell such that
$\chi_{\fC}\in \Irr(\bH_K)$. Writing $\fC=\{x_1,\ldots,x_{d}\}$ and
using the notation in (\ref{leftrep}), we have 
\[ E_\lambda^{ij}=\frac{1}{c_\lambda}\,\bC_{x_i}\,\bD_{x_j^{-1}} \qquad 
\mbox{for $1\leq i,j\leq d$},\]
where $\lambda \in \Lambda$ is such that $\chi_{\lambda}=\chi_{\fC}$ and
where we take $\fX_\lambda=\fX_{\fC}$. In particular, we have 
$\fX_\lambda(\bC_{x_i}\bD_{x_j})\neq 0$ and $\fX_\mu(\bC_{x_i}\bD_{x_j})=0$ 
for any $\mu \in \Lambda \setminus \{\lambda\}$.
\end{lem}

\begin{proof} Since Neunh\"offer only considers the case of the symmetric
group, we give a proof here. By the formula in Proposition~\ref{wedder}, we
have 
\[ E_\lambda^{ij}=\frac{1}{c_\lambda}\sum_{w\in W}
\fX_{\lambda}^{ji}(\bC_w) \bD_{w^{-1}}\qquad \mbox{where $\fX_\lambda:=
\fX_{\fC}$}.\]
We have observed in (\ref{leftrep}) that $\fX_{\lambda}^{ji}(\bC_w)=
h_{w,x_i,x_j}$.  Hence we have 
\[ E_\lambda^{ij}=\frac{1}{c_\lambda}\sum_{w\in W} h_{w,x_i,x_j}\, 
\bD_{w^{-1}}=\frac{1}{c_\lambda}\, \bC_{x_i}\, \bD_{x_j^{-1}},\]
as desired. The remaining statements are clear by Proposition~\ref{wedder}.
\end{proof}

\begin{lem}[Neunh\"offer] \label{neun2} Let $\fC,\fC_1$ be two left cells
of $W$ such that $\fC \approx \fC_1$. Let $\fC \stackrel{\sim}{\rightarrow} 
\fC_1$, $x \mapsto x_1$, be a bijection such that condition $(\heartsuit)$ in 
Definition~\ref{eqcells} holds. Then we have 
\[ \bC_x\bD_{y^{-1}}=\bC_{x_1}\bD_{y_1^{-1}} \qquad \mbox{for all 
$x,y \in \fC$}.\]
\end{lem}

\begin{proof} Condition ($\heartsuit$) means that $h_{w,x,y}=h_{w,x_1,y_1}$
for all $w\in W$ and all $x,y \in \fC$.  Hence we also have 
\[ \bC_x\bD_{y^{-1}}=\sum_{w \in W} h_{w,x,y} \, \bD_{w^{-1}} =
\sum_{w \in W} h_{w,x_1,y_1} \, \bD_{w^{-1}} =\bC_{x_1}\bD_{y_1^{-1}},\]
as required.
\end{proof}

\begin{prop} \label{uniqueeq} In the above setting, let $\fC,\fC_1$ be 
left cells such that $\chi_{\fC}=\chi_{\fC_1}\in \Irr(\bH_K)$ 
and $\fC \approx \fC_1$. Let $\fC \stackrel{\sim}{\rightarrow} \fC_1$, 
$x\mapsto x_1$, be a bijection such that condition $(\heartsuit)$ in 
Definition~\ref{eqcells} holds.  Then we have $x\sim_{\cR} x_1$ for 
any $x\in \fC$.
\end{prop}

\begin{proof}  Let $x\in \fC$. We show that $x_1\leq_{\cR} x$. To see
this, we argue as follows. Choose an enumeration of the elements in 
$\fC$ where $x$ is the first element. Consider the corresponding matrix 
representation $\fX_{\fC}$. By Lemma~\ref{neun1}, $\fX_{\fC}
(\bC_x \bD_{x^{-1}})$ is a matrix with a non-zero coefficient at position 
$(1,1)$ and coefficient $0$ otherwise. Consequently, some coefficient 
in the first row of $\fX_{\fC}(\bC_x)$ must be non-zero.  Using 
(\ref{leftrep}) we see that there exists some $y \in \fC$ such that 
$h_{x,y,x} \neq 0$. Then, by ($\heartsuit$), we have $h_{x,y_1,x_1}=
h_{x,y,x}\neq 0$ and so $x_1 \leq_{\cR} x$, as claimed. 

We now apply a similar discussion to the left cell $\fC_1$ and the
element $x_1$. Working with the representation $\fX_{\fC_1}$, we see
that there exists some $z_1\in \fC_1$ such that $h_{x_1,z_1,x_1}\neq 0$.
But then we have $h_{x_1,z,x}=h_{x_1,z_1,x_1}\neq 0$ and so 
$x \leq_{\cR} x_1$. Hence we conclude that $x \sim_{\cR} x_1$.
\end{proof}

\begin{exmp} \label{expA1} Let us consider once more the case where
$W=\fS_n$, as in Example~\ref{expA}. It is shown by Kazhdan--Lusztig
\cite{KaLu} that 
\[ \chi_{\fC}\in \Irr(\bH(\fS_n)_K) \qquad \mbox{for any left cell
$\fC\subseteq \fS_n$}.\]
Now let $\fC,\fC_1$ be left cells such that $\fC \approx \fC_1$; see
Example~\ref{expA}(b) for a characterisation of this condition. Let 
$\fC \stackrel{\sim}{\rightarrow} \fC_1$, $x\mapsto x_1$, be a bijection
such that condition $(\heartsuit)$ in Definition~\ref{eqcells} holds. Then,
by Proposition~\ref{uniqueeq}, we have $x\sim_{\cR} x_1$ for any $x\in \fC$.
However, the Robinson--Schensted correspondence shows that two elements
which lie in the same left cell and in the same right  cell must be
equal. Hence the element $x_1 \in \fC_1$ is uniquely determined by
the condition that $x \sim_{\cR} x_1$.
\end{exmp}

\section{On the induction of Kazhdan--Lusztig cells} \label{sec-indu}

In \cite{myind}, it is shown that the Kazhdan--Lusztig basis of $\bH$
behaves well with respect to parabolic subalgebras. One of the aims of 
this section is to show that the relation ``$\approx$'' in 
Definition~\ref{eqcells} also behaves well. Corollary~\ref{ind1} (obtained 
at the end of this section) will play a crucial role in the proof of 
Theorem~\ref{eqcellsB}. In a different direction, the techniques 
developed in this section lay the foundations for the discussion of the  
relative version of ($\spadesuit$).

We keep the basic set-up of the previous section. Let us fix a subset $I 
\subseteq S$ and consider the corresponding parabolic subgroup $W_I=
\langle I\rangle \subseteq W$.  Let $\bH_I=\langle T_w \mid w \in W_I 
\rangle_{A}$ be the parabolic subalgebra corresponding to $W_I$. It is clear
by the definition that, for any $w\in W_I$, we have that $\bC_w$ computed 
inside $\bH_I$ is the same as $\bC_w$ computed in $\bH$. 

The following definitions already appear, in a somewhat different form, 
in the work of Barbasch--Vogan \cite[\S 3]{BV}.

\subsection{Relative Kazhdan--Lusztig pre-orders} \label{relpre}
Given $x,y \in W$, we write $x \leftarrow_{\cL,I} y$ if there exists some 
$s\in I$ such that $h_{s,y,x}\neq 0$, that is, $\bC_x$ occurs in $\bC_s\, 
\bC_y$ (when expressed in the $\bC$-basis). Let $\leq_{\cL,I}$ be the 
pre-order relation on $W$ generated by $\leftarrow_{\cL,I}$, that is, we have 
$x\leq_{\cL,I} y$ if there exists a sequence $x=x_0,x_1,\ldots,x_k=y$ of 
all elements in $W$ such that $x_{i-1} \leftarrow_{\cL,I} x_i$ for all~$i$.  
The equivalence relation associated with $\leq_{\cL,I}$ will be denoted by 
$\sim_{\cL,I}$ and the corresponding equivalence classes are called the 
{\em relative left cells} of $W$ with respect to $I$. Note that the 
restriction of $\leq_{\cL,I}$ to $W_I$ is nothing but the usual left 
pre-order on $W_I$.

Similarly, we can define a pre-order $\leq_{\cR,I}$ by considering 
multiplication by $\bC_s$ ($s\in I$) on the right in the defining condition.
The equivalence relation associated with $\leq_{\cR,I}$ will be denoted 
by $\sim_{\cR,I}$ and the corresponding equivalence classes are called the 
{\em relative right cells} of $W$ (with respect to $I$).  We have 
\[ x \leq_{\cR,I} y \quad \Leftrightarrow \quad x^{-1} \leq_{\cL,I} y^{-1}.\]
This follows, as before, by using the antiautomorphism $\flat\colon \bH
\rightarrow \bH$ given by $T_w^\flat=T_{w^{-1}}$. Finally, we define a 
pre-order $\leq_{\cLR,I}$ by the condition that $x\leq_{\cLR,I} y$ if there 
exists a sequence $x=x_0,x_1,\ldots, x_k=y$ such that, for each $i \in 
\{1,\ldots,k\}$, we have $x_{i-1} \leq_{\cL,I} x_i$ or $x_{i-1}\leq_{\cR,I} 
\leq x_i$. The equivalence relation associated with $\leq_{\cLR,I}$ will be 
denoted by $\sim_{\cLR,I}$ and the corresponding equivalence classes are 
called the {\em relative two-sided cells} of $W$.  

Let $X_I$ be the set of distinguished left coset representatives; we have
\[ X_I=\{w \in W \mid \mbox{$w$ has minimal length in $wW_I$}\}.\]
Furthermore, the map $X_I \times W_I \rightarrow W$, $(x,u) \mapsto xu$, 
is a bijection and we have $l(xu)=l(x)+l(u)$ for all $u\in W_I$ and all
$x\in X_I$. We define a relation ``$\sqsubseteq$'' as follows. Let $x,y\in
X_I$ and $u,v\in W_I$. We write $xu\sqsubset yv$ if $x<y$ (Bruhat--Chevalley 
order) and $u\leq_{\cL,I} v$ (Kazhdan--Lusztig pre-order). We write
$xu\sqsubseteq yv$ if $xu\sqsubset yv$ or $x=y$ and $u=v$. With this
notation, we have the following result.

\addtocounter{thm}{1}
\begin{prop}[See \protect{\cite[Prop.~3.3]{myind}}] \label{mythm} For
any $y\in X_I$, $v\in W_I$ we have
\[ \bC_{yv}=\sum_{\atop{x\in X_I,u\in W_I}{xu\sqsubseteq yv}} 
p_{xu,yv}^*\, T_x\,\bC_u\]
where $p_{yv,yv}^*=1$ and $p_{xu,yv}^*\in A_{<0}$ for $ux\sqsubset yv$.
\end{prop}

For later use, we have to recall the basic ingredients in the construction 
of the polynomials $p_{xu,yv}^*$; we also prove some refinements of the
results in \cite[\S 3]{myind}. Let $y \in X_I$ and $v \in W_I$. Then we 
can write uniquely 
\[ \overline{T_y\bC_v}=T_{y^{-1}}^{-1}\,\bC_v=\sum_{\atop{x\in X_I}{u\in W_I}}
\overline{r}_{xu,yv}\,T_x\,\bC_u\qquad\mbox{where $r_{xu,yv}\in A$}\]
and where only finitely many terms $r_{xu,yv}$ are non-zero. 

\begin{lem} \label{rpol} Let $x,y \in X_I$ and $u,v\in W_I$. Then we
have $r_{xu,yv}=0$ unless $l(xu)<l(yv)$ or $xu=yv$. Furthermore, we have
\[ \overline{r}_{xu,yv}=\sum_{\atop{w\in W_I}{xw\leq y}}
\sum_{\atop{w'\in W_I}{w'\leq w}} \overline{R}_{xw,y}^*\,\tilde{p}_{w',w}\, 
h_{w',v,u}\]
where $\tilde{p}_{w',w} \in A$ are independent of $x,y,u,v$ and the 
$R_{z,y}^*\in A$ are the ``absolute'' $R$-polynomials defined in 
\cite[\S 1]{Lusztig83}.
\end{lem}

\begin{proof} First we establish the above identity. Let us fix $y\in X_I$
and $v \in W_I$.  We can write
\[{T}_{y^{-1}}^{-1}=\sum_{\atop{z\in W}{z\leq y}} \overline{R}_{z,y}^*\,T_z
\qquad (R_{y,y}^*=1).\]
Now let $z\in W$ be such that $T_z$ occurs in the above expression. Then 
we can write $z=xw$ where $x\in X_I$ and $w\in W_I$; note that $x\leq z
\leq y$. Since $l(xw)=l(x)+ l(w)$, we have $T_z=T_x\,T_w$ and so
\[ T_{y^{-1}}^{-1}\,\bC_v=\sum_{\atop{x\in X_I}{x\leq y}}
\sum_{\atop{w\in W_I}{xw\leq y}} \overline{R}_{xw,y}^*\, T_x\,T_w\,\bC_v.\]
Now, by \cite[Theorem~5.2]{Lusztig03}, $\bC_{w}$ is a linear combination
of terms $T_{w'}$ where $w'\leq w$ and the coefficient of $T_w$ is~$1$. 
Hence we can also write $T_w=\sum_{w'} \tilde{p}_{w',w} \bC_{w'}$ where 
$\tilde{p}_{w,w}=1$ and $\tilde{p}_{w',w}=0$ unless $w'\leq w$. Thus, we have 
\begin{align*}
T_{y^{-1}}^{-1}\,\bC_v&=\sum_{\atop{x\in X_I}{x\leq y}}
\sum_{\atop{w\in W_I}{xw\leq y}} \overline{R}_{xw,y}^*
\sum_{\atop{w' \in W_I}{w' \leq w}} \tilde{p}_{w',w}\, T_x\,\bC_{w'}\,\bC_v\\ 
&=\sum_{\atop{x \in X_I}{x\leq y}} \sum_{u \in W_I} \Bigl(
\sum_{\atop{w\in W_I}{xw\leq y}} \sum_{\atop{w'\in W_I}{w'\leq w}} 
\overline{R}_{xw,y}^*\,\tilde{p}_{w',w}\, h_{w',v,u}\Bigr)T_x\,\bC_u.
\end{align*}
This yields the desired identity. Now assume that $r_{xu,yv}\neq 0$.
Then there exist $w,w'\in W_I$ such that $w'\leq w$, $xw\leq y$ and
$h_{w',v,u}\neq 0$. The latter condition certainly implies that
$l(u) \leq l(w')+l(v)$; see (\ref{multrule}). Combining this with the 
inequalities $l(w')\leq l(w)$ and $l(xw)\leq l(y)$, we obtain $l(xu)
\leq l(yv)$, as desired. Furthermore, if equality holds, then equality 
holds in all intermediate inequalities, and so we must have $w'=w$, 
$xw=y$ and, hence, $w'=w=1$. Since $h_{1,v,u}\neq 0$, this also yields
$u=v$, as desired.  
\end{proof}

Now the arguments in the proofs of Lemma~3.2 and Proposition~3.3 in 
\cite{myind} (which themselves are an adaptation of the proof of Lusztig 
\cite[Prop.~2]{Lusztig83}) show that the family  of elements
\[ \{p_{xu,yv}^* \mid x,y \in X_I, u,v \in W_I,xu \sqsubseteq yv\}\]
is  uniquely determined by the following three conditions:
\begin{align*}
p_{yv,yv}^* &= 1, \tag{KL1}\\
p_{xu,yv}^* & \in A_{<0} \qquad \mbox{if $xu\sqsubset yv$},\tag{KL2}\\
\overline{p}^{*}_{xu,yv}-p^{*}_{xu,yv} &= \sum_{\atop{z\in X_I,
w\in W_I}{xu\sqsubset zw\sqsubseteq yv}} r_{xu,zw}\,p^{*}_{zw,yv}
\qquad \mbox{if $xu\sqsubset yv$}.\tag{KL3}
\end{align*}
The arguments in [{\em loc.\ cit.}] provide an inductive procedure for
solving the above system of equations.

The following result yields a further property of the elements $p_{xu,yv}^*$. 

\begin{lem} \label{lem-rel1} Let $x,y\in X_I$ and $u,v\in W_I$. Then 
$p_{xu,yv}^*=0$ unless $xu \leq yv$ (Bruhat--Chevalley order).
\end{lem}

\begin{proof} First we claim that $p_{xu,yv}^*=0$ unless $xu=yv$ or $l(xu)<
l(yv)$. To prove this, we argue as follows. We have seen in Lemma~\ref{rpol} 
that $r_{xu,yv}=0$ unless $xu=yv$ or $l(xu)<l(yv)$. Following the inductive 
procedure for solving the system of equations given by (KL1)--(KL3) above, 
we see that we also must have $p_{xu,yv}^*=0$ unless $xu=yv$ or $l(xu)<l(yv)$.

Now let $x,y\in X_I$ and $u,v \in W_I$ be such that $l(xu) \leq l(yv)$ (with 
equality only for $xu=yv$). We want to prove that $p_{xu,yv}^* =0$ unless 
$xu \leq yv$. We proceed by induction on $l(yv)-l(xu)$. If $l(xu)=l(yv)$, 
then $xu=yv$ and $p_{yv,yv}^*=1$. Now assume that $l(xu)<l(yv)$ and 
$p_{xu,yv}^*\neq 0$. By the proof of \cite[Prop.~3.3]{myind}, we have 
\[ 0\neq p_{xu,yv}^*=P_{xu,yv}^*-\sum_{u<u_1} p_{xu_1,yv}^*\,P_{u,u_1}^*.\]
Now, if $P_{xu,yv}^*\neq 0$, then it is well-known that $xu\leq yv$,
as required. On the other hand, if there is some $u_1\in W$ such that
$u<u_1$ and $p_{xu_1,yv}^*\,P_{u,u_1}^* \neq 0$, then we have $xu_1\leq yv$ 
by induction, and so $xu\leq xu_1\leq yv$.
\end{proof}

\begin{cor} \label{cor-rel2} Let $y\in X_I$ and $v\in W_I$. 
\begin{itemize}
\item[(a)] $\bC_{yv}$ is a linear combination of $T_y\,\bC_v$ and terms
$T_x\,\bC_u$ where $x\in X_I$ and $u\in W_I$ are such that $x<y$, 
$u \leq_{\cL,I} v$ and $xu<yv$. 
\item[(b)] Conversely, $T_y\,\bC_v$ is a linear combination of $\bC_{yv}$ and 
terms $\bC_{xu}$ where $x\in X_I$ and $u\in W_I$ are such that $x<y$, $u 
\leq_{\cL,I} v$ and $xu<yv$.
\end{itemize}
\end{cor}

\begin{proof} (a) This is just a restatement of Proposition~\ref{mythm}, 
taking into account the additional information in Lemma~\ref{lem-rel1}.

(b) Let $w\in W$ and set $B_w:=T_y\bC_v$ where $y \in X_I$ and $v \in W_I$
are such that $w=yv$. Then $\{B_{w}\mid w \in W\}$ is a basis of $\bH$ 
and the formula in Proposition~\ref{mythm} describes the base change 
from the $\bC_w$-basis to the $B_w$-basis. By an easy induction on $l(w)$,
we can invert these formulas. (Note that the base change takes place inside 
the finite sets $\{w \in W \mid l(w)\leq n\}$ for $n=0,1,2,\ldots$.)
Hence we obtain expressions for the elements in the $B_w$-basis in terms of 
the $\bC_w$-basis. The terms arising in these expressions must satisfy 
conditions which are analogous to those in (a). 
\end{proof}

%
%
%

Recall that, if $V$ is any $\bH_I$-module, then 
\[\Ind_I^S(V):=\bH \otimes_{\bH_I} V\]
is an $\bH$-module, called the {\em induced module}; see, for example,
\cite[\S 9.1]{ourbuch}. If $V$ is free over $A$ with basis $\{v_\alpha
\mid \alpha \in \cA\}$, then $\Ind_I^S(V)$ is free with basis 
$\{T_x \otimes v_\alpha \mid x \in X_I,\alpha \in \cA\}$.

\begin{thm}[See \protect{\cite[Theorem~1]{myind}}] \label{inducell}
Let $\fC$ be a left cell of $W_I$. Then the set $X_I\fC$ is a union of 
left cells of $W$. We have an isomorphism of $\bH$-modules
\[ [X_I\fC]_A\stackrel{\sim}{\rightarrow} \Ind_I^S([\fC]_A), 
\qquad c_{yv} \mapsto \sum_{\atop{x \in X_I,u\in \fC}{xu\sqsubseteq yv}}
p_{xu,yv}^*\,\big(T_x \otimes c_u\big),\]
where $\{c_{yv} \mid y \in X_I,v \in \fC\}$ is the standard basis of
$[X_I\fC]_A$ and $\{c_u \mid u \in \fC\}$ is the standard basis
of $[\fC]_A$. 
\end{thm}

\begin{proof} The fact that $X_I\fC$ is a union of left cells is proved
in \cite[\S 4]{myind}. Since the statement concerning $[X_I\fC]_A$ is
not explicitly mentioned in [{\em loc. cit.}], let us give the details here.
Recall that $[X_I\fC]_A=\fI_{X_I\fC}/\hat{\fI}_{X_I\fC}$ where
\begin{align*}
\fI_{X_I\fC}&=\Bigl\langle \bC_{xu} \;\Big|\begin{array}{l} x \in X_I,
u \in W_I, ux \leq_{\cL} vy\\ \mbox{for some $y\in X_I,v\in\fC$}\end{array}
\Bigr \rangle_A,\\
\hat{\fI}_{X_I\fC}&=\Bigl\langle \bC_{xu}\;\Big|\begin{array}{l} x \in X_I,
u \in W_I, ux \not\in X_I\fC, ux \leq_{\cL} vy, \\ 
\mbox{for some $y\in X_I,v\in\fC$}\end{array}\Bigr \rangle_A.
\end{align*}
Now, for any $x,y \in X_I$ and $u,v\in W_I$, we have the implication
\[ xu \leq_{\cL} yv \quad \Rightarrow \quad u \leq_{\cL,I} v;\]
see \cite[\S4]{myind}. On the other hand, we have $xu\leq_{\cL} u$ for
any $x\in X_I$ and $u \in W_I$ (since $l(xu)=l(x)+l(u)$). These two 
relations readily imply that we have 
\begin{align*}
\fI_{X_I\fC}&=\langle \bC_{xu} \mid x \in X_I,u\in W_I,u \leq_{\cL,I} 
v \mbox{ for some $v \in \fC$}\rangle_A,\\
\hat{\fI}_{X_I\fC}&=\langle \bC_{xu}\mid x \in X_I,u \in W_I, u\not\in\fC,
u \leq_{\cL,I} v \mbox{ for some $v \in \fC$}\rangle_A.
\end{align*}
By \cite[Cor.~3.4]{myind}, this yields
\begin{align*}
\fI_{X_I\fC}&=\langle T_x\bC_u \mid x \in X_I,u\in W_I,u \leq_{\cL,I} 
v \mbox{ for some $v \in \fC$}\rangle_A,\\
\hat{\fI}_{X_I\fC}&=\langle T_x\bC_u\mid x \in X_I,u \in W_I, u\not\in\fC,
u \leq_{\cL,I} v \mbox{ for some $v \in \fC$}\rangle_A.
\end{align*}
Thus, we see that the $\bH$-module $[X_I\fC]_A$ has two $A$-bases: 
firstly, the standard basis $\{c_{xu} \mid x\in X_I,u \in \fC\}$ where 
$c_{xu}$ is the residue class of $\bC_{xu}$ and, secondly, the basis 
$\{f_{xu} \mid x\in X_I,u \in \fC\}$ where $f_{xu}$ denotes the residue 
class of $T_x\bC_u$. The change of basis is given by the equations:
\[ c_{yv}=\sum_{\atop{x \in X_I,u \in \fC}{xu \sqsubseteq yv}}
p_{xu,yv}^* \, f_{xu} \qquad \mbox{for any $y \in X_I$, $v\in \fC$}.\]
Furthermore, recalling the definition of $f_{xu}$, it is obvious that 
the map 
\[\bH \otimes_{\bH_I} [\fC]_A\rightarrow [X_I\fC]_A,\qquad
T_x \otimes c_u \mapsto f_{xu} \quad (x\in X_I,u \in \fC),\]
is an isomorphism of $\bH$-modules, where $\{c_u\mid u\in \fC\}$
is the standard basis of $[\fC]_A$ as in (\ref{leftrep}). 
\end{proof}

\begin{rem} \label{otherind} In the above setting, we also have an
isomorphism of $\bH$-modules
\[ [X_I\fC]_A^\delta\stackrel{\sim}{\rightarrow}\Ind_I^S([\fC]_A^{\delta_I}), 
\qquad c_{yv} \mapsto \sum_{\atop{x \in X_I,u\in \fC}{xu \sqsubseteq yv}}
(-1)^{l(x)} p_{xu,yv}^*\,\big(T_x \otimes c_u\big),\]
where  $\delta_I$ denotes the restriction of $\delta$ to $\bH_I$.
Indeed, applying Remark~\ref{otherC} to the formula in 
Proposition~\ref{mythm} yields 
\begin{align*}
\delta(\bC_{yv})&=j(\bC_{yv})=\sum_{\atop{x\in X_I,u\in W_I}{xu
\sqsubseteq yv}} (-1)^{l(x)}p_{xu,yv}^*\, T_x\,j(\bC_u)\\
&=\sum_{\atop{x\in X_I,u\in W_I}{xu\sqsubseteq yv}} 
(-1)^{l(x)}p_{xu,yv}^*\, T_x\,\delta_I(\bC_u)
\end{align*}
for any $y \in X_I$ and $v \in W_I$. We can now argue as in the above
proof, using the fact that $[\fC]_A^{\delta_I}=\delta_I(\fI_{\fC})/
\delta_I(\hat{\fI}_{\fC})$ and $[X_I\fC]_A^\delta=\delta(\fI_{X_I\fC})/
\delta(\hat{\fI}_{X_I\fC})$.
\end{rem}

Our aim is to show that the relation ``$\approx$'' in Definition~\ref{eqcells} 
behaves well with respect to the induction of cells. We begin with the
following result.

\begin{lem} \label{indep} Assume that $\fC,\fC_1$ are two left cells in 
$W_I$ such that $\fC \approx \fC_1$. Let $\fC\stackrel{\sim}{\rightarrow} 
{\fC}_1$, $u \mapsto u_1$, be a bijection such that the property 
$(\heartsuit)$ in Definition~\ref{eqcells} holds. Then we have 
\[p_{xu,yv}^*=p_{xu_1,yv_1}^* \qquad \mbox{for all $x,y \in X_I$
and all $u,v \in \fC$}.\]
\end{lem}

\begin{proof} First we claim that 
\begin{equation*}
r_{xu,yv}=r_{xu_1,yv_1} \qquad \mbox{for all $x,y \in X_I$
and all $u,v \in \fC$}.\tag{$*$}
\end{equation*}
To see this, consider the expression of $r_{xu,yv}$ in Lemma~\ref{rpol}
and note that the coefficients $\overline{R}_{xw,y}^*$ and $\tilde{p}_{w',w}$
do not depend on $u$ or $v$. Hence our assumption ($\heartsuit$) implies that 
($*$) holds. Now, following once more the inductive procedure for solving 
the system of equations given by (KL1)--(KL3) above, we see that we also have 
$p_{xu,yv}^*=p_{xu_1,yv_1}^*$ for all $x,y \in X_I$ and all $u,v \in C$. 
Just note that, for $u,v\in \fC$, the condition $xu\sqsubset zw 
\sqsubseteq vy$ implies that $u \leq_{\cL,I} w \leq_{\cL,I} v$ and 
so $w \in \fC$.
\end{proof}

\begin{prop} \label{induind0} Let $\fC,\fC_1$ be two left cells in $W_I$
such that $\fC \approx \fC_1$. Then we also have $X_I \fC \approx X_I\fC_1$.
More precisely, let $\fC\stackrel{\sim}{\rightarrow} {\fC}_1$, 
$u \mapsto u_1$, be a bijection satisfying $(\heartsuit)$. Then the 
bijection $X_I\fC\stackrel{\sim}{\rightarrow} X_I{\fC}_1$, $xu 
\mapsto xu_1$, satisfies $(\heartsuit)$. 
\end{prop}

\begin{proof} We have seen in Theorem~\ref{inducell} that there is an
isomorphism of $\bH$-modules
\[\bH \otimes_{\bH_I} [\fC]_A\rightarrow [X_I\fC]_A,
\qquad T_x \otimes c_u \mapsto f_{xu} \quad (x\in X_I,u \in \fC),\]
where $\{c_u\mid u\in \fC\}$ is the standard basis of $[\fC]_A$ 
as in (\ref{leftrep}) and $f_{xu}$ denotes the residue class of $T_x\bC_u$ in 
$[X_I\fC]_A$. The base change is given by the equations
\[ c_{yv}=\sum_{\atop{x \in X_I,u \in \fC}{xu \sqsubseteq yv}}
p_{xu,yv}^*\, f_{xu} \qquad \mbox{for any $y \in X_I$, $v \in \fC$}.\]
Similarly, we have an isomorphism of $\bH$-modules
\[\bH \otimes_{\bH_I} [\fC_1]_A\rightarrow [X_I\fC_1]_A,\qquad
T_x \otimes c_{u_1} \mapsto f_{xu_1} \quad (x\in X_I,u_1 \in \fC_1),\]
where $\{c_{u_1}\mid u_1\in \fC_1\}$ is the standard basis of 
$[\fC_1]_A$ as in (\ref{leftrep}) and $f_{xu_1}$ denotes the 
residue class of $T_x\bC_{u_1}$ in $[X_I\fC_1]_A$. The base 
change is given by the equations
\[ c_{yv_1}=\sum_{\atop{x \in X_I,u_1 \in \fC_1}{xu_1 \sqsubseteq yv_1}}
p_{xu_1,yv_1}^*\, f_{xu_1}\qquad\mbox{for any $y\in X_I$, $v_1\in\fC$}.\]
Now, the fact that $(\heartsuit)$ holds for the bijection $\fC 
\stackrel{\sim}{\rightarrow} \fC_1$ means that any $\bC_s$ (where $s\in I$
is a generator of $W_I$) acts in the same way on the standard bases of 
$[\fC]_A$ and of $[\fC_1]_A$, respectively. Hence, by the 
definition of the induced module (see also the explicit formulas in
\cite[\S 9.1]{ourbuch}), it is clear that any $\bC_s$ (where $s\in S$ is a 
generator of $W$) will act in the same way on the bases $\{T_x \otimes 
c_u\}$ and $\{T_x\otimes c_{u_1}\}$ of $\bH \otimes_{\bH_I} [\fC]_A$ and 
$\bH\otimes_{\bH_I} [\fC_1]_A$, respectively. Then the above two 
isomorphisms show that any $\bC_s$ ($s\in S$) acts in the same way on the 
bases $\{f_{xu}\}$ and $\{f_{xu_1}\}$ of $[X_I\fC]_A$ and of 
$[X_I\fC_1]_A$, respectively. Finally, by Lemma~\ref{indep}, the 
two base changes are performed by using exactly the same coefficients. 
Hence, any $\bC_s$ ($s \in W$) will also act in the same way on the standard 
bases $\{c_{xu}\}$ and $\{c_{xu_1}\}$ of $[X_I\fC]_A$ and of 
$[X_I\fC_1]_A$, respectively. 
\end{proof}

\begin{cor} \label{ind1} In the setting of Proposition~\ref{induind0}, 
assume that the partitions of $X_I\fC$ and $X_I\fC_1$ into left cells 
of $W$ are given by 
\[X_I\fC=\coprod_{\alpha \in \cA} \fC^{(\alpha)} \qquad \mbox{and} 
\qquad X_I{\fC}_1=\coprod_{\beta \in \cB} \fC_1^{(\beta)},\]
respectively, where $\cA$ and $\cB$ are some indexing sets. Then there 
exists a bijection $f\colon \cA \rightarrow \cB$ such that 
$\fC^{(\alpha)} \approx \fC_1^{(f(\alpha))}$ for all $\alpha\in \cA$.
\end{cor}

\begin{proof} We have seen in Proposition~\ref{induind0} that the
bijection $X_I\fC\stackrel{\sim}{\rightarrow} X_I{\fC}_1$, $xu 
\mapsto xu_1$, satisfies $(\heartsuit)$, that is, we have 
\[ h_{s,xu,yv}=h_{s,xu_1,yv_1} \qquad \mbox{for $s\in S$, $x,y\in X_I$
and $u,v\in \fC$}.\]
By the definition of left cells, this immediately implies that the 
bijection $X_I\fC\stackrel{\sim}{\rightarrow} X_I{\fC}_1$ preserves the 
partition of the sets $X_I\fC$ and $X_I\fC_1$ into left cells, and that 
corresponding left cells are related by ``$\approx$''.
\end{proof}

\section{Relative left, right and two-sided cells} \label{sec1}
We preserve the setting of the previous sections, where we consider a 
parabolic subgroup $W_I$. In this section, we pursue the study of the 
relative pre-orders $\leq_{\cL,I}$, $\leq_{\cR,I}$ etc.\ introduced in 
(\ref{relpre}). Our Conjecture~\ref{conja22} predicts that we have an 
analogue of ($\spadesuit$) (see Section~1) in this relative setting. The 
main result of this section shows that the conjecture is true in the equal 
parameter case. This will play an essential role in our proof of property 
($\spadesuit$) for groups of type $B_n$ in the ``asymptotic case''.

\begin{rem} \label{righthand} Recall that $X_I$ is the set of distinguished
left coset representatives of $W_I$ in $W$. Applying the anti-automorphism 
$\flat \colon \bH \rightarrow \bH$ such that $T_w^\flat=T_{w^{-1}}$ for
all $w\in W$, we also obtain ``right-handed'' versions of the results in
Section~\ref{sec-indu}. First of all, the set $Y_I:=X_I^{-1}$ is the set of 
distinguished right coset representatives of $W_I$ in $W$. Thus, we can 
write any $w\in W$ uniquely in the form $w=ux$ where $u\in W_I$, $x\in Y_I$ 
and $l(ux)=l(u)+l(x)$. Since this will play a crucial role in the proof
of Lemma~\ref{relcellA0}, let us explicitly state the analogue of 
Corollary~\ref{cor-rel2}. Let $y\in Y_I$ and $v\in W_I$. 
\begin{itemize}
\item[(a)] $\bC_{vy}$ is a linear combination of $\bC_v\,T_y$ and terms
$\bC_u\,T_x$ where $x\in Y_I$ and $u\in W_I$ are such that $x<y$, 
$u \leq_{\cR,I} v$ and $ux<vy$. More precisely, by Proposition~\ref{mythm},
we have
\[ \bC_{vy}=\sum_{x\in X_I,u\in W_I} a_{ux,vy}\, \bC_uT_x\]
where the coefficients satisfy the following conditions:
\begin{alignat*}{2}
a_{vy,vy} &=1 \qquad &&\mbox{if $ux=vy$}\\
a_{ux,vy} & \in A_{<0}\qquad &&\mbox{if $u\leq_{\cR,I} v$ and 
$x<y$},\\ a_{ux,vy}&=0 \qquad &&\mbox{otherwise}.
\end{alignat*}
(We have $a_{ux,vy}=p_{(ux)^{-1},(vy)^{-1}}^*$ in the notation of
Proposition~\ref{mythm}.)
\item[(b)] Conversely, $\bC_v\,T_y$ is a linear combination of $\bC_{vy}$ and 
terms $\bC_{ux}$ where $x\in X_I$ and $u\in W_I$ are such that $x<y$, $u 
\leq_{\cR,I} v$ and $ux<vy$. More precisely, arguing as in the proof
of Corollary~\ref{cor-rel2}, we have
\[ T_v\bC_y=\sum_{x\in X_I,u\in W_I} b_{ux,vy}\, \bC_{ux}\]
where the coefficients satisfy the following conditions:
\begin{alignat*}{2}
b_{vy,vy} &=1 \qquad &&\mbox{if $ux=vy$}\\
b_{ux,vy} & \in A_{<0} \qquad &&\mbox{if $u\leq_{\cR,I} v$ and 
$x<y$},\\ b_{ux,vy}&=0 \qquad &&\mbox{otherwise}.
\end{alignat*}
\end{itemize}
Using the above relations, we obtain the following formula.
\end{rem}

\begin{lem} \label{right1} Let $u,v,w \in W_I$ and $x,y \in Y_I$. Then we
have 
\[ h_{w,vy,ux}=\sum_{\atop{x_1\in Y_I}{u',u_1\in W_I}} 
a_{u'x_1,vy}\, h_{w,u',u_1}\, b_{ux,u_1x_1}.\]
In the above sum, we can assume that $u \leq_{\cLR} u_1 \leq_{\cLR}
u'\leq_{\cLR} v$ and $x \leq x_1 \leq y$.
\end{lem}

\begin{proof} Using the formulas in Remark~\ref{righthand}, we compute:
\begin{align*}
\bC_w\bC_{vy}&=\sum_{x_1\in Y_I,u'\in W_I} a_{u'x_1,vy} \bC_w\bC_{u'}T_{x_1}\\
&=\sum_{\atop{x_1\in Y_I}{u',u_1\in W_I}} a_{u'x_1,vy}\,h_{w,u',u_1}\,
\bC_{u_1}T_{x_1}\\ &=\sum_{\atop{x,x_1\in Y_I}{u,u',u_1\in W_I}}
a_{u'x_1,vy}\, h_{w,u',u_1}\, b_{ux,u_1x_1}\, \bC_{ux}.
\end{align*}
This yields the above formula. Now let $x_1\in Y_I$ and $u',u_1\in W_I$
be such that the corresponding term in the expression for $h_{w,vy,ux}$ is 
non-zero. Then $a_{u'x_1,vy}\neq 0$ and $b_{ux,u_1x_1} \neq 0$. This implies 
$x\leq x_1\leq y$, $u\leq_{\cR,I} u_1$ and $u' \leq_{\cR,I} v$; see 
the conditions in Remark~\ref{righthand}. Furthermore, if $h_{w,u',u_1} 
\neq 0$, then $u_1 \leq_{\cL,I} u'$. In particular, we have $u 
\leq_{\cLR} u_1 \leq_{\cLR} u'\leq_{\cLR} v$.
\end{proof}

\begin{lem} \label{lema23} Let $u,v\in W_I$ and $y\in Y_I$. Then we have
\[ uy \leq_{\cL,I} vy \qquad \Leftrightarrow \qquad u \leq_{\cL,I} v.\]
\end{lem}

\begin{proof} For the implication ``$\Leftarrow$'', see 
\cite[Prop.~9.11]{Lusztig03}. To prove the implication ``$\Rightarrow$'', 
we may assume without loss of generality that $u\neq v$ and $uy 
\leftarrow_{\cL,I} vy$, that is, we have $h_{s,vy,uy} \neq 0$ for some 
$s\in I$. Then we have $sv>v$, $su<u$ and the formula in (\ref{multrule})
shows that there are two cases:  If $svy=uy$, then $u=sv>v$ and so 
$u\leq_{\cL,I} v$. If $suy<uy<vy<svy$ and $M_{uy,vy}^s \neq 0$, then 
$su<u<v<sv$ and \cite[Lemma~9.10]{Lusztig03} shows that $M_{u,v}^s= 
M_{uy,vy}^s \neq 0$.  Again, we have $u\leq_{\cL,I} v$.
\end{proof}

\begin{prop} \label{propa21} Let $u,v\in W_I$ and $x,y\in Y_I$. Then
we have the following implication:
\[ ux \leq_{\cL,I} vy \quad \Rightarrow \quad u\leq_{\cLR,I} v 
\quad \mbox{and} \quad x \leq y.\]
In particular, if  $ux \sim_{\cL,I} vy$, then we necessarily have $x=y$ and
$u\sim_{\cL,I} v$.
\end{prop}
 
\begin{proof} We may assume without loss of generality that $ux 
\leftarrow_{\cL,I} vy$, that is, $h_{s,vy,ux}\neq 0$ for some $s \in I$. 
Then the assertion follows from Lemma~\ref{right1}.
\end{proof}

\begin{conj}[Relative version of $(\spadesuit)$] \label{conja22} Let 
$u,v\in W_I$ and $x,y\in Y_I$. Then we have the following implication:
\[ ux \leq_{\cL,I} vy \quad \mbox{and} \quad u \sim_{\cLR,I} v \quad 
\Rightarrow \quad u \sim_{\cL,I} v \quad \mbox{and} \quad x =y.\]
Note that $u \sim_{\cLR,I} v$ and $u \sim_{\cL,I} y$ just mean the
usual Kazhdan--Lusztig relations inside $W_I$.
\end{conj}

\begin{rem} \label{conja221} Assume that $I=S$; then $W_I=W$ and  
$Y_I=\{1\}$. In this case, the above conjecture reads:
\[ u \leq_{\cL} v \quad \mbox{and} \quad u \sim_{\cLR} v \quad 
\Rightarrow \quad u \sim_{\cL} v \]
(for any $u,v\in W$).  Thus, Conjecture~\ref{conja22} can be seen as a 
generalization of the implication ($\spadesuit$) stated in the introduction. 
Using computer programs written in the {\sf GAP} programming language, we 
have verified that Conjecture~\ref{conja22} holds for $W$ of type $F_4$,
all choices of $I$ and all choices of integer-valued weight functions on 
$W$ (using the techniques in \cite{my04}). In Theorem~\ref{relcellA} we 
will show that this is also true in the case of equal parameters.
\end{rem}

For the remainder of this section, we assume that $W$ is bounded and 
integral in the sense of \cite[1.11 and 13.2]{Lusztig03}. Furthermore, we
assume that $q_s=q_t$ for all $s,t\in S$ (the ``equal parameter'' case). 
Let $q:=q_s$ ($s\in S$). Then our hypotheses imply that
\[ P_{x,y}^*\in q^{-1}{\N}[q^{-1}] \qquad \mbox{and}\qquad 
h_{x,y,z} \in {\N}[q,q^{-1}]\]
for all $x,y,z\in W$, where $\N=\{0,1,2,\ldots\}$. See Lusztig \cite{Lu1},
\cite[15.1]{Lusztig03} and Springer \cite{Spr}.  We shall need some 
properties of Lusztig's function $\ba_I\colon W_I\rightarrow\N$ defined by 
\[ \ba_I(w)=\min\{n \in \N \mid q^nh_{u,v,w} \in {\Z}[q] \mbox{ for all
$u,v\in W_I$}\}.\]
Note that $h_{u,v,w}=\overline{h}_{u,v,w}$. So, if $\ba_I(w)=n$, then 
$q^nh_{u,v,w} \in {\Z}[q]$ and $q^{-n}h_{u,v,w} \in {\Z}[q^{-1}]$. 
Furthermore, both $q^n$ and $q^{-n}$ occur with non-zero coefficient
in $h_{u,v,w}$. In \cite[Chap.~15]{Lusztig03}, the following three
properties are established:
\begin{itemize}
\item[{\bf (P4)}] The function $\ba_I\colon W_I \rightarrow \N$ is constant 
on two-sided cells. 
\item[{\bf (P8)}] Let $u,v,w\in W_I$ be such that $q^{\ba_I(w)}h_{u,v,w}$ 
has a non-zero constant term. Then $v \sim_{\cL,I} w$, $u \sim_{\cR,I} w$ 
and $u \sim_{\cL,I} v^{-1}$.
\item[{\bf (P9)}] Let $u,v\in W_I$ be such that $u \leq_{\cL,I} v$ and 
$\ba_I(u)=\ba_I(v)$. Then $u \sim_{\cL,I} v$. 
\end{itemize}
(There even is a list of $15$ properties, but we only need the above 
three.) Note that {\bf (P4), (P9)} together imply that ($\spadesuit$) holds 
for $W_I$.

\begin{lem} \label{relcellA0} In the above setting, let $u,v,w \in W_I$ and 
$x,y \in Y_I$. Then the coefficient $h_{w,vy,ux}$ has the following properties.
\begin{itemize}
\item[(a)] If $h_{w,vy,ux} \neq 0$, then $u \leq_{\cLR,I} v$ and $x\leq y$. 
\item[(b)] If $x=y$, then $h_{w,vy,uy}=h_{w,v,u}$.
\item[(c)] Assume that $u \sim_{\cLR,I} v$ and let $n:=\ba_I(u)=\ba_I(v)$;
see {\bf (P4)}. If the coefficient of $q^n$ in $h_{w,vy,ux}$ is non-zero, 
then $x=y$.
\end{itemize}
\end{lem}

\begin{proof} If $h_{w,vy,ux} \neq 0$ then $ux\leq_{\cL,I} vy$ and (a)
follows from Proposition~\ref{propa21}. To prove (b) and (c), we use
the formula in Lemma~\ref{right1}:
\[ h_{w,vy,ux}=\sum_{\atop{x_1\in X_I}{u',u_1\in W_I}} 
a_{u'x_1,vy}\, h_{w,u',u_1}\, b_{ux,u_1x_1},\]
where the sum runs over all $x_1,u',u_1$ such that 
\begin{gather*}
x\leq x_1\leq y, \tag{$*_1$}\\
u \leq_{\cLR,I}u_1\leq_{\cLR,I} u'\leq_{\cLR,I} v.\tag{$*_2$}
\end{gather*}
Now, if $x_1=x$, then $b_{ux,u_1x}=0$ unless $u=u_1$ (in which case the 
result is~$1$; see the conditions in Remark~\ref{righthand}). Similarly, 
if $x_1=y$, then $a_{u'y,vy}=0$ unless $u'=v$ (in which case the 
result is~$1$).  Hence, if $x=y$, the above sum reduces to
\[ h_{w,vy,ux}= a_{vy,vy}\, h_{w,v,u}\, b_{uy,uy}= h_{w,v,u}.\]
Thus, (b) is proved. Finally, to prove (c), assume that $x<y$ and
that the coefficient of $q^n$ in $h_{w,vy,ux}$ is non-zero, where 
$n=\ba_I(u)=\ba_I(v)$. We must show that $u,v$ cannot be in the same
two-sided cell. Splitting the above sum into three pieces according to 
$x_1=x$, $x_1=y$ and $x<x_1<y$, we obtain
\begin{align*}
 h_{w,vy,ux}&=\sum_{u'\in W_I} a_{u'x,vy}\, h_{w,u',u} +
\sum_{u_1\in W_I} h_{w,v,u_1}\, b_{ux,u_1y} \\&+ 
\sum_{u',u_1\in W_I} \Big(\sum_{\atop{x_1\in X_I}{x<x_1<y}} 
a_{u'x_1,vy}\, b_{ux,u_1x_1}\Big)\, h_{w,u',u_1}. 
\end{align*}
Note that, since $x<x_1<y$, all the coefficients $a_{u'x,vy}$, $b_{ux,u_1y}$, 
$a_{u'x_1,vy}$ and $b_{ux,u_1x_1}$ occuring in the above expression
lie in $q^{-1}{\Z}[q^{-1}]$; see once more the conditions in
Remark~\ref{righthand} and recall that $q=q_s$ (for all $s\in S$). Hence 
we can re-write the above expression as follows:
\[ h_{w,vy,ux}=\sum_{u_1,u' \in W_I} f_{u_1,u'} \, h_{w,u',u_1} 
\qquad \mbox{where $f_{u_1,u'} \in q^{-1}{\Z}[q^{-1}]$},\]
where we can assume that ($*_2$) holds. 

Now, we are assuming that the coefficient of $q^n$ in $h_{w,vy,ux}$ is 
non-zero. So there exist some $u',u_1\in W_I$ such that the coefficient 
of $q^n$ in $f_{u_1,u'} \, h_{w,u',u_1}$ is non-zero. Since $f_{u_1,u'} \in 
q^{-1}{\Z}[q^{-1}]$, we deduce that there exists some $m>n$ such $q^m$ 
has a non-zero coefficient in $h_{w,u',u_1}$. By the definition of
the $\ba$-function, this means that $\ba_I(u_1)\geq m>n$. Now, if we had
$u \sim_{\cLR,I} v$, then ($*_2$) would imply $u \sim_{\cLR}u_1\sim_{\cLR} u'
\sim_{\cLR} v$, yielding the contradiction 
\[ \ba(u_1)=\ba(u')=\ba(u)=\ba(v)=n; \qquad \mbox{see {\bf (P4)}}.\]
Consequently, $u$ and $v$ cannot lie in the same two-sided cell.
\end{proof}

\begin{thm} \label{relcellA} Assume that $W$ is bounded, integral in 
the sense of \cite{Lusztig03} and that $q_s=q_t$ for all $s,t\in S$. Then 
Conjecture~\ref{conja22} holds for all parabolic subgroups $W_I 
\subseteq W$.
\end{thm}

\begin{proof} Let us fix a subset $I \subseteq S$. Let $u,v \in W_I$ and 
$x,y \in Y_I$ be such that $ux \leq_{\cL,I} vy$ and $u \sim_{\cLR,I} v$.
We want to show that $x=y$ and $u \sim_{\cL,I} v$. Suppose we already know 
that $x=y$. Then, since $uy \leq_{\cL,I} vy$, we can apply Lemma~\ref{lema23}
and this yields $u \leq_{\cL,I} v$. Thus, we have $u \leq_{\cL,I} v$ and
$u \sim_{\cLR,I} v$.  So {\bf (P4), (P9)} imply that $u \sim_{\cL,I} v$, as
desired. Hence, it is sufficient to prove that $x=y$. First of all, using 
Proposition~\ref{propa21}, we may assume without loss of generality that 
$ux \neq vy$ and $ux\leftarrow_{\cL,I} vy$, that is, $\bC_{ux}$ occurs in 
$\bC_s\bC_{vy}$ for some $s \in I$ such that $svy>vy$. Since $s \in I$, this 
implies $sv>v$, and the multiplication rule for the Kazhdan--Lusztig basis 
(see Section~\ref{sec0a}) shows that we must have $su<u$ and $u\neq v$. We 
shall now try to imitate the proof of ($\tilde{P}$) in \cite[15.5]{Lusztig03}.

Since $u \sim_{\cLR,I} v$, we have $n:=\ba_I(u)=\ba_I(v)$ by {\bf (P4)}. For 
any Laurent polynomial $f \in {\Z}[q,q^{-1}]$, we denote by $\pi_n(f)$ the 
coefficient of $q^n$ in $f$, where we write $q:=q_s$ ($s\in S$) as above. Now 
we argue as follows. By the definition of the $\ba$-function, there exist 
some $w,v'\in W_I$ such that $q^nh_{w,v',v}$ has a non-zero constant term. 
Since $h_{w,v',v}= \overline{h}_{w,v',v}$, this means that the coefficient
of $q^n$ in $h_{w,v',v}$ is non-zero. Thus, using {\bf (P8)}, we have  
\begin{equation*}
\pi_n(h_{w,v',v})\neq 0 \qquad \mbox{and}\qquad v' \sim_{\cL,I} v.\tag{1}
\end{equation*}
We can express the product $\bC_s(\bC_w\bC_{v'y})$ as a linear combination of 
terms $\bC_{wz}$ where $w\in W_I$ and $z\in Y_I$. Denote by $\kappa_{wz}$ the 
coefficient of $\bC_{wz}$ in that product. We have 
\[ \kappa_{wz}=\sum_{w_1 \in W_I,z_1\in Y_I} h_{w,v'y,w_1z_1}\, 
h_{s,w_1z_1,wz}.\]
In particular,
\begin{align*}
\kappa_{ux}&=\sum_{w_1\in W_I,z_1\in Y_I} h_{w,v'y,w_1z_1}\,h_{s,w_1z_1,ux}\\
&= h_{w,v'y,vy}\, h_{s,vy,ux} +\sum_{\atop{w_1 \in W_I,z_1\in Y_I}{w_1z_1 
\neq vy}} h_{w,v'y,w_1z_1}\, h_{s,w_1z_1,ux}.
\end{align*}
Since $svy>vy$,  the multiplication rule for the Kazhdan--Lusztig 
basis shows that $h_{s,vy,ux}$ equals $1$ or $M_{ux,vy}^s$, and the latter
is an integer by \cite[6.5]{Lusztig03}. Hence we have $h_{s,vy,ux}\in \Z$ 
in both cases and so 
\[ \pi_n(h_{w,v'y,vy}\,h_{s,vy,ux})= \pi_n(h_{w,v'y,vy})\, h_{s,vy,ux}=
\pi_n(h_{w,v',v})\, h_{s,vy,ux},\]
where the last equality holds by Lemma~\ref{relcellA0}(b). We are assuming
that $h_{s,vy,ux}\neq 0$. In combination with (1) and the above identity,
we conclude that 
\begin{equation*}
\pi_n(h_{w,v'y,vy}\,h_{s,vy,ux})=\pi_n(h_{w,v',v})\, h_{s,vy,ux} 
\neq 0. \tag{2} 
\end{equation*}
Since all polynomials involved in the expression for $\kappa_{ux}$ have 
non-negative coefficients (thanks to the assumption that $W$ is integral), 
the non-zero coefficient of $q^n$ arising from (2) will not cancel out 
with the coefficients of $q^n$ from the remaining terms in $\kappa_{ux}$. So 
we can conclude, as in the proof of Lusztig \cite[15.5]{Lusztig03}, that 
\[ \pi_n(\kappa_{ux}) \neq 0.\]
On the other hand, since $\bC_s(\bC_w\bC_{v'y})=(\bC_s\bC_w)\bC_{v'y}$, we 
also have the following expression for $\kappa_{ux}$:
\[ \kappa_{ux}=\sum_{w' \in W_I} h_{s,w,w'} h_{w',v'y,ux}.\]
Since $\pi_n(\kappa_{ux})\neq 0$, there exists some $w' \in W_I$ such that 
\begin{equation*}
 \pi_n(h_{s,w,w'}\,h_{w',v'y,ux})\neq 0.\tag{3}
\end{equation*}
By (1), we have $h_{w,v',v} \neq 0$ and so $v \leq_{\cR,I} w$. Hence the 
left descent set of $w$ is contained in the left descent set of~$v$; see 
\cite[8.6]{Lusztig03}. So, since $sv>v$, we also have $sw>w$. Then the 
multiplication rule for the Kazhdan--Lusztig basis and \cite[6.5]{Lusztig03}
show that $h_{s,w,w'}\in \Z$. Hence (3) implies that 
\[\pi_n(h_{w',v'y,ux})\neq 0 \qquad \mbox{where} \qquad w' \in W_I.\] 
By (1), we also have $v' \sim_{\cLR,I} v \sim_{\cLR,I} u$. Hence
Lemma~\ref{relcellA0}(c) yields $x=y$, as desired. 
\end{proof}

\begin{exmp} \label{conjA} Let $W=\fS_n$ be the symmetric group. Then
Conjecture~\ref{conja22} holds for all parabolic subgroups $W_I\subseteq W$.

Indeed, $\fS_n$ is finite, hence bounded. Since the product of any two 
generators has order $2$ or $3$, the group is integral.  Furthermore, since 
all generators are conjugate, all the parameters are equal. Hence the 
hypotheses of Theorem~\ref{relcellA} are satisfied. 
\end{exmp}

\section{On the left pre-order $\leq_{\cL}$ in type $B_n$} \label{sec-cellu}
In this and the subsequent sections, we let $W=W_n$ be a Coxeter group of 
type $B_n$ ($n \geq 2$). We assume that the generators, relations and the 
weight function $L \colon W_n \rightarrow \Gamma$ are given by the 
following diagram:
\begin{center}
\begin{picture}(250,50)
\put(  3, 25){$B_n$}
\put(  0, 05){$\{q_s\}$:}
\put( 40, 25){\circle{10}}
\put( 44, 22){\line(1,0){33}}
\put( 44, 28){\line(1,0){33}}
\put( 81, 25){\circle{10}}
\put( 86, 25){\line(1,0){29}}
\put(120, 25){\circle{10}}
\put(125, 25){\line(1,0){20}}
\put(155, 22){$\cdot$}
\put(165, 22){$\cdot$}
\put(175, 22){$\cdot$}
\put(185, 25){\line(1,0){20}}
\put(210, 25){\circle{10}}
\put( 37, 37){$t$}
\put( 36, 05){$e^b$}
\put( 76, 37){$s_1$}
\put( 78, 05){$e^a$}
\put(116, 37){$s_2$}
\put(118, 05){$e^a$}
\put(203, 37){$s_{n-1}$}
\put(208, 05){$e^a$}
\end{picture}
\end{center}
where $a,b \in \Gamma$ are such that $a>0$ and $b>0$. Let $\bH_n$ be the 
corresponding Iwahori--Hecke algebra over $A={\Z}[\Gamma]$, where we set 
\[ Q:=q_t=e^b \qquad \mbox{and}\qquad q:=q_{s_1}=\cdots =q_{s_{n-1}}=e^a.\]
Let $K$ be the field of fractions of $A$ and set $\bH_{n,K}=K \otimes_A 
\bH_n$. Throughout this and the subsequent sections, we assume that 
$b/a$ is ``large'' with respect to $n$, more precisely:
\begin{center} \fbox{$b > (n-1) a$}\end{center}
(Here, $(n-1)a$ means $a+\cdots +a$ in $\Gamma$, with $n-1$ summands.) We 
refer to this hypothesis as the {\bf ``asymptotic case''} in type $B_n$. 

The main results of this section are:
\begin{itemize}
\item Theorem~\ref{strongL}, which gives a strengthening of the results of
Bonnaf\'e--Iancu \cite{BI} concerning the left cells of $W_n$ (and, as a
bi-product, also yields a new proof of Bonnaf\'e's result \cite{BI2} 
on the two-sided cells);
\item Theorem~\ref{relsn}, which shows that ($\spadesuit$) holds in $W_n$.
\end{itemize}

\begin{rem} \label{asymB} Let us consider the abelian group $\Gamma^\circ=
\Z^2$ and let $\leq$ be the usual lexicographic order on $\Gamma^\circ$. 
Thus, we have $(i,j)<(i',j')$ if $i<i'$ or if $i=i'$ and $j<j'$. Let $L^\circ
\colon W_n \rightarrow {\Z}^2$ be the weight function such that 
\[ L(t)=(1,0)\qquad \mbox{and}\qquad L(s_1)= \cdots = L(s_{n-1})=(0,1).\]
Then $A^\circ={\Z}[\Gamma^\circ]$ is nothing but the ring of Laurent 
polynomials in two independent indeterminates $V=e^{(1,0)}$ and 
$v=e^{(0,1)}$.  This is the ``asymptotic case'' originally 
considered by Bonnaf\'e--Iancu \cite{BI}. We may refer to this case as 
the {\bf ``generic asymptotic case''} in type $B_n$. Let us denote the 
corresponding Iwahori--Hecke algebra by $\bH_n^\circ$; let $\{\bC_w^\circ
\mid w\in W_n\}$ be the Kazhdan--Lusztig basis of $\bH_n^\circ$ and write
\[ \bC_x^\circ \, \bC_y^\circ=\sum_{z\in W_n} h_{x,y,z}^\circ\,\bC_z^\circ
\quad \mbox{where} \quad h_{x,y,z}^\circ\in A^\circ={\Z}[V^{\pm 1},
v^{\pm 1}].\] 
Now, given an abelian group $\Gamma$ as above and two elements $a,b>0$,
we have a unique ring homomorphism 
\[\theta\colon A^\circ\rightarrow A,\qquad V^iv^j\mapsto e^{ib+ja}.\]
Bonnaf\'e \cite[\S 5]{BI2} has shown that, if $b>(n-1)a$, then the 
Kazhdan--Lusztig basis of $\bH_n$ (with respect to $L\colon W_n 
\rightarrow \Gamma$) is obtained by ``specialisation'' from the 
Kazhdan--Lusztig of $\bH_n^\circ$ and that we have 
\begin{equation*}
h_{x,y,z}=\theta(h_{x,y,z}^\circ) \qquad \mbox{for all $x,y,z\in W_n$}.
\tag{a}
\end{equation*}
In particular, denoting by $\leq_{\cL}^\circ$, $\sim_{\cL}^\circ$, 
$\leq_{\cR}^\circ$, $\sim_{\cR}^\circ$, $\leq_{\cLR}^\circ$, 
$\sim_{\cLR}^\circ$ the pre-order relations on $W_n$ with respect to 
$L^\circ$, we have the implications:
\begin{equation*}
x \leq_{\cL} y \Rightarrow x \leq_{\cL}^\circ y, \qquad
 x \leq_{\cR} y \Rightarrow x \leq_{\cR}^\circ y, \qquad
 x \leq_{\cLR} y \Rightarrow x \leq_{\cLR}^\circ y.\tag{b}
\end{equation*}
These results show that it is usually sufficient to prove identities 
concerning the Kazhdan--Lusztig basis in the ``generic asymptotic case''; 
the analogous identity in the general ``asymptotic case'' then follows by 
specialisation, assuming that $b>(n-1)a$. (In this and the following 
sections, we make an explicit remark at places where we use this kind of 
argument.)
\end{rem}

We shall need some notation from \cite{BI}. Given $w\in W_n$, we denote
by $l_t(w)$ the number of occurences of the generator $t$ in a reduced
expression for $w$, and call this the ``$t$-length''  of $w$.

The parabolic subgroup $\fS_n:=\langle s_1,\ldots,s_{n-1}\rangle$ is
naturally isomorphic to the symmetric group on $\{1,\ldots,n\}$, where
$s_i$ corresponds to the basic transposition $(i,i+1)$. Let $1 \leq l 
\leq n-1$. Then we set $\Sigma_{l,n-l}:=\{s_1,\ldots,s_{n-1}\}\setminus 
\{s_l\}$.  For $l=0$ or $l=n$, we also set $\Sigma_n:=\Sigma_{0,n}=
\Sigma_{n,0}=\{s_1,\ldots,s_{n-1}\}$. Let $X_{l,n-l}$ be the set of 
distinguished left coset representatives of the Young subgroup $\fS_{l,n-l}
:=\langle \Sigma_{l,n-l}\rangle$ in $\fS_n$. We have the parabolic subalgebra 
$\bH_{l,n-l}=\langle T_\sigma \mid \sigma \in \fS_{l,n-l}\rangle_A\subseteq
\bH_n$. Given $x,y\in W_n$, we write
\[ x\leq_{\cL,l}y\qquad\stackrel{\text{def}}{\Longleftrightarrow}
\qquad x\leq_{\cL,\Sigma_{l,n-l}} y\qquad 
\mbox{(see Section~\ref{sec-indu})}.\]
Furthermore, as in \cite[\S 4]{BI}, we set $a_0=1$ and 
\[ a_l:=t(s_1t)(s_2s_1t) \cdots (s_{l-1}s_{l-2} \cdots s_1t) \qquad
\mbox{for $l>0$}.\]
Then, by \cite[Prop.~4.4]{BI}, the set $X_{l,n-l}a_l$ is precisely
the set of distinguished left coset representatives of $\fS_n$ in $W_n$
whose $t$-length equals $l$. Furthermore, every element $w\in W_n$ has
a unique decomposition
\[w=a_wa_l\sigma_w b_w^{-1} \qquad \mbox{where $l=l_t(w)$, $\sigma_w \in 
\fS_{l,n-l}$ and $a_w,b_w\in X_{l,n-1}$};\]
see \cite[4.6]{BI}.  On a combinatorial level, Bonnaf\'e and Iancu 
\cite[\S 3]{BI} define a generalized Robinson--Schensted correspondence which 
associates with each element $w\in W_n$ a pair of $n$-standard bi-tableaux 
$(A(w), B(w))$ such that $A(w)$ and $B(w)$ have the same shape. Here, a 
standard $n$-bitableau is a pair of standard tableaux with a total number
of $n$ boxes (filled with the numbers $1,\ldots,n$), and the shape of such 
a bitableau is a pair of partitions $\lambda=(\lambda_1,\lambda_2)$ such 
that $n=|\lambda_1|+ |\lambda_2|$. With this notation, we have the 
following result.  

\begin{thm}[Bonnaf\'e--Iancu \protect{\cite{BI}} and Bonnaf\'e
\protect{\cite[\S 5]{BI2}}] \label{mainbi} In the above setting, let 
$x,y\in W_n$. Then the following conditions are equivalent:
\begin{itemize}
\item[($\text{a}_1$)] $x\sim_{\cL} y$;
\item[($\text{a}_2$)] $x\sim_{\cL}^\circ y$ (see Remark~\ref{asymB});
\item[(b)] $l:=l_t(x)=l_t(y)$, $b_x=b_y$ and $\sigma_x\sim_{\cL,l}\sigma_y$;
\item[(c)] $B(x)=B(y)$.
\end{itemize}
\end{thm}

(This is the first example where the discussion in Remark~\ref{asymB} 
applies: the equivalences between ($\text{a}_2$), (b) and (c) are proved
in \cite[Theorem~7.7]{BI}; the equivalence between ($\text{a}_1$) and  
($\text{a}_2$) is proved in \cite[Cor.~5.2]{BI2}.)

Note that the equivalence ``($\text{a}_1$) $\Leftrightarrow$ (c)'' is in 
complete formal analogy to the situation in the symmetric group $\fS_n$; see 
Example~\ref{expA}(a). 

Let $\Lambda_n$ be set of all pairs of partitions of total size $n$. We set
\[ \fR_\lambda:=\{w\in W_n\mid\mbox{ $A(w)$, $B(w)$ have shape 
$\lambda$}\} \qquad \mbox{for $\lambda \in \Lambda_n$}.\]
Thus, we have a partition $W_n=\coprod_{\lambda \in \Lambda_n} \fR_\lambda$.
The above result and the properties of the generalized Robinson--Schensted 
correspondence in \cite[\S 3]{BI} immediately imply the following statement:

\begin{cor}[Bonnaf\'e--Iancu \protect{\cite{BI}}] \label{unique} In the 
above setting, let $\lambda \in \Lambda_n$ and denote by $\fT_\lambda$ the 
set of $n$-standard bitableaux of shape $\lambda$. Then the generalized 
Robinson--Schensted correspondence defines a bijection 
\[w_\lambda \colon \fT_\lambda \times \fT_\lambda
\stackrel{\sim}{\rightarrow} \fR_\lambda, \qquad (T,T')\mapsto 
w_\lambda(T,T'),\]
with the following property:
\begin{itemize}
\item[(a)] For a fixed $T'$, the elements $\{w_\lambda(T,T')\mid 
T \in \fT_\lambda\}$ form a left cell.
\item[(b)] For a fixed $T$, the elements $\{w_\lambda(T,T')\mid 
T' \in \fT_\lambda\}$ form a right cell.
\item[(c)] We have $w_\lambda(T,T')^{-1}=w_\lambda(T',T)$ for 
all $T,T'\in \fT_\lambda$.
\end{itemize}
In particular, any left cell contained in $\fR_\lambda$ meets any right
cell contained in $\fR_\lambda$ in exactly one element. Furthermore,
every left cell contains a unique element of the set $\cD_n:=\{z \in W_n\mid 
z^2=1\}$.
\end{cor}

In order  to prove the main results of this section, we need a number
of preliminary steps. We shall frequently use the following result.

\begin{prop}[Bonnaf\'e--Iancu \protect{\cite[Cor.~6.7]{BI}} and Bonnaf\'e
\protect{\cite[\S 5]{BI2}}] \label{maintl} In the above setting, let 
$x,y\in W_n$ be such that $x\leq_{\cLR} y$. Then $l_t(y)\leq l_t(x)$.
In particular, if $x \sim_{\cLR} y$, then $l_t(x)=l_t(y)$.
\end{prop}

(The above result was first proved in \cite{BI} for the
weight function $L^\circ\colon W_n \rightarrow {\Z}^2$; then 
Remark~\ref{asymB}(b) immediately yields the analogous statement in the 
general ``asymptotic case''.) The following two results give some information 
about certain elements of the Kazhdan--Lusztig basis of $\bH_n$.

\begin{lem}[Bonnaf\'e \protect{\cite[\S 2]{BI2}}] \label{lem-c1} For any 
$\sigma \in \fS_n$ and any $0\leq l\leq n$, we have 
\[ \bC_{\sigma}\,\bC_{a_l}=\bC_{\sigma a_l} \qquad \mbox{and}\qquad
\bC_{a_l}\,\bC_{\sigma}=\bC_{a_l \sigma}.\]
Furthermore, if $\sigma \in \fS_{l,n-l}$, then 
\[ \bC_{\sigma}\,\bC_{a_l}=\bC_{\sigma a_l}=\bC_{a_l}\bC_{a_l\sigma a_l} 
\qquad \mbox{where $a_l\sigma a_l\in \fS_{l,n-l}$}.\]
\end{lem}

\begin{proof} By Remark~\ref{asymB}, it is sufficient to prove the
equality $\bC_{\sigma}\,\bC_{a_l}=\bC_{\sigma a_l}$ (for $\sigma \in \fS_n$) in 
the original setting of \cite{BI} where we consider the weight function 
$L^\circ\colon W_n \rightarrow {\Z}^2$. In this case, the statement is proved 
in \cite[Prop.~2.3]{BI2}. The equality $\bC_{a_l}\,\bC_{\sigma}=
\bC_{a_l \sigma}$ is proved similarly. 

Finally, since $a_l=a_l^{-1}$ stabilizes $\Sigma_{l,n-l}$, we have $a_l
\sigma a_l \in \fS_{l,n-l}$ for any $\sigma \in \fS_{l,n-l}$, which yields 
the second statement.
\end{proof}

The following result plays an essential role in the proof of 
Lemma~\ref{lem-c4}. 

\begin{lem} \label{lem-c2} For any $0\leq l \leq n-1$, we have 
\[ T_{ts_1\cdots s_{l}}\,\bC_{a_l}=\bC_{a_{l+1}}+h(a_l)\,\bC_{a_l}\]
where $h(a_l)\in \bH_n$ is an $A$-linear combination of basis elements 
$T_w$ with $w\leq s_1s_2\cdots s_l$. (For $l=0$, we have $a_0=1$, $a_1=t$
and $h(a_0)=-Q^{-1}T_1$.) 
\end{lem}

\begin{proof} Following Dipper--James \cite[3.2]{DiJa92}, we define
\[ u_k^+=(T_{t_1}+Q^{-1}T_1)(T_{t_2}+Q^{-1}T_1)\cdots (T_{t_k}+Q^{-1}T_1)\]
for any $1\leq k\leq n$, where $t_1=t$ and $t_{i+1}=s_it_is_i$ for 
$i\geq 1$. The factors in the definition of $u_k^+$ commute with each other 
and we have
\[ u_k^+\,T_{s_i}= T_{s_i}\,u_k^+ \quad \mbox{for $1 \leq i \leq k-1$};\]
see \cite[\S 3]{DiJa92}. By Bonnaf\'e \cite[Prop.~2.5]{BI2}, we have 
\[ \bC_{a_k}=u_k^+\,T_{\sigma_k}^{-1}=T_{\sigma_k}^{-1}\,u_k^+,\]
where $\sigma_k$ is the longest element in $\fS_k$. (Again, this is 
first proved in the ``generic asymptotic case''; the general case follows 
from the argument in Remark~\ref{asymB}.) Now let $k=l+1$ and note that
\[ T_{\sigma_{l+1}}=T_{\sigma_l} T_{s_l\cdots s_2s_1} \qquad \mbox{and}
\qquad u_{l+1}^+=(T_{t_{l+1}}+Q^{-1}T_1)\,u_l^+.\]
Since $T_{t_{l+1}}$ commutes with $T_{s_i}$ for  $1\leq i\leq l$, we conclude
that 
\begin{align*}
\bC_{a_{l+1}}& =T_{s_l\cdots s_2s_1}^{-1}\,T_{\sigma_l}^{-1}\,
(T_{t_{l+1}}+Q^{-1}T_1)\,u_l^+\\
& =T_{s_l\cdots s_2s_1}^{-1}\, (T_{t_{l+1}}+Q^{-1}T_1)\,T_{\sigma_l}^{-1}
\, u_l^+\\ & =T_{s_l\cdots s_2s_1}^{-1}\, (T_{t_{l+1}}+Q^{-1}T_1)\,\bC_{a_l}
\end{align*}
and so $T_{ts_1s_2\cdots s_l}\bC_{a_l}=\bC_{a_{l+1}}-
Q^{-1}T_{s_l\cdots s_1}^{-1} \bC_{a_l}$, as required.
\end{proof}

The following definitions are inspired by Bonnaf\'e's construction
in \cite[\S 3]{BI2}. Let $w\in W_n$ and write $w=a_wa_l \sigma_w b_w^{-1}$ 
as usual, where $l:=l_t(w)$. We set 
\[E_w:=T_{a_w}\,\bC_{a_l}\,\bC_{\sigma_w b_w^{-1}}=T_{a_w}\, \bC_{a_l\sigma_w
b_w^{-1}},\]
where the second equality holds by Lemma~\ref{lem-c1}. One easily shows 
that the elements $\{E_w \mid w\in W_n\}$ form a basis of $\bH_n$. We will
be interested in the base change from the Kazhdan--Lusztig basis to this
new basis. 

For $y,w\in W_n$, we write $y \preceq w$ if the following conditions 
are satisfied:
\begin{itemize}
\item[(1)] $l:=l_t(y)=l_t(w)$,
\item[(2)] $\sigma_y b_y^{-1} \leq_{\cL,l} \sigma_w b_w^{-1}$, and
\item[(3)] $l(y)<l(w)$ or $y=w$.
\end{itemize}
We write $y\prec w$ if $y \preceq w$ and $y\neq w$. 
Since $\{E_w\}$ is a basis of $\bH_n$, we can write uniquely
\[ \overline{E}_w=T_{a_w^{-1}}^{-1}\, \bC_{a_l}\, \bC_{\sigma_w b_w^{-1}}=
\sum_{y \in W_n} \overline{\lambda}_{y,w}\, E_y \qquad \mbox{where 
$\lambda_{y,w} \in A$}.\]

\begin{lem} \label{qpol} We have $\lambda_{w,w}=1$ and $\lambda_{y,w}=0$ 
unless $y \preceq w$. Furthermore, we have $\lambda_{y,w}\in {\Z}[q,q^{-1}]$.
\end{lem}

\begin{proof} We argue as in the proof of Lemma~\ref{rpol}. Let $w\in W_n$ 
and $l=l_t(w)$. We have 
\[{T}_{a_w^{-1}}^{-1}=\sum_{z\in W_n} \overline{R}_{z,a_w}^*\,T_z\]
where $R_{z,a_w}^*\in A$ are the ``absolute'' $R$-polynomials defined in
\cite[\S 1]{Lusztig83}. We have $R_{a_w,a_w}^*=1$ and $R_{z,a_w}^*=0$
unless $z \leq a_w$. Since $a_w\in \fS_n$, we have $R_{z,a_w}^*\in 
{\Z}[q,q^{-1}]$.

Now let $z\in W_n$ be such that $T_z$ occurs in the above expression. 
Then we can write $z=c\sigma$ where $c\in X_{l,n-l}$ and $\sigma\in 
\fS_{l,n-l}$. Since $l(c\sigma)=l(c)+ l(\sigma)$, we have $T_z=T_c\,
T_\sigma$ and so
\[ \overline{E}_w=\sum_{c,\sigma} \overline{R}_{c\sigma,a_w}^* \, 
T_c\, T_\sigma \, \bC_{a_l} \, \bC_{\sigma_w b_w^{-1}},\]
where the sum runs over all $c\in X_{l,n-l}$ and $\sigma \in \fS_{l,n-l}$.
Now we can also write $T_\sigma=\sum_{\sigma'} \tilde{p}_{\sigma',\sigma} 
\bC_{\sigma'}$ where $\tilde{p}_{\sigma',\sigma} \in {\Z}[q,q^{-1}]$ and 
the sum runs over all $\sigma'\in \fS_{l,n-l}$. Note that
$\tilde{p}_{\sigma,\sigma}=1$ and $\tilde{p}_{\sigma',\sigma}=0$ unless 
$\sigma' \leq \sigma$. Thus, we have
\[\overline{E}_w =\sum_{c,\sigma,\sigma'} \overline{R}_{c\sigma,a_w}^*
\tilde{p}_{\sigma',\sigma}\, T_c\,\bC_{\sigma'} \,\bC_{a_l} \,\bC_{\sigma_w 
b_w^{-1}},\]
where the sum runs over all $c\in X_{l,n-l}$ and all $\sigma,\sigma'\in
\fS_{l,n-l}$.  Now Lemma~\ref{lem-c1} shows that 
\[\bC_{\sigma'}\,\bC_{a_l}\,\bC_{\sigma_w b_w^{-1}}=\bC_{a_l}\,
\bC_{a_l\sigma' a_l}\, \bC_{\sigma_w b_w^{-1}}.\]
Since $a_l\sigma' a_l \in\fS_{l,n-l}$, we can write
\[ \bC_{a_l\sigma' a_l}\bC_{\sigma_w b_w^{-1}}=\sum_{\sigma''\in \fS_{l,n-l}} 
h_{a_l\sigma'a_l, \sigma_wb_w^{-1},\sigma''}\, \bC_{\sigma''},\]
where $h_{a_l\sigma'a_l, \sigma_wb_w^{-1},\sigma''} \in {\Z}[q,q^{-1}]$.
So we conclude that 
\[ \overline{E}_w=\sum_{c,\sigma,\sigma',\sigma''}
\overline{R}_{c\sigma,a_w}^* \tilde{p}_{\sigma',\sigma}h_{a_l\sigma'a_l,
\sigma_w b_w^{-1},\sigma''} \, T_c \, \bC_{a_l}\,\bC_{\sigma''},\]
where the sum runs over all $c\in X_{l,n-l}$ and all $\sigma,\sigma',\sigma''
\in \fS_{l,n-l}$. Now every term $T_c\bC_{a_l}\bC_{\sigma''}$ in the above 
sum is of the form $E_y$ for a unique $y\in W_n$ where $l=l_t(y)$, $a_y=c$, 
$\sigma_y b_y^{-1}=\sigma''$. So we can re-write the above expression as 
\[ \overline{E}_w= \sum_{\atop{y \in W_n}{l_t(y)=l}}
\overline{\lambda}_{y,w}\,E_y\]
where 
\[ \lambda_{y,w}=\sum_{\sigma,\sigma'\in\fS_{l,n-l}} R_{a_y\sigma,a_w}^*\,
\overline{\tilde{p}}_{\sigma', \sigma}\, h_{a_l\sigma'a_l,\sigma_w 
b_w^{-1},\sigma_yb_y^{-1}} \in {\Z}[q,q^{-1}].\]
Assume that $\lambda_{y,w}\neq 0$. We must show that $y \preceq w$. First of
all, we certainly have $l=l_t(w)=l_t(y)$. Furthermore, there exist 
$\sigma,\sigma'\in \fS_{l,n-l}$ such that 
\[ R_{a_y\sigma,a_w}^*\neq 0, \qquad \tilde{p}_{\sigma',\sigma}\neq 0,\qquad
h_{a_l\sigma'a_l,\sigma_wb_w^{-1}, \sigma_y b_y^{-1}}\neq 0.\]
The first condition implies $a_y\sigma\leq a_w$ and so $l(a_y\sigma)\leq
l(a_w)$. The second condition implies $l(\sigma')\leq l(\sigma)$, while the
third condition implies that $\sigma_y b_y^{-1} \leq_{\cL,l} \sigma_w
b_w^{-1}$ and $l(\sigma_y b_y^{-1})\leq l(\sigma')+l(\sigma_wb_w^{-1})$.
(See (\ref{multrule}) and note that $l(a_l\sigma'a_l)=l(\sigma')$.) Hence 
we also have $l(y) \leq l(w)$. Altogther, this means that $y \preceq w$.
Finally, if $y=w$, it is readily checked that $\lambda_{w,w}=1$. 
\end{proof}

The above result shows that, for any $w\in W_n$, we have 
\[ \overline{E}_w=E_w+\sum_{\atop{y \in W_n}{y\prec w}} 
\overline{\lambda}_{y,w}\, E_y \qquad \mbox{where $\lambda_{y,w} \in 
{\Z}[q,q^{-1}]$}.\]
We can now use exactly the same arguments as in the proofs of
Lemma~3.2 and Proposition~3.3 in \cite{myind} (which themselves are an
adaptation of the proof of Lusztig \cite[Prop.~2]{Lusztig83}) to conclude 
that 
\[\bC_w=E_w+\sum_{\atop{y \in W_n}{y\prec w}} \pi_{y,w}\, E_y \]
where $\pi_{y,w} \in q^{-1}{\Z}[q^{-1}]$ for any $y \prec w$. Indeed, 
the family of elements
\[ \{\pi_{y,w} \mid y,w\in W_n,y \preceq w\}\]
is  uniquely determined by the following three conditions:
\begin{align*}
\pi_{w,w} &= 1, \tag{KL1'}\\
\pi_{y,w} & \in A_{<0}\qquad \mbox{if $y\prec w$},\tag{KL2'}\\
\overline{\pi}_{y,w}-\pi_{y,w} &= \sum_{\atop{z\in W_n}{y \prec z 
\preceq w}} \lambda_{y,z}\,\pi_{z,w} \qquad \mbox{if $y \prec w$}.\tag{KL3'}
\end{align*}
Since $\lambda_{y,w}\in {\Z}[q,q^{-1}]$, it then follows that 
$\pi_{y,w}\in q^{-1}{\Z}[q^{-1}]$ if $y\prec w$.

\begin{cor} \label{baseEC} Let $w\in W_n$.
\begin{itemize}
\item[(a)] $\bC_w$ can be written as an $A$-linear combination of $E_w$ and
terms $E_y$ where $y \prec w$.
\item[(b)] $E_w$ can be written as an $A$-linear combination of $\bC_w$
and terms $\bC_y$ where $y \prec w$.
\end{itemize}
\end{cor}

\begin{proof} (a) See the above expression for $\bC_w$. (b) Argue as
in the proof of Corollary~\ref{cor-rel2}.
\end{proof}

The next two results describe the action of $\bC_t$ and $\bC_{s_i}$ on
$E_w$.

\begin{lem} \label{lem-c3} Let $w\in W_n$ and $s=s_i$ for some $1\leq i 
\leq n-1$. Then $\bC_sE_w$ is an $A$-linear combination of terms $E_z$ where 
$l:=l_t(z)= l_t(w)$ and $\sigma_z b_z^{-1}\leq_{\cL,l} \sigma_w b_w^{-1}$. 
\end{lem}

\begin{proof} Recall that $E_w=T_{a_w}\,\bC_{a_l}\,\bC_{\sigma_wb_w^{-1}}$. 
Now $\bC_{s}=T_{s}+q^{-1}T_1$ and so 
\[ \bC_{s}E_w=T_{s}E_w+q^{-1}E_w.\]
By Deodhar's Lemma (see \cite[2.1.2]{ourbuch}), there are three
cases to consider.

(i) $sa_w\in X_{l,n-l}$ and $l(sa_w)>l(a_w)$. Then 
\[T_s\,E_w= T_s\,T_{a_w}\,\bC_{a_l}\,\bC_{\sigma_wb_w^{-1}}=
T_{sa_w}\,\bC_{a_l}\,\bC_{\sigma_wb_w^{-1}}= E_{sw}\]
and so $\bC_sE_w=E_{sw}+q^{-1}E_w$. Since $sw=(sa_w)a_l\sigma_wb_w^{-1}$, 
the required conditions are satisfied.

(ii) $sa_w\in X_{l,n-l}$ and $l(sa_w)<l(a_w)$. Then $T_sT_{a_w}=
T_{sa_w}+(q-q^{-1})T_{a_w}$ and so 
\[T_sE_w=E_{sw}+(q-q^{-1})E_{w}.\]
This yields $\bC_sE_w=E_{sw}+qE_w$. Since, again, $sw=(sa_w)a_l
\sigma_wb_w^{-1}$, the required conditions are satisfied.

(iii) $sa_w=a_ws'$ for some $s'\in \Sigma_{l,n-l}$. Then $l(sa_w)=l(a_w)+1=
l(a_ws')$ and so $T_sT_{a_w}=T_{sa_w}=T_{a_ws'}=T_{a_w}T_{s'}$. This yields
\[T_sE_w=T_{a_w}\,T_{s'}\,\bC_{a_l}\,\bC_{\sigma_wb_w^{-1}}=
T_{a_w}\,\bC_{s'}\,\bC_{a_l}\,\bC_{\sigma_wb_w^{-1}} -q^{-1}E_w\]
and so 
\[ \bC_sE_w=T_{a_w}\,\bC_{s'}\,\bC_{a_l}\,\bC_{\sigma_wb_w^{-1}}.\]
Now Lemma~\ref{lem-c1} shows that 
\[ \bC_{s'}\,\bC_{a_l}\,\bC_{\sigma_w b_w^{-1}}=
\bC_{a_l}\,\bC_{a_ls'a_l}\bC_{\sigma_w b_w^{-1}}.\]
Since $a_ls'a_l \in\fS_{l,n-l}$, we can express $\bC_{a_ls'a_l}
\bC_{\sigma_w b_w^{-1}}$ as an $A$-linear combination of terms 
$\bC_{\rho b^{-1}}$ where $\rho \in \fS_{l,n-l}$ and $b \in X_{l,n-l}$ are 
such that $\rho b^{-1} \leq_{\cL,l} \sigma_w b_w^{-1}$. We conclude that 
$\bC_sE_w$ is an $A$-linear combination of terms $E_y$ where $a_y=a_w$, 
$l:=l_t(y)=l_t(w)$ and $\sigma_y b_y^{-1}\leq_{\cL,l} \sigma_w b_w^{-1}$.
\end{proof} 

\begin{lem} \label{lem-c4} Let $w \in W_n$ and $l=l_t(w)$. Then 
$\bC_tE_w$ is an $A$-linear combination of terms $E_z$ where $l_t(z)>l$
or where $l_t(z)=l$ and $\sigma_z b_z^{-1} \leq_{\cL,l} \sigma_w b_w^{-1}$.
\end{lem}

\begin{proof} The following argument is inspired from the proof of 
Dipper--James--Murphy \cite[Lemma~4.9]{DiJa95}. Write $w=a_wa_l\sigma_w 
b_w^{-1}$. We distinguish three cases.

{\em Case 1.} We have $l=0$. Then $a_w=1$ and so $E_w= \bC_{\sigma_w 
b_w^{-1}}$. By Proposition~\ref{maintl}, $\bC_tE_w$ is a linear combination 
of terms $\bC_z$ where $l_t(z)\geq 1$. Using Corollary~\ref{baseEC}(a), we 
see that $\bC_tE_w$ can also be written as a linear combination of term 
$E_{z'}$ where $l_t(z')\geq 1$.

{\em Case 2.} We have $l\geq 1$ and the element $a_w$ fixes the number~$1$.
(Here, we regard $a_w$ as an element of $\fS_n$.) Then $T_t$ commutes with 
$T_{a_w}$. Since $l(ta_l)< l(a_l)$, we have $T_t \bC_{a_l}= -Q^{-1}\bC_{a_l}$. 
So $\bC_tE_w$ is a multiple of $E_w$ and we are done in this case.

{\em Case 3.} We have $l\geq 1$ and the element $a_w$ does not fix the 
number~$1$. Then we consider the Young subgroup $\fS_{1,n-1} \subset \fS_n$. 
We can write $a_w\in \fS_n$ as a product of an element of $\fS_{1,n-1}$ 
times a distinguished right coset representative of $\fS_{1,n-1}$ in $\fS_n$. 
These coset representatives are given by 
\[ \{1,\quad s_1,\quad s_1s_2,\quad s_1s_2s_3,\quad \ldots,\quad s_1s_2s_3
\cdots s_{n-1}\}.\]
Thus, we have $a_w=\sigma s_1s_2 \cdots s_m$ for some $m\in \{0,1,\ldots,
n-1\}$ where $l(a_w)=m+l(\sigma)$. Now, the fact that $a_w\in X_{l,n-l}$ 
implies that we must have $m=l$ and so 
\[a_w=\sigma s_1s_2 \cdots s_l \quad \mbox{for some $\sigma \in 
\fS_{1,n-1}$ such that $l(a_w)=l+l(\sigma)$}.\]
This yields
\[ T_tT_{a_w}=T_tT_\sigma T_{s_1s_2\cdots s_l}=T_\sigma T_t
T_{s_1s_2\cdots s_l}=T_\sigma T_{ts_1s_2\cdots s_l}.\]
Using the expression in Lemma~\ref{lem-c2}, we obtain
\[T_tT_{a_w}\bC_{a_l}=T_\sigma \,T_{ts_1s_2\cdots s_l}\,
\bC_{a_l}=T_\sigma \,\bC_{a_{l+1}} +T_\sigma \,h(a_l)\, \bC_{a_l},\]
where $h(a_l)$ is an $A$-linear combination of basis elements $T_\pi$ with
$\pi \leq s_1\cdots s_l$. This yields
\[ T_tE_w=T_\sigma \,\bC_{a_{l+1}}\bC_{\sigma_wb_w^{-1}} +T_\sigma \,h(a_l)\, 
\bC_{a_l}\bC_{\sigma_w b_w^{-1}}.\]
Now Lemma~\ref{lem-c3} shows that $T_\sigma\,h(a_l)\, \bC_{a_l}\,
\bC_{\sigma_wb_w^{-1}}$ is a linear combination of terms $E_z$ 
where $l=l_t(z)$ and $\sigma_zb_z^{-1} \leq_{\cL,l}\sigma_wb_w^{-1}$. On 
the other hand, by Proposition~\ref{maintl}, $T_\sigma \bC_{a_{l+1}} 
\bC_{\sigma_wb_w^{-1}}$ is a linear combination of terms $\bC_{w'}$ where 
$l_t(w')\geq l+1$. Hence this is also a linear combination of terms 
$E_{z'}$ where $l_t(z')\geq l+1$.
\end{proof}

\begin{thm} \label{strongL} Let $x,y \in W_n$ be such that $l:=l_t(x)=
l_t(y)$. Then we have $x\leq_{\cL} y$ if and only if $\sigma_x b_x^{-1}
\leq_{\cL,l} \sigma_y b_y^{-1}$. 
\end{thm}

\begin{proof} First assume that $x\leq_{\cL} y$. We must show that 
$\sigma_x b_x^{-1} \leq_{\cL,l} \sigma_y b_y^{-1}$. Now, by definition, 
there exists a sequence $x=x_0,x_1,\ldots,x_k=y$ such that $x_{i-1} 
\leftarrow_{\cL} x_i$ for all $i$. By Proposition~\ref{maintl}, we have 
$l_t(x_{i-1} )\geq l_t(x_i)$ for all $i$. 
Since $l_t(x)=l_t(y)$, we conclude that all $x_i$ have the same $t$-length. 
Thus, it is enough to consider the case where $x \leftarrow_{\cL} y$, that 
is, we have that $\bC_x$ occurs in $\bC_s\bC_y$, for some $s \in \{t,s_1,
\ldots, s_{n-1}\}$. 

Assume first that $s=s_i$ for some $i \in \{1,\ldots,n-1\}$. By 
Corollary~\ref{baseEC}(a), we can write $\bC_y$ as an $A$-linear combination 
of $E_w$ where $w \preceq y$. So $\bC_s\bC_y$ is an $A$-linear combination of  
terms of the form $\bC_sE_w$ where $w \preceq y$. Now consider such a 
term. By Lemma~\ref{lem-c3}, $\bC_sE_w$ is a linear combination of terms 
$E_z$ where $\sigma_z b_z^{-1} \leq_{\cL,l} \sigma_w b_w^{-1}$. 
Consequently, by Corollary~\ref{baseEC}(b), $\bC_sE_w$ is a linear 
combination of terms $\bC_z$ where $\sigma_z b_z^{-1}\leq_{\cL,l} \sigma_w 
b_w^{-1}$, as required.

Now assume that $s=t$. By Corollary~\ref{baseEC}(a), we can write $\bC_y$ as 
an $A$-linear combination of $E_w$ where $w \preceq y$. So $\bC_t\bC_y$ is 
an $A$-linear combination of terms of the form $\bC_tE_w$ where $w \preceq y$. 
By Lemma~\ref{lem-c4} and Corollary~\ref{baseEC}(b), we can write any
such term as a linear combination of terms $\bC_z$ where $l_t(z)>l$ or
$l_t(z)=l$ and $\sigma_zb_z^{-1}\leq_{\cL,l}\sigma_wb_w^{-1}$.  

Summarizing, we have shown that $\bC_t\bC_y$ is a linear combination of
terms $\bC_z$ where $l=l_t(z)=l_t(y)$ and $\sigma_zb_z^{-1} \leq_{\cL,l} 
\sigma_yb_y^{-1}$, and terms $\bC_{w'}$ where $l_t(w')>l$. Hence, since 
$l_t(x)=l$, we must have $\sigma_x b_x^{-1} \leq_{\cL,l} \sigma_yb_y^{-1}$,
as required.

Conversely, let us assume that $\sigma_x b_x^{-1} \leq_{\cL,l} \sigma_y 
b_y^{-1}$. We must show that $x\leq_{\cL} y$. Again, it is enough to
consider the case where  $\sigma_x b_x^{-1} \leftarrow_{\cL,l}
\sigma_y b_y^{-1}$, that is, $\bC_{\sigma_x b_x^{-1}}$ occurs in $\bC_s
\bC_{\sigma_y b_y^{-1}}$ for some $s=s_i$ where $i \neq l$. Thus, writing 
\[ \bC_{s_i}\bC_{\sigma_y b_y^{-1}}=\sum_{\pi \in \fS_{l,n-l}} 
\sum_{z \in X_{l,n-l}} h_{s_i,\sigma_yb_y^{-1},\pi z^{-1}} \bC_{\pi z^{-1}},\]
we have $h_{s_i,\sigma_yb_y^{-1},\sigma_xb_x^{-1}} \neq 0$. 
Multiplying the above equation on the left by $\bC_{a_l}$ and using 
Lemma~\ref{lem-c1}, we conclude that 
\begin{align*}
\bC_{s'}\bC_{a_l\sigma_y b_y^{-1}}&=\bC_{s'}\bC_{a_l}\bC_{\sigma_y b_y^{-1}}=
\bC_{a_l}\bC_{s_i}\bC_{\sigma_y b_y^{-1}}\\
&=\sum_{\pi \in \fS_{l,n-l}} \sum_{z \in X_{l,n-l}} 
h_{s_i,\sigma_y b_y^{-1},\pi z^{-1}} \bC_{a_l\pi z^{-1}}, 
\end{align*}
where $s'=a_ls_ia_l \in \fS_{l,n-l}$. Considering the term corresponding
to $\pi=\sigma_x$ and $z=b_x$, we see that $a_l\sigma_x b_x^{-1}
\leq_{\cL} a_l\sigma_y b_y^{-1}$. Finally, this yields
\[x=a_x a_l\sigma_x b_x^{-1} \sim_{\cL} a_l\sigma_x b_x^{-1} 
\leq_{\cL} a_l\sigma_y b_y^{-1}\sim_{\cL} a_xa_l\sigma_x b_x^{-1}=y, \]
by Theorem~\ref{mainbi}.
\end{proof} 

The above result has two immediate applications.

Firstly, it provides a refinement of Theorem~\ref{mainbi}.  Indeed, if 
we have $x\sim_{\cL} y$, then Theorem~\ref{strongL} shows that $\sigma_x 
b_x^{-1} \sim_{\cL,l} \sigma_y b_y^{-1}$ and, hence, $b_x=b_y$ and 
$\sigma_x \sim_{\cL,l} \sigma_y$ (by Proposition~\ref{propa21} and
Lemma~\ref{lema23}). 

Secondly, it refines the methods that Bonnaf\'e used in \cite{BI2}.
Indeed, we obtain a new proof of the following statement concerning 
the two-sided Kazhdan--Lusztig pre-order.

\begin{cor}[See Bonnaf\'e \protect{\cite{BI2}}] \label{mainbi2} Let 
$x,y\in W_n$. Then the following hold.
\begin{itemize}
\item[(a)] If $l:=l_t(x)=l_t(y)$ and $x\leq_{\cLR} y$, then $\sigma_x 
\leq_{\cLR,l} \sigma_y$. 
\item[(b)] If $x\sim_{\cLR} y$, then $l:=l_t(x)=l_t(y)$ and $\sigma_x 
\sim_{\cLR,l} \sigma_y$. 
\end{itemize}
\end{cor}

\begin{proof} (a) Assume that $l:=l_t(x)=l_t(y)$. To prove the implication 
``$x\leq_{\cLR} y \Rightarrow \sigma_x \leq_{\cLR,l} \sigma_y$'', we may 
assume without loss of generality that $x\leq_{\cL} y$ or $x^{-1}
\leq_{\cL} y^{-1}$ (since these are the elementary steps in the definition 
of $\leq_{\cLR}$.) If $x \leq_{\cL} y$, then Theorem~\ref{strongL} and 
Proposition~\ref{propa21} immediately yield $\sigma_x \leq_{\cLR,l} 
\sigma_y$, as required. Assume now that $x^{-1}\leq_{\cL} y^{-1}$. We have 
\[ x^{-1}=(a_xa_l\sigma_xb_x^{-1})^{-1}=b_xa_l (a_l\sigma_x^{-1}a_l)
a_x^{-1}\]
and so 
\[a_{x^{-1}}=b_x,\quad  \sigma_{x^{-1}}=\sigma_l\sigma_x^{-1} \sigma_l,
\quad b_{x^{-1}}=a_x\]
where $\sigma_l$ is the longest element of $\fS_l$.  Note that $a_l=w_l
\sigma_l$ where $w_l$ is the longest element in $W_l$, and that $w_l$ 
commutes with all elements of $\fS_{l,n-l}$; see \cite[\S 4]{BI}. A similar 
remark applies to $y=a_ya_l\sigma_yb_y^{-1}$. 
Now Theorem~\ref{strongL} and Proposition~\ref{propa21} imply $\sigma_l 
\sigma_x^{-1} \sigma_l \leq_{\cLR,l} \sigma_l\sigma_y^{-1}\sigma_l$.
Furthermore, conjugation with $\sigma_l$ defines a Coxeter group 
automorphism of $\fS_{l,n-l}$ and, hence, preserves the Kazhdan--Lusztig 
pre-order relations $\leq_{\cL,l}$, $\leq_{\cR,l}$ and $\leq_{\cLR,l}$; see 
\cite[Cor.~11.7]{Lusztig03}. Consequently, we have $\sigma_x^{-1} 
\leq_{\cLR,l} \sigma_y^{-1}$. Finally, note that inversion certainly 
preserves the two-sided pre-order $\leq_{\cLR,l}$. Hence we have 
$\sigma_x \leq_{\cLR,l} \sigma_y$, as desired.

(b) If $x\leq_{\cLR} y$, then $l_t(y)\leq l_t(x)$ by 
Proposition~\ref{maintl}. Hence, if $x\sim_{\cLR} y$, then we automatically 
have $l:=l_t(x)=l_t(y)$ and (a) yields $\sigma_x\sim_{\cLR,l} \sigma_y$.
\end{proof}

Now our efforts will be rewarded. Combining Example~\ref{conjA} with
Theorem~\ref{mainbi}, Theorem~\ref{strongL} and Corollary~\ref{mainbi2},
we obtain:

\begin{thm} \label{relsn} Recall that we are in the ``asymptotic case''
in type $B_n$. Then the following implication holds for all $x,y \in W_n$:
\begin{equation*}
 x\leq_{\cL} y \quad \mbox{and} \quad x\sim_{\cLR} y \quad\Rightarrow 
\quad x \sim_{\cL} y.\tag{$\spadesuit$}
\end{equation*}
\end{thm}

\begin{proof} Let $x,y \in W_n$ be such that $x \leq_{\cL} y$ and 
$x \sim_{\cLR} y$. First of all, Corollary~\ref{mainbi2} implies that 
$l:=l_t(x)=l_t(y)$ and $\sigma_x \sim_{\cLR,l} \sigma_y$. Furthermore, 
Theorem~\ref{strongL}  implies that $\sigma_x b_x^{-1} \leq_{\cL,l} 
\sigma_y b_y^{-1}$. Thus, the hypotheses of Conjecture~\ref{conja22}
are satisfied for the elements $\sigma_x b_x^{-1}$ and $\sigma_y b_y^{-1}$
in the symmetric group $\fS_n$, where we consider the parabolic subgroup
$\fS_{l,n-l}$. Hence Example~\ref{conjA} implies that $b_x=b_y$ and
$\sigma_x \sim_{\cL,l} \sigma_y$. Then Theorem~\ref{mainbi} yields 
$x \sim_{\cL} y$, as desired.
\end{proof}

\begin{cor} \label{2cells} The sets $\{\fR_\lambda\mid \lambda \in 
\Lambda_n\}$ are precisely the two-sided cells of $W_n$.
\end{cor}

\begin{proof} Once ($\spadesuit$) is known to hold, two elements $x,y\in 
W_n$ lie in the same two-sided cell if and only if there exists a sequence
$x=x_0,x_1,\ldots,x_k=y$ of elements in $W_n$ such that, for each $i$,
we have $x_{i-1}\sim_{\cL} x_i$ or $x_{i-1}\sim_{\cR} x_i$. Hence the
assertion is an immediate consequence of Corollary~\ref{unique}. 
\end{proof}

\section{On the left cell representations in type $B_n$} \label{sec-repB}
We keep the set-up of the previous  section, where $W_n$ is a Coxeter group 
of type $B_n$ and where we consider the Kazhdan--Lusztig cells in the
``asymptotic case''. Recall the partition
\[W_n=\coprod_{\lambda \in \Lambda_n} \fR_\lambda,\]
where $\Lambda_n$ is set of all pairs of partitions of total size~$n$. An
element $w\in W_n$ belongs to $\fR_\lambda$ if and only if $w$ corresponds
to a pair of bitableaux of shape $\lambda$ under the generalized
Robinson--Schensted correspondence. By Corollary~\ref{2cells}, each set 
$\fR_\lambda$ is a two-sided cell. 

Recall that we denote by $\Irr(\bH_{n,K})$ the set of irreducible 
characters of $\bH_{n,K}$. For any left cell $\fC$, we denote by 
$\chi_{\fC}$ the character afforded by the $\bH_{n,K}$-module $[\fC]_K=
K\otimes_A [\fC]_A$. 

\begin{thm}[Bonnaf\'e--Iancu \protect{\cite{BI}} and Bonnaf\'e
\protect{\cite[\S 5]{BI2}}] \label{mainbichar} In the above setting, we 
have $\chi_{\fC}\in \Irr(\bH_{n,K})$ for any left cell $\fC$ in 
$W_n$.  Furthermore, let $\fC,\fC_1$ be left cells and assume that 
$\fC\subseteq \fR_\lambda$, $\fC_1\subseteq \fR_\mu$ where $\lambda,\mu\in
\Lambda_n$.  Then the characters $\chi_{\fC}$ and 
$\chi_{\fC_1}$ are equal if and only if $\lambda=\mu$. 
\end{thm}

(This is another example where the discussion in Remark~\ref{asymB} applies: 
the above statements were first proved in \cite[\S 7]{BI} for the weight 
function $L^\circ\colon W_n \rightarrow {\Z}^2$. Using Remark~\ref{asymB}(a),
one easily shows that $[\fC]_A=A \otimes_{A^\circ} [\fC]_{A^\circ}$ where 
$A$ is regarded as an $A^\circ$-module via the map $\theta\colon A^\circ 
\rightarrow A$.)

The main result of this section is Theorem~\ref{eqcellsB} which shows that 
we even have $\fC \approx \fC_1$ for any two left cells $\fC,\fC_1 
\subseteq \fR_\lambda$, where ``$\approx$'' is the relation introduced in 
Definition~\ref{eqcells}. 

Let us fix a pair of partitions $\lambda=(\lambda_1, \lambda_2)\in 
\Lambda_n$ and let $\fC \subseteq \fR_\lambda$ be a left cell. We
set $l:=|\lambda_2|$. By \cite[Prop.~4.8]{BI}, we have $l_t(w)=l$ for 
all $w\in \fR_\lambda$. In particular, we have $l_t(w)=l$ for all 
$w\in \fC$. Now recall the decomposition $w=a_wa_l\sigma_wb_w^{-1}$ for 
any element $w\in W_n$, where $l=l_t(w)$.  We set
\[\overline{\fC}:=\{\sigma \in \fS_{l,n-l} \mid \sigma=\sigma_w \mbox{ for 
some $w\in \fC$}\}.\]
By Theorem~\ref{mainbi}, $\overline{\fC}$ is a left cell in $\fS_{l,n-l}$. 
Next recall that $\fS_{l,n-l}=\fS_l \times \fS_{[l+1,n]}$ where $\fS_{[l+1,n]}
\cong \fS_{n-l}$. It is well-known and easy to check that the Kazhdan--Lusztig
pre-order relations are compatible with direct products; in particular, every 
left cell in $\fS_{l,n-l}$ is a product of a left cell in $\fS_l$ and a left 
cell in $\fS_{[l+1,n]}$. Thus, we can write 
\[ \overline{\fC}=\overline{\fC}^{\,(l)}\cdot \overline{\fC}^{\,(n-l)}\]
where $\overline{\fC}^{\,(l)}$ is a left cell in $\fS_l$ and 
$\overline{\fC}^{\,(n-l)}$ is a left cell in $\fS_{[l+1,n]}$. 
We use the explicit dot to indicate that the lengths of elements add
up in this product: we have $l(\sigma\tau)=l(\sigma)+l(\tau)$ for $\sigma 
\in \overline{\fC}^{\,(l)}$ and $\tau\in \overline{\fC}^{\,(n-l)}$. By 
Theorem~\ref{mainbi}, we have $b_x=b_y$ for all $x,y\in \fC$. Let us 
denote $b=b_w$ for $w\in \fC$. Then we have 
\[ \fC=X_{l,n-l}\cdot  a_l \cdot \overline{\fC}\cdot  b^{-1}=
\{c a_l\sigma b^{-1} \mid c \in X_{l,n-l},\sigma\in \overline{\fC}\}.\]
A first reduction is provided by the following result:

\begin{lem}[Bonnaf\'e--Iancu \protect{\cite[Prop.~7.2]{BI}} and
Remark~\ref{asymB}] \label{redbi} In the above setting, $\fC b$ is a left 
cell and we have 
\[ \fC \approx \fC b=X_{l,n-l}\cdot  a_l \cdot \overline{\fC}.\]
\end{lem}

Now we can state the main result of this section. Again, this is in
complete formal analogy to the situation in the symmetric group 
$\fS_n$; see Example~\ref{expA}(b). 

\begin{thm} \label{eqcellsB} Let $\lambda \in \Lambda_n$. Then we have 
$\fC \approx \fC_1$ for all left cells $\fC,\fC_1 \subseteq \fR_\lambda$.
Recall that this means that there exists a bijection 
$\fC\stackrel{\sim}{\rightarrow} {\fC}_1$, $x \mapsto x_1$, 
such that $h_{w,x,y}=h_{w,x_1,y_1}$ for all $w\in W_n$ and all $x,y\in\fC$.

The bijection $x \mapsto x_1$ is uniquely determined by the condition that 
$x_1\in \fC_1$ is the unique element in the same right cell as~$x\in \fC$.
\end{thm}

\begin{proof} First note that the second statement (concerning the uniqueness
of the bijection) is a consequence of the first. Indeed, if there exists
a bijection $\fC \stackrel{\sim}{\rightarrow} \fC_1$, $x\mapsto x_1$,
satisfying ($\heartsuit$), then Proposition~\ref{uniqueeq} shows that
$x \sim_{\cR} x_1$ for any $x\in \fC$. But Corollary~\ref{unique} shows
that two elements which are in the same right cell and in the same left cell
are equal.  Hence the element $x_1$ is uniquely determined by the condition
that $x_1 \sim_{\cR} x$.

To establish the existence of such a bijection, let $\lambda=(\lambda_1,
\lambda_2)$ and set $l:=|\lambda_2|$. Let $\fC,\fC_1 \subseteq \fR_\lambda$
be two left cells. We set
\begin{align*}
\overline{\fC}&:=\{\sigma \in \fS_{l,n-l} \mid \sigma=\sigma_w \mbox{ for 
some $w\in \fC$}\},\\
\overline{\fC}_1&:=\{\sigma \in \fS_{l,n-l} \mid \sigma=\sigma_w \mbox{ for
some $w\in \fC_1$}\};
\end{align*}
by the above discussion, these are left cells in $\fS_{l,n-l}$. Furthermore,
we can write 
\[ \overline{\fC}=\overline{\fC}^{\,(l)}\cdot \overline{\fC}^{\,(n-l)}
\qquad \mbox{and}\qquad \overline{\fC}_1=\overline{\fC}_1^{\,(l)}\cdot 
\overline{\fC}_1^{\,(n-l)}\]
where $\overline{\fC}^{\,(l)}$, $\overline{\fC}_1^{\,(l)}$ are left cells
in $\fS_l$ and $\overline{\fC}^{\,(n-l)}$,  $\overline{\fC}_1^{\,(n-l)}$  
are left cells in $\fS_{[l+1,n]}$. We claim that 
\begin{equation*}
\overline{\fC}^{\,(l)} \approx \overline{\fC}_1^{\,(l)}, \qquad
\overline{\fC}^{\,(n-l)}\approx \overline{\fC}_1^{\,(n-l)},
\qquad \overline{\fC}\approx \overline{\fC}_1. \tag{$*$}
\end{equation*}
Indeed, by Example~\ref{expA}, the classical Robinson--Schensted 
correspondence associates to a left cell of $\fS_l$ a partition of $l$ and 
to a left cell in $\fS_{[l+1,n]}$ a partition of $n-l$. Thus, we can 
associate a pair of partitions to $\overline{\fC}$. By \cite[4.7]{BI}, that 
pair of partitions is given by $(\lambda_2,\lambda_1)$. A similar remark 
applies to $\overline{\fC}_1$, where we obtain the same pair of partitions.  
Now ($*$) follows from Example~\ref{expA}(b) and the compatibility of left 
cells with direct products.

To continue the proof it is sufficient, by Lemma~\ref{redbi}, to 
consider the case where 
\[\fC=X_{l,n-l}\cdot  a_l\overline{\fC}\qquad \mbox{and}\qquad 
\fC_1=X_{l,n-l}\cdot a_l\overline{\fC}_1.\]
In this situation, we note that the sets $a_l\overline{\fC}$ and $a_l
\overline{\fC}_1$ are contained in the parabolic subgroup 
\[ W_{l,n-l}=W_l \times \fS_{[l+1,n]} \quad \mbox{where}\quad  
W_l=\langle t,s_1,\ldots,s_{l-1}\rangle \quad \mbox{(type $B_l$)}.\]
By \cite[4.1]{BI}, we have $a_l=w_l\sigma_l$ where $w_l$ is 
the longest element in $W_l$ and $\sigma_l$ is the longest element in
$\fS_l$. Since multiplication with the longest element preserves left 
cells, the sets $\sigma_l\overline{\fC}^{\,(l)}$ and $\sigma_l
\overline{\fC}_1^{\,(l)}$ are left cells in $\fS_l$. Hence ($*$) and 
Lemma~\ref{cellw0} show that 
\[\sigma_l\overline{\fC}^{\,(l)}\approx\sigma_l\overline{\fC}_1^{\,(l)}.\] 
Applying Theorem~\ref{mainbi} to the group $W_l$, we notice that every
left cell in $\fS_l$ also is a left cell in $W_l$. Hence the sets 
$\sigma_l \overline{\fC}^{\,(l)}$ and $\sigma_l \overline{\fC}_1^{\,(l)}$
are left cells in $W_l$. Then multiplication with the longest element 
$w_l\in W_l$ and Lemma~\ref{cellw0} yield that 
\[ a_l \overline{\fC}^{\,(l)}= w_l(\sigma_l\overline{\fC}^{\,(l)})\approx
w_l(\sigma_l \overline{\fC}_1^{\,(l)})=a_l\overline{\fC}_1^{\,(l)},\]
where the above sets are left cells in $W_l$. Using the compatibility of 
left cells with direct products, we obtain 
\[a_l\overline{\fC}=(a_l\overline{\fC}^{\,(l)}) \cdot 
\overline{\fC}^{\,(n-l)}\approx (a_l\overline{\fC}_1^{\,(l)}) 
\cdot \overline{\fC}_1^{\,(n-l)}=a_l\overline{\fC}_1,\]
where the above sets are left cells in $W_{l,n-l}$. 

Thus, we have two left cells in the parabolic subgroup $W_{l,n-l}$ which
are related by ``$\approx$''. Now let $\hat{X}_{l,n-l}$ be the set of 
distinguished left coset representatives of $W_{l,n-l}$ in $W_n$. We
certainly have $X_{l,n-l} \subseteq \hat{X}_{l,n-l}$ and so
\[\begin{array}{ccc}
\fC =X_{l,n-l}\cdot a_l\overline{\fC}&\subseteq&\hat{X}_{l,n-l}\cdot
a_l\overline{\fC},\\
\fC_1=X_{l,n-l}\cdot a_l\overline{\fC}_1&\subseteq&\hat{X}_{l,n-l}\cdot
a_l\overline{\fC}_1.
\end{array}\]
By Theorem~\ref{inducell}, the sets $\hat{X}_{l,n-l} \cdot a_l
\overline{\fC}$ and $\hat{X}_{l,n-l} \cdot a_l\overline{\fC}_1$ are
both unions of left cells in $W_n$ and we have 
\[\hat{X}_{l,n-l} \cdot a_l \overline{\fC}\approx
\hat{X}_{l,n-l} \cdot a_l\overline{\fC}_1;\]
see Proposition~\ref{induind0}. Furthermore, by 
Corollary~\ref{ind1}, there is a left cell
\[ \tilde{\fC} \subseteq \hat{X}_{l,n-l} \cdot a_l \overline{\fC}
\qquad \mbox{such that}\qquad  \tilde{\fC} \approx \fC_1.\]
It remains to show that $\fC=\tilde{\fC}$. This can be seen as follows.
Since $\tilde{\fC} \approx \fC_1$, we have $\tilde{\fC} \subseteq 
\fR_\lambda$. In particular, all elements in $\tilde{\fC}$ must have 
$t$-length~$l$. Now we leave it as an exercice to the reader to check that 
\[ X_{l,n-l}\cdot W_{l,n-l}=\{ w \in W_n \mid l_t(w)\leq l\}.\]
Hence we must have $\tilde{\fC}\subseteq X_{l,n-l} \cdot W_{l,n-l}$. 
On the other hand, we also have $\tilde{\fC} \subseteq \hat{X}_{l,n-l}
\cdot a_l \overline{\fC}$. Since $X_{l,n-l} \subseteq \hat{X}_{l,n-l}$
and $a_l\overline{\fC} \subseteq W_{l,n-l}$, we conclude that 
\[\tilde{\fC}\subseteq \big(X_{l,n-l}\cdot W_{l,n-l}\big) \cap 
\big(\hat{X}_{l,n-l} \cdot a_l \overline{\fC}\big)=X_{l,n-l}
\cdot a_l \overline{\fC}=\fC\]
and so $\tilde{\fC}=\fC$, as required.
\end{proof}

Following Graham--Lehrer \cite[Definition~1.1]{GrLe}, a quadruple 
$(\Lambda,M,C,*)$ is called a ``cell datum'' for $\bH_n$ if the following 
conditions are satisfied.
\begin{itemize}
\item[(C1)] $\Lambda$ is a partially ordered set, $\{M(\lambda) \mid 
\lambda \in \Lambda\}$ is a collection of finite sets  and 
\[ C \colon \coprod_{\lambda \in \Lambda} M(\lambda)\times M(\lambda)
\rightarrow \bH_n \]
is an injective map whose image is an $A$-basis of $\bH_n$;
\item[(C2)] If $\lambda \in \Lambda$ and $S,T\in M(\lambda)$, write 
$C(S,T)=C_{S,T}^\lambda \in \bH_n$. Then $* \colon \bH_n \rightarrow \bH_n$
is an $A$-linear anti-involution such that $(C_{S,T}^\lambda)^*=
C_{T,S}^\lambda$. 
\item[(C3)] If $\lambda \in \Lambda$ and $S,T\in M(\lambda)$, then for any
element $h \in \bH_n$ we have 
\[ hC_{S,T}^\lambda\equiv \sum_{S'\in M(\lambda)} r_h(S',S)\, 
C_{S',T}^\lambda\quad \bmod \bH_n(<\lambda),\]
where $r_h(S',S) \in A$ is independent of $T$ and where $\bH_n(<\lambda)$
is the $A$-submodule of $\bH_n$ generated by $\{C_{S'',T''}^\mu
\mid \mu <\lambda; S'',T''\in M(\mu)\}$.
\end{itemize}
In this case, we call the basis $\{C_{S,T}^\lambda\}$ a 
{\bf ``cellular basis''} of $\bH_n$. 

One reason for the importance of a cellular structure lies in the fact 
that it leads to a general theory of ``Specht modules'' and various 
applications concerning modular representations; see \cite{GrLe} for more 
details. Graham and Lehrer \cite[\S 5]{GrLe} already showed that $\bH_n$ has 
a cellular structure, where they use a mixture of the Kazhdan--Lusztig basis
and the standard basis. The point of the following result is that the
Kazhdan--Lusztig basis in the ``asymptotic case'' directly gives a cellular 
structure. (The relation between the two structures will be discussed 
elsewhere.)

\begin{cor} \label{grle} Recall that we are in the ``asymptotic case'' in
type $B_n$. Then the Kazhdan--Lusztig basis $\{\bC_w\mid w\in W_n\}$ is a 
``cellular basis'' of $\bH_n$. 
\end{cor}

\begin{proof} We specify a ``cell datum'' as follows. First of all, let 
$\Lambda:=\Lambda_n$, the set of all pairs of partitions of total size~$n$. 
By Corollary~\ref{2cells}, these parametrize the two-sided cells of $W_n$. 
Hence we can define a partial order ``$\leq$'' on $\Lambda$ by
\[ \lambda \leq \mu \quad \mbox{if} \quad x \leq_{\cLR} y \mbox{ for some
$x\in \fR_\lambda$, $y \in \fR_\mu$}.\]
(More explicitly, we could use the dominance order on bipartitions;
see \cite[Prop.~5.3]{GeIa05}.)
Next, for each $\lambda \in \Lambda_n$, let $M(\lambda):=\fT_\lambda$,
the set of $n$-standard bitableaux of shape $\lambda$. By 
Corollary~\ref{unique}, we have a bijection
\[ \fT_\lambda \times \fT_\lambda \rightarrow \fR_\lambda,
\qquad (T,T')\mapsto w_\lambda(T,T')\]
such that $w_\lambda(T,T')^{-1}=w_\lambda(T',T)$ for all $T,T'\in 
\fT_\lambda$. We set 
\[ C_{S,T}^\lambda:=\bC_{w_{\lambda}(S,T)} \qquad \mbox{for $\lambda
\in \Lambda_n$ and $S,T\in \fT_\lambda$}.\]
Then the map 
\[ C\colon \coprod_{\lambda \in \Lambda_n} \fT_\lambda \times \fT_\lambda
\rightarrow \bH_n, \qquad (S,T) \mapsto C_{S,T}^\lambda \qquad (S,T \in
\fT(\lambda)),\]
satisfies the requirements in (C1).

We define $*\colon \bH_n\rightarrow \bH_n$ by $T_w^*=T_w^\flat=
T_{w^{-1}}$ for all $w\in W_n$. This is an $A$-linear anti-involution
such that $\bC_w^*=\bC_{w^{-1}}$ for all $w\in W_n$; see the remarks in
(\ref{klpre}). Thus, we have 
\[ (C_{S,T}^\lambda)^*=\bC_{w_\lambda(S,T)}^*=\bC_{w_\lambda(S,T)^{-1}}=
\bC_{w_\lambda(T,S)}=C_{T,S}^\lambda\]
for all $\lambda \in \Lambda_n$ and $S,T\in \fT_\lambda$. Hence condition
(C2) is satisfied.

In order to check (C3), it is sufficient to assume that $h=\bC_w$ for some
$w\in W_n$. Let $\lambda \in \Lambda_n$ and $T\in \fT_\lambda$. For any
$S,S'\in \fT_\lambda$, we define
\[ r_w(S',S):=h_{w,x,x'} \quad \mbox{where} \quad \left\{\begin{array}{c}
x:=w_\lambda(S,T),\\ x':=w_\lambda(S',T).\end{array}\right.\]
Now consider the product
\[ \bC_wC_{S,T}^\lambda=\bC_w\bC_x=\sum_{y \in W_n} h_{w,x,y} \bC_y \qquad
\mbox{where $x=w_\lambda(S,T)$}.\]
If $h_{w,x,y}\neq 0$, then $y \leq_{\cL} x$. Hence Theorem~\ref{relsn}
shows that either $x\sim_{\cL} y$ or $x <_{\cLR} y$. So we can write
\[ \bC_wC_{S,T}^\lambda=\sum_{\atop{y \in W_n}{x \sim_{\cL} y}} h_{w,x,y}
\bC_y \quad \bmod \bH_n(<\lambda).\]
Using Corollary~\ref{unique}, every $y\in W_n$ such that $x\sim_{\cL} y$
has the form $y=w_\lambda(S',T)$ for some $S'\in \fT_\lambda$. So we can
rewrite the above relation as follows.
\[ \bC_wC_{S,T}^\lambda=\sum_{S' \in \fT_\lambda} r_w(S',S) C_{S',T}^\lambda
\quad \bmod \bH_n(<\lambda).\]
Finally, we must check that $r_w(S',S)$ is independent of $T$. To see
this, let $T_1\in \fT_\lambda$ and define $r_w^1(S',S):=h_{w,x_1,x_1'}$
where $x_1=w_\lambda(S,T_1)$ and $x_1'=w_\lambda(S',T_1)$. Arguing as 
above, we see that 
\[ \bC_wC_{S,T_1}^\lambda=\sum_{S' \in \fT_\lambda} r_w^1(S',S) C_{S',
T_1}^\lambda \quad \bmod \bH_n(<\lambda).\]
Hence we have 
\[ r_w(S',S)=r_w^1(S,S')\qquad \Leftrightarrow \qquad h_{w,x,x'}=
h_{w,x_1,x_1'}.\]
Now, Corollary~\ref{unique} shows that $x \sim_{\cR} x_1$, $x'\sim_{\cR} 
x_1'$, $x\sim_{\cL} x'$ and $x_1\sim_{\cL} x_1'$. Hence the desired
equality follows from Theorem~\ref{eqcellsB}.
\end{proof}

\begin{rem} \label{grrrr} The above proof is modeled on the discussion
of the Iwahori--Hecke algebra of the symmetric group $\fS_n$ in 
\cite[Example~1.2]{GrLe}. In that case, Graham and Lehrer state that (C3) 
is already implicit in Kazhdan--Lusztig \cite{KaLu} (or the work of 
Barbasch--Vogan and Vogan), which is not really the case. In fact, as the 
above proof shows, (C3) relies on the validity of both ($\heartsuit$) 
and ($\spadesuit$), and the latter was first proved by Lusztig \cite{Lu0}
(even for the symmetric group $\fS_n$). 
\end{rem}

\section{Lusztig's homomorphism from $\bH_n$ to the ring $J$} \label{sec-jB}
As a further application of the results of the previous section, we will 
now construct a new basis of $\bH_n$ with integral structure constants.
First of all, this will lead to an analogue of Lusztig's ``canonical'' 
isomorphism from $\bH_{n,K}$ onto the group algebra $KW_n$; see
Theorem~\ref{caniso}. At the end of this section, we will see that the 
subring generated by that new basis is nothing but Lusztig's ring $J$. 
To establish that identification, we will rely on the recent results of 
Iancu and the author \cite{GeIa05} concerning Lusztig's $\ba$-function. 

Recall the basic set-up from the previous sections. In particular, recall 
the partition
\[ W_n=\coprod_{\lambda \in \Lambda_n} \fR_\lambda,\]
where $\Lambda_n$ is the set of all pairs of partitions of total size~$n$.
Let us fix $\lambda \in \Lambda_n$. In the following discussion, we will 
make repeated use of the bijection
\[\fT_\lambda \times \fT_\lambda \stackrel{\sim}{\rightarrow} 
\fR_\lambda, \qquad (T,T')\mapsto w_\lambda(T,T');\]
see Corollary~\ref{unique}. Recall that this implies, in particular, that 
every left cell contains a unique element from the set 
\[\cD_n:=\{z \in W_n\mid z^2=1\}.\] 
For $z\in W_n$, we denote  by $d_z$ the unique element in $\cD_n$ such that
$z\sim_{\cL} d_z$. 

By Theorem~\ref{mainbichar}, we have $\chi_{\fC}\in\Irr(\bH_{n,K})$
for all left cells in $W_n$. This allows us to make the following 
construction.  Let $\lambda\in \Lambda_n$. We fix one left cell in
$\fR_\lambda$ and denote its elements by $\{x_1,\ldots,x_{d_\lambda}\}$. 
We have a corresponding matrix representation 
\[ \fX_{\lambda} \colon \bH_{n,K} \rightarrow M_{d_\lambda}(K), \quad 
\mbox{where} \quad \fX_\lambda^{ij}(\bC_w)=h_{w,x_j,x_i}\]
for $1\leq i,j\leq d_\lambda$; see (\ref{leftrep}). Let $\chi_\lambda
\in \Irr(\bH_{n,K})$ be the character afforded by $\fX_\lambda$. Now, if 
we vary $\lambda$, we get all irreducible characters of $\bH_{n,K}$ exactly 
once. Thus, we have a labelling
\[ \Irr(\bH_{n,K})=\{\chi_\lambda \mid \lambda \in \Lambda_n\}.\]
As in Section~\ref{sec0a}, denote by $\{\bD_w\mid w \in W_n\}$ the basis 
which is dual to the basis $\{\bC_w\mid w\in W_n\}$ with respect to the 
symmetrizing trace~$\tau$. We have the following formula:
\[ \tau=\sum_{\lambda \in \Lambda_n} \frac{1}{c_\lambda}\,\chi_{\lambda},
\qquad \mbox{where $c_\lambda \in A$ for all $\lambda \in \Lambda_n$}.\]
(In the present case, we do have $c_\lambda \in A$; see 
\cite[Theorem~9.3.5]{ourbuch}.) The main idea in this section is to apply 
Neunh\"offer's results from the end of Section~2, concerning the explicit 
Wedderburn decomposition of $\bH_{n,K}$ in terms of the products 
$\bC_x\bD_{y^{-1}}$ where $x\sim_{\cL} y$. 

The following result is inspired by an analogous result for the
symmetric group; see Neunh\"offer \cite[Kap.~VI, \S 4]{max}. It crucially
relies on Theorem~\ref{eqcellsB}.

\begin{prop} \label{maxcor2} The elements $\{\bC_z\bD_{d_z} \mid z \in W_n\}$ 
form a $K$-basis of $\bH_{n,K}$. We have the following identity for any  
$w\in W$:
\[ \bC_w=\sum_{z\in W} h_{w,d_z,z}\, c_{\lambda_z}^{-1}\,\bC_z\bD_{d_z},\]
where $\lambda_z \in \Lambda_n$ is defined by the condition that
$z \in \fR_{\lambda_z}$. Furthermore, we have the following multiplication 
rule: If $w,z\in W_n$ satisfy $w\sim_{\cR} z^{-1}$ then
\[ \bC_z\bD_{d_z} \cdot \bC_w\bD_{d_w}=c_{\lambda_z} \,\bC_u\bD_{d_u},\]
where $u \in W_n$ is the unique element such that $z\sim_{\cR} u
\sim_{\cL} w$ (see Corollary~\ref{unique}). Otherwise, we have 
$\bC_z\bD_{d_z} \cdot \bC_w\bD_{d_w}=0$.
\end{prop}

\begin{proof} First we prove the multiplication rule, by using a
representation-theoretic argument. Let $w,z \in W_n$ and suppose that 
$\bC_z\bD_{d_z}\cdot \bC_w\bD_{d_w}\neq 0$. Since $\bH_{n,K}$ is split 
semisimple, we have
\[ \fX_\lambda(\bC_z\bD_{d_z}\cdot \bC_w\bD_{d_w}) \neq 0 \qquad \mbox{for 
some $\lambda \in \Lambda_n$}.\]
Let $\fC$ be the left cell containing $z$ and $\fC_1$ be the left cell 
containing~$w$. We claim that 
\begin{equation*}
\fC, \fC_1 \subseteq \fR_\lambda \qquad \mbox{and} \qquad 
\fC \approx \fC_1 \approx \{x_1,\ldots,x_{d_\lambda}\},\tag{$*$}
\end{equation*}
where $\{x_1,\ldots, x_{d_\lambda}\}$ is our chosen left cell in 
$\fR_\lambda$. Indeed, since $\fX_{\fC}$ is irreducible, there exists 
some $\mu\in\Lambda$ such that $\fX_{\fC}$ is equivalent to $\fX_\mu$.
If we had $\lambda\neq \mu$, then Lemma~\ref{neun1} would imply $\fX_{\mu}
(\bC_z\bD_{d_z})=0$, contradicting the choice of $\lambda$. Thus, we have 
$\mu=\lambda$ and so $\chi_{\fC}=\chi_{\lambda}$. A similar argument shows 
that we also have  $\chi_{\fC_1}=\chi_{\lambda}$. But then 
Theorem~\ref{mainbichar} implies that $\fC, \fC_1 \subseteq \fR_\lambda$, 
and Theorem~\ref{eqcellsB} yields the second statement in~($*$).

Let $i,j\in \{1,\ldots,d_\lambda\}$ be such that $z \sim_{\cR} x_i$ and 
$d_z \sim_{\cR} x_j$. (These indices exist and are unique by 
Corollary~\ref{unique}.) Then Theorem~\ref{eqcellsB} and 
Lemma~\ref{neun2} imply that
\[ \bC_z\bD_{d_z}=\bC_{x_i}\bD_{x_j^{-1}}.\]
Similarly, if $k,l\in \{1,\ldots,d_\lambda\}$ are such that 
$w \sim_{\cR} x_k$ and $d_w \sim_{\cR} x_l$, then 
\[ \bC_w\bD_{d_w}=\bC_{x_k}\bD_{x_l^{-1}}.\]
Now, by Lemma~\ref{neun1}, the above elements are multiples of matrix
units with respect to the representation $\fX_\lambda$. (Recall that this
is the representation afforded by the left cell $\{x_1,\ldots,
x_{d_\lambda}\}$.) Hence the usual multiplication rules for matrix units 
imply that $j=k$ and
\[\bC_{z}\bD_{d_z} \cdot \bC_{w}\bD_{d_w} =\bC_{x_i}\bD_{x_j^{-1}} 
\cdot \bC_{x_j}\bD_{x_l^{-1}}=c_\lambda\,\bC_{x_i}\bD_{x_l^{-1}}.\]
Finally, let $d \in \cD_n$ be the unique element such that $d \sim_{\cR}
x_l$. Then, by Corollary~\ref{unique}, there is a unique element $u\in 
\fR_\lambda$ such that  $d=d_u$ and $u \sim_{\cR} x_i$. In particular,
this means that 
\[ u \sim_{\cR} x_i \sim_{\cR} z \quad \mbox{and}\quad u \sim_{\cL}
d_u=d_u^{-1} \sim_{\cL} x_l^{-1} \sim_{\cL} d_w \sim_{\cL} w.\]
Furthermore, the condition $j=k$ means that $w \sim_{\cR} x_k=x_j
\sim_{\cR} d_z=d_z^{-1} \sim_{\cR} z^{-1}$. Thus, if $\bC_z\bD_{d_z}\cdot
\bC_w\bD_{d_w}\neq 0$, we have established the desired multiplication rule.
Conversely, by following the above arguments backwards, one readily 
checks that $\bC_z\bD_{d_z}\cdot \bC_w\bD_{d_w}$ has the desired result if 
$w,z \in \fR_\lambda$ for some $\lambda \in \Lambda_n$, $w\sim_{\cR} z^{-1}$ 
and $z \sim_{\cR} u \sim_{\cL} w$ where $u \in \fR_\lambda$. Thus, the 
multiplication rule is proved.

Now let $x\in W_n$. Then we have 
\[\bC_{d_x}\bD_{x^{-1}}=\bC_w\bD_{d_w} \qquad \mbox{where $d_x\sim_{\cR} w$ 
and $x\sim_{\cR} d_w$}.\]
Hence, for any $z \in W_n$, we obtain
\[ \bC_z\bD_{d_z} \cdot \bC_{d_x}\bD_{x^{-1}}=\left\{\begin{array}{cl}
c_{\lambda_z}\,\bC_{u}\bD_{d_u}&\quad\mbox{if $z\sim_{\cR} u
\sim_{\cL}w\sim_{\cR}z^{-1}$},\\0&\quad\mbox{otherwise}. \end{array}\right.\]
Using the fact that we are dealing with a pair of dual basis, this yields
\[ \tau(\bC_z\bD_{d_z} \cdot \bC_{d_x}\bD_{x^{-1}})=\left\{\begin{array}{cl}
c_{\lambda_z}&\quad\mbox{if $z\sim_{\cR} u=d_u \sim_{\cL}w\sim_{\cR}z^{-1}$},
\\0&\quad\mbox{otherwise}. \end{array}\right.\]
Now, if the above condition on $u,z,w$ is satisfied, then we have 
$w \sim_{\cR} z^{-1}$ and $w\sim_{\cL} d_u=d_u^{-1} \sim_{\cL} z^{-1}$;
so we must have $w=z^{-1}$ by Corollary~\ref{unique}. But then we have
$x\sim_{\cL} d_x=d_x^{-1} \sim_{\cL} w^{-1}=z$ and $x^{-1}\sim_{\cL} d_w^{-1}
=d_w\sim_{\cL} w^{-1}=z$, which yields $x=z$. Thus, we have shown that
\[ \tau(\bC_z\bD_{d_z} \cdot \bC_{d_x}\bD_{x^{-1}})=\left\{\begin{array}{cl}
c_{\lambda_z}&\quad\mbox{if $x=z$}, \\0&\quad\mbox{otherwise}. 
\end{array}\right.\]
In order to prove the identity $\bC_w=\sum_{z\in W} h_{w,d_z,z}\,
c_{\lambda_z}^{-1}\,\bC_z\bD_{d_z}$, we just multiply both sides by $\bC_{d_x}
D_{x^{-1}}$ and note that, upon applying $\tau$, we obtain the same result. 
Once this identity is established, it follows that the elements $\{\bC_z
D_{d_z} \mid z \in W_n\}$  generate $\bH_{n,K}$. Since this generating 
set has the correct cardinality, it forms a basis.
\end{proof}

\begin{cor} \label{invert} The matrix $\bigl(h_{w,d_z,z}\bigr)_{w,z\in W_n}$ 
is invertible over $K$.
\end{cor}

\begin{proof} By Proposition~\ref{maxcor2},  the above matrix describes
the base change between two basis of $\bH_{n,K}$.
\end{proof}

\begin{defn} \label{maxcor3} In the above setting, we consider the
$\Z$-submodule 
\[ J_n:=\langle \hat{t}_w \mid w\in W_n\rangle_{\Z} \subseteq 
\bH_{n,K},\]
where we set $\hat{t}_w:=c_{\lambda_w}^{-1}\, \bC_w\bD_{d_w}$ for any 
$w\in W_n$. The multiplication rules in Corollary~\ref{maxcor2} immediately
imply that $J_n$ is a subring of $\bH_{n,K}$; indeed, we have 
\[ \hat{t}_x\hat{t}_y= \left\{\begin{array}{cl} \hat{t}_z & \quad \mbox{if 
$x \sim_{\cL} y^{-1}$ and $x \sim_{\cR} z \sim_{\cL} y$},\\ 0 & \quad 
\mbox{otherwise}. \end{array}\right.\]
Furthermore, we have a decomposition 
\[ J_n=\bigoplus_{\lambda \in \Lambda_n} J_{n,\lambda}\qquad
\mbox{(direct sum of two-sided ideals)},\]
where $J_{n,\lambda}:=\langle \hat{t}_w\mid w \in \fR_\lambda \rangle_{\Z}$
for every $\lambda \in \Lambda_n$. 
\end{defn}

We will see at the end of this section that $J_n$ actually is the ring $J$
introduced by Lusztig in \cite[Chap.~18]{Lusztig03}. However, our basis 
elements $\hat{t}_w$ will not correspond directly to Lusztig's basis 
elements. We have to perform a transformation of the following type.

Let $w \mapsto \hat{n}_w$ be an integer-valued function on $W_n$
satisfying the following two properties:
\begin{itemize}
\item[{\bf (N1)}] we have $\hat{n}_w=\pm 1$ for all $w \in W_n$;
\item[{\bf (N2)}] the function $w \mapsto \hat{n}_w$ is constant on right 
cells.
\end{itemize}
Having fixed a function as above, we set $t_w:=\hat{n}_w \hat{t}_w$ for 
all $w\in W_n$. By {\bf (N1)}, the elements $\{t_w\mid w \in W_n\}$ form a 
new $\Z$-basis of the ring $J_n$. Writing
\[ t_x\, t_y=\sum_{z\in W_n}\hat{\gamma}_{x,y,z^{-1}} \, t_z 
\qquad (x,y\in W_n),\]
the structure constants $\hat{\gamma}_{x,y,z^{-1}}$ are given by:
\begin{itemize}
\item[{\bf (N3)}] $\qquad \hat{\gamma}_{x,y,z^{-1}}=
\left\{\begin{array}{cl} \hat{n}_y=\pm 1 & \quad \mbox{if $x \sim_{\cL} 
y^{-1}$ and $x \sim_{\cR} z \sim_{\cL} y$},\\ 0 & \quad 
\mbox{otherwise}; \end{array}\right.$
\end{itemize}
see Proposition~\ref{maxcor2}. The following results will all be 
formulated in terms of the basis $\{t_w\mid w\in W_n\}$ of $J_n$, where we 
assume throughout that a function satisfying {\bf (N1)}, {\bf (N2)} has been 
fixed. An obvious example is given by the function $\hat{n}_w=1$ for
all $w \in W_n$. (As we will see at the end of this section, we have 
to take a different function in order to identify $t_w$ with the 
corresponding element in Lusztig's construction.)

\begin{cor} \label{maxcor4} The ring $J_n$ introduced above has
unit element
\[ T_1=\sum_{z \in \cD_n} \hat{n}_z\,t_z \qquad \mbox{(where 
$T_1$ is the unit element in $\bH_n$)}.\] 
For every $\lambda \in \Lambda_n$, we have $J_{n,\lambda} \cong 
M_{d_\lambda}(\Z)$. 
\end{cor}

\begin{proof} Let $z  \in \cD_n$ and assume that $z\in \fR_\lambda$.
Since $d_z=z$ is an  involution, the argument in the proof of 
Proposition~\ref{maxcor2} now shows that 
\[ \bC_{z}\bD_{z}=\bC_{x_i}\bD_{x_i^{-1}} \quad \mbox{for some $1\leq i
\leq d_\lambda$},\]
where $\{x_1,\ldots,x_{d_\lambda}\} \subseteq \fR_\lambda$ is our
chosen left cell. Furthermore, we have
\[ \sum_{z \in \cD_n \cap \fR_\lambda} \bC_{z}\bD_{z}=\sum_{i=1}^{d_\lambda}
\bC_{x_i}\bD_{x_i^{-1}}.\]
By Lemma~\ref{neun1}, the image of the above element under $\fX_\mu$ 
is $0$ if $\lambda\neq \mu$, and $c_\lambda$ times the identity matrix
if $\lambda=\mu$. 

Hence we conclude that the image of $\varepsilon:=\sum_{z\in \cD_n} 
\hat{t}_z\in J_n$ under $\fX_\lambda$ (for any $\lambda$) is the 
identity matrix. Since $\bH_{n,K}$ is split semisimple, this implies 
that $\varepsilon$ is the identity element in $\bH_{n,K}$, that is, 
we have $\varepsilon=T_1$.
\end{proof}

We can now establish the following result which is in complete analogy to 
Lusztig \cite[Theorem~18.9]{Lusztig03}.

\begin{cor} \label{maxcor5} Let $J_{n,A}=A\otimes_{\Z} J_n=\langle t_w
\mid w\in W_n\rangle_A\subseteq \bH_{n,K}$. The $A$-linear map ${\phi} 
\colon \bH_n \rightarrow J_{n,A}$ defined by 
\[ {\phi}(\bC_w^{\, \delta})=\sum_{z\in W_n} h_{w,d_z,z}\,\hat{n}_z\,t_z
\qquad (w\in W_n)\]
is a homomorphism of $A$-algebras respecting the unit elements.
\end{cor}

\begin{proof} Since $J_{n,A}\hookrightarrow \bH_{n,K}$, the above
formula actually defines a $K$-linear map $\phi_K \colon \bH_{n,K}
\rightarrow \bH_{n,K}$ whose restriction to $\bH_n$ is $\phi$.
By Proposition~\ref{maxcor2}, the set $\{\hat{t}_w\mid w\in W_n\}$ is a 
basis of $\bH_{n,K}$. Furthermore, the formula in that proposition shows 
that ${\phi_K}(\hat{t}_w^{\,\delta})=\hat{t}_w$ for all $w\in W_n$. 
Thus, we have ${\phi_K} \circ \delta=\mbox{id}$ on $\bH_{n,K}$ and, 
consequently, ${\phi_K}$ is a $K$-algebra homomorphism respecting 
the unit elements.
\end{proof}

The next result can be regarded as a weak version of property (P15) in
Lusztig's list of conjectures in \cite[Chap.~14]{Lusztig03}. 

\begin{prop}[Compare Lusztig \protect{\cite[18.9(b)]{Lusztig03}}] 
\label{weakP15} Let $x,x',y,z,w\in W_n$ be such that $y \sim_{\cL} x' 
\sim_{\cR} x^{-1}$. Then we have 
\[\sum_{z \in W_n} h_{w,z,y}\,\hat{\gamma}_{x,x',z^{-1}}=
\sum_{z\in W_n} h_{w,x,z}\, \hat{\gamma}_{z,x',y^{-1}}.\]
\end{prop}

\begin{proof} By assumption, we have $x\sim_{\cL} x'^{-1}$ and 
$x'\sim_{\cL} y$. Hence we have $x,x',y \in \fR_\lambda$ for some
$\lambda \in \fR_n$. So, by Corollary~\ref{unique}, there exist unique
elements $z_0,z_1\in \fR_\lambda$ such that 
\[ x\sim_{\cR} z_0\sim_{\cL} x' \qquad \mbox{and}\qquad  x'^{-1} 
\sim_{\cL} z_1 \sim_{\cR} y.\]  
Now Proposition~\ref{maxcor2} shows that $\hat{\gamma}_{x,x',z^{-1}}=0$ 
unless $z=z_0$, in which case the result is $\hat{\gamma}_{x,x',z_0^{-1}}=
\hat{n}_{x'}$. Similarly, we have $\hat{\gamma}_{z,x', y^{-1}}=0$ unless 
$z=z_1$, in which case the result is $\hat{\gamma}_{z_1,x',y^{-1}}=
\hat{n}_{x'}$. Hence the desired equality is equivalent to the identity 
\begin{equation*}
h_{w,z_0,y}=h_{w,x,z_1}.\tag{$*$}
\end{equation*}
Suppose that $h_{w,z_0,y}\neq 0$. Then $y \leq_{\cL} z_0$. We conclude that
\[ x' \sim_{\cL} y \leq_{\cL} z_0 \sim_{\cL} x'\]
and so $y \sim_{\cL}z_0$. On the other hand, we have $z_0 \sim_{\cR} x$
and, hence,
\[ z_1 \sim_{\cL} x'^{-1}\sim_{\cL} x \sim_{\cL} z_0.\]
So we can apply Theorem~\ref{eqcellsB} and this yields $h_{w,z_0,y}=
h_{w,x,z_1}$. Thus, ($*$) holds in this case. Conversely, if $h_{w,x,z_1}
\neq 0$, then a similar argument shows that, again, ($*$) holds. Finally, 
this also yields that $h_{w,z_0,y}=0$ if and only if $h_{w,x,z_1}=0$.
Thus, ($*$) holds in all cases.
\end{proof}

Let $\cM$ be the free $A$-module with basis $\{\varepsilon_x \mid x\in 
W_n\}$. Identifying $\bH_n$ and $\cM$ via $\bC_w\mapsto\varepsilon_w$, the 
obvious $\bH_n$-module on $\bH_n$ (given by left multiplication) becomes
the $\bH_{n}$-module on $\cM$ given by 
\[ \bC_w.\varepsilon_x=\sum_{y \in W_n} h_{w,x,y}\, \varepsilon_y \qquad
\qquad  (w,x \in W_n).\]
On the other hand, we can also identify $\cM$ with $J_{n,A}$, via 
$\varepsilon_w \mapsto \hat{n}_w t_w$. Then the obvious $J_{n,A}$-module
structure on $J_{n,A}$ (given by left multiplication) becomes the 
$J_{n,A}$-module structure on $\cM$ given by 
\[ t_w \ast\varepsilon_x=\sum_{y \in W_n} \hat{\gamma}_{w,x,y^{-1}}\, 
\hat{n}_x \hat{n}_y \, \varepsilon_y\qquad\qquad (w,x \in W_n).\]
Now we can state the following result.

\begin{cor}[Compare Lusztig \protect{\cite[18.10]{Lusztig03}}]
\label{l18-10} For any $h \in \bH_n$ and any $x\in W_n$, the difference 
$h.\varepsilon_x-\phi(h^\delta)\star\varepsilon_x$ is an $A$-linear 
combination of elements $\varepsilon_y$ where $y<_{\cLR} x$ (that is,
we have $y\leq_{\cLR} x$ but $y \not\sim_{\cLR} x$).
\end{cor}

\begin{proof} It is enough to prove this for $h=\bC_w$ where $w\in W_n$. 
Then we have 
\begin{align*}
\phi(\bC_w^\delta)\star\varepsilon_x & =
\sum_{z \in W_n} h_{w,d_z,z}\,\hat{n}_z \, t_z \star\varepsilon_x\\
&=\sum_{z \in W_n} \sum_{y \in W_n} h_{w,d_z,z}\,\hat{\gamma}_{z,x,y^{-1}} 
\, \hat{n}_z\hat{n}_x \hat{n}_y \, \varepsilon_y.
\end{align*}
Now, if the term correponding to $y,z$ is non-zero, then we have
$\gamma_{z,x,y^{-1}}\neq0$ and so $z \sim_{\cR} y$. Hence we also
have $\hat{n}_z=\hat{n}_y$ and so $\hat{n}_z\hat{n}_y=1$. By a similar
argument, we can also assume that $x \sim_{\cL} y$ and $x^{-1} \sim_{\cL}
z \sim_{\cL} d_z$. In particular, we have $d_z=d_{x^{-1}}$. Consequently, 
we can rewrite the above sum as follows: 
\begin{align*}
\phi(\bC_w^\delta)\star\varepsilon_x &=
\sum_{\atop{y \in W_n}{x \sim_{\cL} y}}\Bigl(\sum_{z \in W_n} 
h_{w,d_{x^{-1}},z}\, \hat{\gamma}_{z,x,y^{-1}} \Bigr) \, \hat{n}_x \, 
\varepsilon_y\\ &=\sum_{\atop{y \in W_n}{x \sim_{\cL} y}} 
\Bigl(\sum_{z \in W_n} h_{w,z,y}\, \hat{\gamma}_{d_{x^{-1}},x,z^{-1}} 
\Bigr) \, \hat{n}_x \, \varepsilon_y
\end{align*}
where the second equality holds by Proposition~\ref{weakP15}. Now, by
{\bf (N3)}, we have $\gamma_{d_{x^{-1}},x,z^{-1}}=0$ unless $x=z$ in which 
case the result equals $\hat{n}_x$. Hence the above sum reduces to:
\[ \phi(\bC_w^\delta)\star\varepsilon_x =
\sum_{\atop{y \in W_n}{x \sim_{\cL} y}} h_{w,x,y}\, \varepsilon_y.\]
On the other hand, we know that ($\spadesuit$) holds by Theorem~\ref{relsn}.
So, for any $y'\in W_n$, we have $h_{w,x,y'}=0$ unless $y' \sim_{\cL} x$ or 
$y' <_{\cLR} x$. Hence we see that, indeed, the difference $h.\varepsilon-
\phi(\bC_w^\delta)\star \varepsilon_x$ has the required form.
\end{proof}

We now apply the above results to construct a ``canonical'' algebra 
isomorphism from $\bH_{n,K}$ onto $KW_n$, the group algebra of $W_n$ 
over $K$. Let $R={\Q}[\Gamma]={\Q} \otimes_{\Z} A$ and set $\bH_{n,R}=
R \otimes_A \bH_n$, $J_{n,R}:=R \otimes_A J_A$. The previously defined 
modules structures of $\bH_n$ and $J_{n,A}$ on $\cM$ naturally extend to 
module structures of $\bH_{n,R}$ and $J_{n,R}$, respectively, on $\cM_R=
R \otimes_A \cM$. Now we also describe an $RW_n$-module structure on 
$\cM_R=R \otimes_A \cM$, as follows. We have a ring homomorphism 
\[ \theta \colon R\rightarrow R, \qquad e^\gamma \mapsto 1 \qquad
(\gamma \in \Gamma).\]
We can regard $R$ as an $R$-module via $\theta$; then we obtain 
$R\otimes_R \bH_n=RW_n$. We denote $c_w=1 \otimes \bC_w \in RW_n$ for 
any $w \in W_n$. Hence, we may also regard $\cM_R$ as an $RW_n$-module, 
where $c_w$ ($w\in W_n)$ acts by
\[ c_w \diamond \varepsilon_x=\sum_{y \in W_n} \theta(h_{w,x,y})\, 
\varepsilon_y \qquad \mbox{for any $x \in W_n$}.\]
Note that this $RW_n$-module structure on $\cM_R$ coincides with 
the obvious structure (given by left multiplication), where we identity 
$RW_n$ and $\cM_R$ via $c_w \mapsto \varepsilon_w$. 

\begin{thm}[See Lusztig \protect{\cite[Theorem~3.1]{Lu0}} in the case
of equal parameters] \label{caniso} There is a unique homomorphism of
$R$-algebras $\Phi \colon \bH_{n,R} \rightarrow RW_n$ such that, for 
any $h \in \bH_{n,R}$ and any $x \in W_n$, the difference $h.\varepsilon_x
-\Phi(h)\diamond\varepsilon_x$ is a linear combination of elements 
$\varepsilon_y$ with $y <_{\cLR} x$. Furthermore, writing
\[ \Phi(\bC_w)=\sum_{z \in W_n} \Phi_{w,z} \, z \qquad \mbox{where
$\Phi_{w,z} \in R$},\]
we have $\Phi_{w,z}=\overline{\Phi}_{w,z}$ and $\theta(\Phi_{w,z})=
\delta_{wz}$ for all $w,z\in W_n$. Finally, the induced map $\Phi_K 
\colon \bH_{n,K} \rightarrow KW_n$ is an isomorphism.
\end{thm}

(Here, $\delta_{wz}$ denotes the Kronecker delta, and $r\mapsto \bar{r}$ 
is the ring involution such that $e^\gamma \mapsto e^{-\gamma}$ for 
all $\gamma\in \Gamma$).

\begin{proof} First we show the uniqueness statement. Let $\Phi_i\colon 
\bH_{n,R}\rightarrow RW_n$ ($i=1,2$) be two homomorphisms such that, 
for any $h \in \bH_{n,R}$ and any $x\in W_n$, the difference 
$h.\varepsilon_x - \Phi_i(h)\diamond \varepsilon_x$ is a linear 
combination of elements $\varepsilon_y$ with $y <_{\cLR} x$. Then the 
difference $(\Phi_1(h)- \Phi_2(h))\diamond \varepsilon_x$ is a linear 
combination of elements $\varepsilon_y$ with $y<_{\cLR} y$. Consequently,
$\Phi_1(h)-\Phi_2(h)\in RW_n \subseteq KW_n$ acts as a nilpotent
operator on $\cM_K=K \otimes_{R} \cM$.  But, as we already noted above,
$\cM_K$ is the left regular $KW_n$-module, hence we must have $\Phi_1(h)-
\Phi_2(h)=0$.  

So it remains to show that an $R$-algebra homomorphism $\Phi$ with the 
required properties does exist. In Corollary~\ref{maxcor5}, we extend 
scalars from $A$ to $R$ and obtain a homomorphism of $R$-algebras
\[ \alpha \colon \bH_{n,R} \rightarrow J_R=R \otimes_A J_A, \qquad
\bC_w \mapsto \phi(\bC_w^\delta).\]
Explicitly, $\alpha$ is given by the formula
\[ \alpha(\bC_w)=\sum_{z \in W_n} h_{w,d_z,z}\,\hat{n}_z\, t_z 
\qquad \mbox{for any $w\in W_n$}.\]
(We can take $\hat{n}_z=1$ for all $z\in W_n$.) By Corollary~\ref{l18-10}, 
the above homomorphism has the property that, for any $h \in \bH_{n,R}$ and 
any $x \in W_n$, the difference $h.\varepsilon_x-\alpha(h)\star
\varepsilon_x$ is an $R$-linear combination of elements $\varepsilon_y$ 
where $y<_{\cLR} x$.

Now, as before, we regard $R$ as an $R$-module via $\theta$ and extend
scalars. Since the structure constants of $J_n$ with respect to the 
basis $\{t_w\}$ are integers, they are not affected by $\theta$. Hence we 
obtain an induced homomorphism of $R$-algebras 
\[ \beta\colon RW_n\rightarrow J_R\]
such that 
\[ \beta(c_w)=\sum_{z \in W_n} \theta(h_{w,d_z,z})\,\hat{n}_z\, t_z 
\qquad \mbox{for any $w\in W_n$}.\]
Now the identity in Proposition~\ref{maxcor2} ``specializes'' to an analogous 
identity in $RW_n$. (Note that $\theta(c_\lambda)=|W_n|/d_\lambda \neq 0$ 
for each $\lambda \in \Lambda_n$; see \cite[\S 8.1]{ourbuch}.) We deduce from 
this that the matrix
\[ \bigl(\theta(h_{w,d_z,z})\bigr)_{w,z\in W_n}\]
is invertible over $K$. Since the coefficients of that matrix lie in $\Q$,
so do the coefficients of its inverse. Consequently, $\beta$ is an 
isomorphism of $R$-algebras. Furthermore, a computation analogous to that
in the proof of Corollary~\ref{l18-10} shows that we have 
\[\beta(c_w) \star \varepsilon_x=\sum_{\atop{y \in  W_n}{x \sim_{\cL} y}}
\theta(h_{w,x,y})\, \varepsilon_y \qquad \mbox{for any $x,w \in W_n$},\]
and that $\beta(c_w)\star \varepsilon_x-c_w \diamond \varepsilon_x$ is an 
$R$-linear combination of elements $\varepsilon_y$ where $y<_{\cLR} x$. 
Consequently, since $\beta$ is an isomorphism, we also have that, for any 
$\iota\in J_{n,R}$ and any $x \in W_n$, the difference $\iota\star 
\varepsilon_x -\beta^{-1}(\iota) \diamond \varepsilon_x$ is an $R$-linear 
combination of elements $\varepsilon_y$ where $y<_{\cLR} x$. Now we set  
\[ \Phi:=\beta^{-1} \circ \alpha \colon \bH_{n,R} \rightarrow RW_n.\]
Let $h \in \bH_{n,R}$ and $x \in W_n$. Setting $\iota:=\alpha(h) \in 
J_{n,R}$, we obtain that  
\begin{align*}
h.\varepsilon_x-\Phi(h)\diamond \varepsilon_x&=
h.\varepsilon_x-\alpha(h) \star \varepsilon_x+
\alpha(h)\star \varepsilon_x- \Phi(h)\diamond \varepsilon_x\\& = 
(h.\varepsilon_x-\alpha(h) \star \varepsilon_x)+ (\iota\star \varepsilon_x- 
\beta^{-1}(\iota)\diamond \varepsilon_x)
\end{align*}
is an $R$-linear combination of elements $\varepsilon_y$ where $y<_{\cLR} x$, 
as required. 

Finally, $\phi_K$ is an isomorphism since $\alpha_K$ is invertible (see 
Corollary~\ref{invert}) and $\beta$ is an isomorphism.  Furthermore, 
the coefficients $\Phi_{w,z}$ have the  stated properties, since $\Phi$ 
is defined as the composition of $\alpha$ (whose matrix is given by the 
coefficients $h_{w,d_z,z}$) and the inverse of $\beta$ (whose matrix is 
given by the inverse of the matrix with coefficients $\theta(h_{w,d_z,z})$).
\end{proof}

Note that the above proof relies on the existence of the homomorphism
$\phi \colon \bH_n \rightarrow J_n$ and Corollaries~\ref{invert}, 
\ref{l18-10}. We could not follow Lusztig's original proof in \cite{Lu0} 
since, in the present case, the constants $M_{y,w}^s$ appearing in the
multiplication formula for the Kazhdan--Lusztig basis are no longer integers.

\begin{exmp} \label{expB2} Let us consider the case $n=2$, where $W_2=
\langle t,s_1\rangle$ is the dihedral group of order~$8$. We set $s_0=t$. 
The coefficients $h_{s,y,z}$ (for $s=s_0,s_1$) and the left cells have 
already been determined by an explicit computation in \cite[\S 6]{Lusztig83}.
The left cells are
\[ \{1\},\quad \{s_1\}, \quad \{s_0,s_1s_0\}, \quad \{s_1s_0s_1,s_0s_1\},
\quad \{s_0s_1s_0\}, \quad \{w_0\}\]
where $w_0=s_1s_0s_1s_0$ is the unique element of maximal length. For each 
left cell, the first element listed is the unique element from $\cD_2$ in 
that left cell. From the information in \cite[\S 6]{Lusztig83}, we know
$h_{w,d_z,z}$ for $w \in \{s_0, s_1\}$. This yields the following formulas 
for the homomorphism $\phi \colon \bH_2 \rightarrow J_{2,A}$:
\begin{align*}
\phi(\bC_{s_0}) &= (Q{+}Q^{-1})t_{s_0}+ (Qq^{-1}{+}Q^{-1}q)t_{s_0s_1}\\
& \qquad \qquad +(Q{+}Q^{-1})t_{s_0s_1s_0}+(Q{+}Q^{-1})t_{w_0},\\
\phi(\bC_{s_1}) &= (q{+}q^{-1})t_{s_1}+t_{s_1s_0}+(q{+}q^{-1})t_{s_1s_0s_1}+
(q{+}q^{-1})t_{w_0},
\end{align*}
where we take the function $\hat{n}_w=1$ for all $w \in W_2$. Using the
multiplication formula for the Kazhdan--Lusztig basis, we can deduce 
explicit expressions of $\phi(\bC_w)$, for any $w \in W_2$; this yields the 
whole matrix of coefficients $(h_{w,d_z,z})_{w,z}$. In order to construct 
$\Phi \colon \bH_{2,R} \rightarrow RW_2$, we follow the proof of 
Theorem~\ref{caniso}. First, we apply the ring homomorphism $\theta 
\colon R \rightarrow R$, that is, we specialise $Q,q\mapsto 1$. The matrix 
of all coefficients $\theta (h_{w,d_z,z})$ is given in Table~\ref{tab1}. 
Composing the matrix of $\phi$ with the inverse of the matrix in
Table~\ref{tab1} and expressing the basis $\{c_w\}$ of $RW_2$ in terms 
of the standard basis consisting of group elements, we obtain the following 
explicit description of the homomorphism $\Phi \colon \bH_{2,R}
\rightarrow RW_2$:
\begin{align*}
\Phi(T_{s_0}) &= \frac{1}{2}(Q-Q^{-1}) \cdot 1+\frac{1}{2}(Q+Q^{-1})
\cdot s_0\\&\quad+\frac{1}{4}(Q{-}Qq^{-1}{-}Q^{-1}q{+}Q^{-1})\cdot 
(-s_1+s_1s_0-s_0s_1+s_0s_1s_0),\\
\Phi(T_{s_1}) &= \frac{1}{2}(q-q^{-1}) \cdot 1+\frac{1}{2}(q+q^{-1})
\cdot s_1\\&\qquad\qquad +\frac{1}{4}(q-2+q^{-1})\cdot 
(-s_0-s_1s_0+s_0s_1+s_1s_0s_1).
\end{align*}
Note that the formulas do not make any reference to the Kazhdan--Lusztig
basis. Further note that the above formulas specialise to $T_{s_0} \mapsto 
s_0$, $T_{s_1} \mapsto s_1$ when we set $Q,q \mapsto 1$. Also, if we 
consider the images of $\bC_{s_0}=T_{s_0}+Q^{-1}T_1$ and $\bC_{s_1}=T_{s_1}+
q^{-1}T_1$, then all the coefficients are seen to be fixed by the 
involution $r \mapsto \bar{r}$, as stated in Theorem~\ref{caniso}.
\end{exmp}

\begin{table}[htbp]
\caption{The coefficients $\theta(h_{w,d_z,z})$ in type $B_2$} 
\label{tab1} \begin{center}
$\renewcommand{\arraystretch}{1.1} 
\begin{array}{|c|c|c|c@{\hspace{1mm}}c|c@{\hspace{1mm}}c|c|c|}\hline  
\theta(h_{w,d_z,z}) &1&s_1&s_0&s_1s_0&s_1s_0s_1&s_0s_1&s_0s_1s_0& w_0\\\hline 
1  &   1&   1&   1&   0&   1&   0&   1&   1 \\ \hline
s_1  &   0&   2&   0&   1&   2&   0&   0&   2 \\ \hline
s_0  &   0&   0&   2&   0&   0&   2&   2&   2 \\
s_1s_0  &   0&   0&   0&   2&   2&   0&   0&   4 \\ \hline
s_1s_0s_1  &   0&   0&   0&   2&   4&   0&   0&   8 \\
s_0s_1  &   0&   0&   2&   0&   0&   4&   0&   4 \\ \hline
s_0s_1s_0  &   0&   0&   0&   0&   0&   0&  -4&   4 \\ \hline
w_0  &   0&   0&   0&   0&   0&   0&   0&   8 \\
\hline \end{array}$
\end{center}
\end{table}

To close this section, we explain how to identify $J_n$ with Lusztig's 
ring $J$. First we need some definitions.

Let $z\in W_n$. Following Lusztig \cite[14.1]{Lusztig03}, we define an 
element $\bDt_n(z) \in \Gamma$ and an integer $0\neq n_z\in \Z$ 
by the condition
\[ e^{\bDt_n(z)}\,P_{1,z}^* \equiv n_z \bmod A_{<0};\]
note that $\bDt_n(z)\geq 0$. Following \cite[13.6]{Lusztig03}, we define a 
function $\ba_n \colon W_n \rightarrow \Gamma$ as follows. Let $z\in W_n$.  
Then we set
\[ \ba_n(z):=\min \{\gamma \geq 0 \mid e^\gamma \, h_{x,y,z} \in 
A_{\geq 0} \mbox{ for all $x,y \in W_n$}\}.\]
Furthermore, for any $x,y,z\in W_n$, we set 
\[ \gamma_{x,y,z^{-1}}=\mbox{ constant term of $e^{\ba_n(z)}\,h_{x,y,z}
\in A_{\geq 0}$}.\]
Following \cite[Chap.~18]{Lusztig03}, we use the constants $\gamma_{x,y,z}$ 
to define a new bilinear pairing on our free abelian group $J_n$ with basis 
$\{t_w\mid w\in W_n\}$ by 
\[ t_x \bullet t_y:=\sum_{z\in W} \gamma_{x,y,z^{-1}} \, t_z\qquad
\mbox{for all $x,y \in W$}.\]
As explained in \cite[Chap.~18]{Lusztig03}, one can show that $(J_n,
\bullet)$ is an associative ring with unit element $1_J=\sum_{d \in \cD_n} 
n_d t_d$, if the following properties from Lusztig's list of conjectures in
\cite[Chap.~14]{Lusztig03} hold, where 
\[ \cD:=\{z\in W_n \mid \ba_n(z)=\bDt_n(z)\}.\]
\begin{itemize}
\item[\bf (P1)] For any $z\in W_n$ we have $\ba_n(z)\leq \bDt_n(z)$.
\item[\bf (P2)] If $d \in \cD_n$ and $x,y\in W_n$ satisfy $\gamma_{x,y,d}
\neq 0$, then $x=y^{-1}$.
\item[\bf (P3)] If $y\in W_n$, there exists a unique $d\in \cD_n$ such that
$\gamma_{y^{-1},y,d}\neq 0$.
\item[\bf (P4)] If $z'\leq_{\cL\cR} z$ then $\ba_n(z')\geq \ba_n(z)$. Hence, 
if $z'\sim_{\cLR} z$, then $\ba_n(z)=\ba_n(z')$.
\item[\bf (P5)] If $d\in \cD_n$, $y\in W_n$, $\gamma_{y^{-1},y,d}\neq 0$, then
$\gamma_{y^{-1},y,d}=n_d=\pm 1$.
\item[\bf (P6)] If $d\in \cD_n$, then $d^2=1$.
\item[\bf (P7)] For any $x,y,z\in W_n$, we have $\gamma_{x,y,z}=
\gamma_{y,z,x}$.
\item[\bf (P8)] Let $x,y,z\in W_n$ be such that $\gamma_{x,y,z}\neq 0$. Then
$x\sim_{\cL} y^{-1}$, $y \sim_{\cL} z^{-1}$, $z\sim_{\cL} x^{-1}$.
\end{itemize}
Now we have the following result:

\begin{thm}[Geck--Iancu \protect{\cite[Theorem~1.3]{GeIa05}}] \label{Pconj}
In the ``asymptotic case'' in type $B_n$, all the properties {\bf (P1)--(P15)}
from Lusztig's list of conjectures in \cite[Chap.~14]{Lusztig03} hold, 
except possibly {\bf (P9)}, {\bf (P10)} and {\bf (P15)}. Furthermore, we 
have $\cD=\cD_n:= \{z\in W_n \mid z^2=1\}$.
\end{thm}

In particular, {\bf (P1)--(P8)} hold and so we do have an associative ring
$(J_n,\bullet$) with unit element.

Following \cite[18.8]{Lusztig03}, we set $\hat{n}_w:=n_d$, where $d 
\in\cD_n$ is the unique element such that $d \sim_{\cL} w^{-1}$ (which
exists and is unique by Corollary~\ref{unique}). By construction, the 
function $w \mapsto \hat{n}_w$ is constant on the right cells of $W_n$.  
Thus, $w \mapsto \hat{n}_w$ is an integer-valued function on $W_n$ which
satisfies {\bf (N2)}. By {\bf (P5)}, we see that {\bf (N1)} also holds.

We shall now take this function in the above discussion. That is, the
formula {\bf (N3)} reads:
\[\hat{\gamma}_{x,y,z^{-1}}=\left\{\begin{array}{cl} n_d & \quad 
\mbox{if $x \sim_{\cL} y^{-1}$, $y \sim_{\cL} z$, $z \sim_{\cR} x$,
$d=d^{-1}\sim_{\cL} y^{-1}$},\\ 0 & \quad \mbox{otherwise}; 
\end{array}\right.\]
Now we can state the following result.

\begin{prop} \label{finalid} For any $x,y,z\in W_n$, we have 
$\hat{\gamma}_{x,y,z^{-1}}=\gamma_{x,y,z^{-1}}$. 
\end{prop}

\begin{proof} Let $x,y,z\in W_n$ be such that $x\sim_{\cL} y^{-1}$, 
$y \sim_{\cL} z$ and $z^{-1}\sim_{\cL} x^{-1}$.  By {\bf (P3)}, there
exists a unique $d \in \cD_n$ such that $\gamma_{x^{-1},x,d}\neq 0$.
We have $d^2=1$. Hence, by {\bf (P8)}, we obtain $d \sim_{\cL} x\sim_{\cL} 
y^{-1}$ and so 
\[ \hat{\gamma}_{x,y,z^{-1}}=\hat{n}_y=n_d;\]
see the formula in Definition~\ref{maxcor3}. This yields the identity
\[ \hat{\gamma}_{x,y,z^{-1}}=n_d=\gamma_{x^{-1},x,d}=\gamma_{x,d,x^{-1}};
\qquad\mbox{see {\bf (P5)}, {\bf (P7)}}.\]
Now, since $d=d^{-1}\sim_{\cR} y$, $x \sim_{\cR} z$ and $y \sim_{\cL} z$, 
we have
\[h_{x,d,x}=h_{x,y,z}; \qquad \mbox{see Theorem~\ref{eqcellsB}}.\] 
By {\bf (P4)}, we have $\ba_n(x)=\ba_n(z)$. Hence the above identity implies
that $\ba_n(x)h_{x,d,x}$ and $\ba(z)_nh_{x,y,z}$ have the same constant term 
and so
\[ \hat{\gamma}_{x,y,z^{-1}}=\gamma_{x,d,x^{-1}}=\gamma_{x,y,z^{-1}},\]
as required. It remains to consider the case where $x,y,z$ do not satisfy
the conditions  $x\sim_{\cL} y^{-1}$, $y \sim_{\cL} z$, $z^{-1}\sim_{\cL}
x^{-1}$. But then we have $\hat{\gamma}_{x,y,z^{-1}}=0$ by 
Proposition~\ref{maxcor2} and $\gamma_{x,y,z^{-1}}=0$ by {\bf (P8)}.
Hence we have $\gamma_{x,y,z^{-1}}=\hat{\gamma}_{x,y,z^{-1}}$ in all cases.
\end{proof}

Combining Theorem~\ref{Pconj} with the results in this paper, we can 
summarize the situation as follows.

\begin{cor} \label{finalcor} In the ``asymptotic case'' in type $B_n$, the
properties {\bf (P1)--(P14)} from Lusztig's list in 
\cite[Chap.~14]{Lusztig03} hold. Furthermore, we have the weak version of
{\bf (P15)} in Proposition~\ref{weakP15}. The ring $J_n$ with its ring 
structure given by $J_n \subseteq \bH_{n,K}$ as in Definition~\ref{maxcor3} 
is Lusztig's ring $(J_n,\bullet)$.
\end{cor}

\begin{proof} The statement concerning $J_n$ follows from 
Proposition~\ref{finalid}. Furthermore, taking into account 
Theorem~\ref{Pconj}, it only remains to consider:
\begin{gather*}
x \leq_{\cL} y \quad \mbox{and} \quad \ba_n(x)=\ba_n(y)\quad \Rightarrow
\quad x\sim_{\cL} y,\tag{\bf P9}\\
x \leq_{\cR} y \quad \mbox{and} \quad \ba_n(x)=\ba_n(y)\quad \Rightarrow
\quad x\sim_{\cR} y.\tag{\bf P10}
\end{gather*}
Now, by \cite[14.10]{Lusztig03}, property {\bf (P10)} is a formal consequence
of {\bf (P9)}. To prove {\bf (P9)}, let $x,y\in W_n$ be such that 
$x\leq_{\cL} y $ and $\ba_n(x)=\ba_n(y)$. By Theorem~\ref{Pconj}, we know
that the following holds:
\begin{equation*}
x \leq_{\cLR} y \quad \mbox{and} \quad \ba_n(x)=\ba_n(y)\quad \Rightarrow
\quad x\sim_{\cLR} y.\tag{\bf P11}
\end{equation*}
Hence we conclude that $x\sim_{\cLR} y$, and ($\spadesuit$)
yields $x\sim_{\cL} y$, as desired. 
\end{proof}

\begin{rem} \label{luform} The defining formula for $\phi$ in
\cite[18.9]{Lusztig03} reads
\[ \phi(\bC_w^{\, \delta})=\sum_{\atop{z \in W_n,d \in \cD_n}{\ba_n(z)=
\ba_n(d)}} h_{w,d,z}\, \hat{n}_z\, t_z.\]
But, once {\bf (P9)} is known to hold, the above formula reduces to the 
one in Corollary~\ref{maxcor5}. Indeed, assume that $z \in W_n$ and $d 
\in \cD_n$ are such that $\ba_n(z)=\ba_n(d)$ and $h_{w,d,z}\neq 0$. Then 
$z \leq_{\cL} d$ and {\bf (P9)} implies that $z \sim_{\cL} d$. Thus,
we must have $d=d_z$.
\end{rem}

\begin{rem} \label{modular} In \cite{my98}, it is shown that the
existence of Lusztig's homomorphism into the ring $J$ has various
applications in the modular representation theory of Iwahori--Hecke
algebras, most notably the fact that there is a natural ``unitriangular'' 
structure on the decomposition matrix associated with a non-semisimple
specialisation. Given Corollary~\ref{finalcor} and the results in this paper,
we can now apply the theory developed in \cite{my98} to $\bH_n$ as well. 
This should lead to new proofs of some results by Dipper--James--Murphy
\cite{DiJa95}.
\end{rem}

\medskip
\noindent {\bf Acknowledgements.} This paper was written while the
author enjoyed the hospitality of the Bernoulli center at the EPFL 
Lausanne (Switzerland), in the framework of the research program 
``Group representation theory'' from january to june 2005.



\begin{thebibliography}{131}

\bibitem{Ar}
{\sc S.~Ariki}, {Robinson--Schensted correspondence and left cells}.
Combinatorial methods in representation theory (Kyoto, 1998), 1--20,
Adv. Stud. Pure Math., {\bf 28}, Kinokuniya, Tokyo, 2000.

\bibitem{BV}
{\sc D.~Barbasch and D.~Vogan}, Primitive ideals and orbital integrals in
complex exceptional groups, J. Algebra {\bf 80} (1983), 350--382.

\bibitem{BI}
{\sc C.~Bonnaf\'e and L.~Iancu}, Left cells in type $B_n$ with unequal
parameters, Represent. Theory {\bf 7} (2003), 587--609.

\bibitem{BI2}
{\sc C.~Bonnaf\'e},  Tow-sided cells in type $B$ in the asymptotic case,
preprint available at  {\tt http://arXiv.org/math.RT/0502086}

\bibitem{DiJa92}
{\sc R.~Dipper and G.~D. James}, Representations of Hecke algebras of 
type $B_n$, J. Algebra {\bf 146} (1992), 454--481.

\bibitem{DiJa95}
{\sc R.~Dipper, G.~D. James and G.~E. Murphy}, {H}ecke algebras of type
${B}_n$ at roots of unity, Proc. London Math. Soc. (3) {\bf 70} (1995),
505--528.


\bibitem{my98}
{\sc M.~Geck}, Kazhdan-Lusztig cells and decomposition numbers, 
Represent.\ Theory {\bf 2} (1998), 264--277 (electronic).

\bibitem{myind}
{\sc M.~Geck}, On the induction of Kazhdan--Lusztig cells, Bull. London
Math. Soc. {\bf 35} (2003), 608--614.

\bibitem{my04}
{\sc M.~Geck}, Computing Kazhdan--Lusztig cells for unequal parameters,
 J. Algebra {\bf 281} (2004), 342--365; section "Computational Algebra".

\bibitem{GeIa05}
{\sc M.~Geck and L.~Iancu},  Lusztig's $a$-function in type $B_n$ in the
asymptotic case, preprint (april 2005); available at
{\tt http://arXiv.org/math.RT/0504213}.

\bibitem{ourbuch}
{\sc M.~Geck and G.~Pfeiffer}, {\em Characters of finite Coxeter
groups and Iwahori--Hecke algebras},  London Math. Soc. Monographs,
New Series {21} (Oxford University Press, 2000).

\bibitem{GrLe}
{\sc J.~J.~Graham and G.~I.~Lehrer}, Cellular algebras,
Invent.~Math. {\bf 123} (1996), 1--34.

\bibitem{KaLu}
{\sc D.~A. Kazhdan and G.~Lusztig}, {Representations of {C}oxeter groups and
{H}ecke algebras}, Invent. Math. {\bf 53} (1979), 165--184.
                                                                   
\bibitem{Lu0}
{\sc G.~Lusztig}, {On a theorem of {B}enson and {C}urtis}, 
J. Algebra {\bf 71} (1981), 490--498.

\bibitem{Lusztig83}
{\sc G.~Lusztig}, Left cells in {W}eyl groups, {\em Lie Group
Representations, I} (eds R.~L. R.~Herb and J.~Rosenberg), Lecture Notes
in Mathematics 1024 (Springer, Berlin, 1983), pp.~99--111.

\bibitem{LuBook}
{\sc G.~Lusztig}, {\em Characters of reductive groups over a finite field},
Annals Math.\ Studies 107 (Princeton University Press, 1984).
 
\bibitem{Lu1}
{\sc G.~Lusztig}, Cells in affine Weyl groups, Advanced Studies in Pure
Math. {\bf 6}, Algebraic groups and related topics, Kinokuniya and
North--Holland, 1985, 255--287.

\bibitem{Lusztig03}
{\sc G.~Lusztig}, Hecke algebras with unequal parameters, CRM Monographs
Ser.~{\bf 18}, Amer. Math. Soc., Providence, RI, 2003.

\bibitem{max}
{\sc M.~Neunh\"offer}, Untersuchungen zu James' Vermutung \"uber
Iwahori--Hecke--Algebren vom Typ $A$, Ph. D. Thesis, RWTH Aachen, 2003.

\bibitem{Spr}
{\sc T. A.~Springer}, {Quelques applications de la cohomologie d'intersection},
S\'{e}minaire Bourbaki (1981/82), exp. 589, Ast\'erisque {\bf 92--93} (1982).

\end{thebibliography}
\end{document}